\documentclass[12pt,letterpaper]{article}

\usepackage{epsfig}
\usepackage{amsmath}
\usepackage{amssymb}
\usepackage{amsthm}
\usepackage{indentfirst}
\usepackage{xspace}
\usepackage{multirow}
\usepackage{hyperref}
\usepackage{xcolor}
\definecolor{darkred}{rgb}{0.5,0.15,0.15}
\hypersetup{colorlinks=true,urlcolor=darkred,linkcolor=darkred,citecolor=darkred}
\usepackage{verbatim}
\usepackage[letterpaper,margin=1in,headheight=15pt]{geometry}
\usepackage{mathpazo}
\usepackage{authblk}
\usepackage{tikz-cd}

\numberwithin{equation}{section}

\newcommand{\cD}{\ensuremath{\mathcal D}}
\newcommand{\cL}{\ensuremath{\mathcal L}}
\newcommand{\cS}{\ensuremath{\mathcal S}}
\newcommand{\cF}{\ensuremath{\mathcal F}}

\newcommand{\cZ}{\ensuremath{\mathcal Z}}
\newcommand{\cM}{\ensuremath{\mathcal M}}
\newcommand{\cO}{\ensuremath{\mathcal O}}
\newcommand{\cP}{\ensuremath{\mathcal P}}

\newcommand{\cX}{\ensuremath{\mathcal X}}

\newcommand{\cW}{\ensuremath{\mathcal W}}
\newcommand{\cG}{\ensuremath{\mathcal G}}
\newcommand{\cT}{\ensuremath{\mathcal T}}

\newcommand{\R}{\ensuremath{\mathbb R}}
\newcommand{\C}{\ensuremath{\mathbb C}}

\newcommand{\Z}{\ensuremath{\mathbb Z}}
\newcommand{\Q}{\ensuremath{\mathbb Q}}

\newcommand{\N}{{\mathcal N}}

\newcommand{\HH}{{\mathcal H}}

\newcommand{\hk}{hyperk\"ahler\xspace}

\newcommand{\I}{{\mathrm i}}

\newcommand{\de}{\mathrm{d}}
\newcommand{\ab}{\mathrm{ab}}

\newcommand{\IP}[1]{\langle#1\rangle}

\newcommand{\eps}{\epsilon}
\newcommand{\fro}{{\overline{\underline{\Omega}}}}
\newcommand{\simarrow}{\overset{\sim}{\rightarrow}}

\newcommand{\ti}[1]{\textit{#1}}

\DeclareMathOperator{\re}{Re}
\DeclareMathOperator{\Tr}{Tr}
\DeclareMathOperator{\End}{End}
\DeclareMathOperator{\Hol}{Hol}

\newcommand{\insfig}[2]{\begin{figure}[htbp] \centering \includegraphics[scale=0.3]{figures/#1-crop.pdf} \caption{#2} \label{fig:#1} \end{figure}}

\newcommand{\insfigscaled}[3]{\begin{figure}[htbp] \centering \includegraphics[scale=#2]{figures/#1-crop.pdf} \caption{#3} \label{fig:#1} \end{figure}}

\begin{document}

\title{Spectral networks and Fenchel-Nielsen coordinates}
\date{}
\author[1,2]{Lotte Hollands}
\author[3]{Andrew Neitzke}
\affil[1]{Mathematical Institute, University of Oxford}
\affil[2]{Department of Mathematics, Heriot-Watt University}
\affil[3]{Department of Mathematics, University of Texas at Austin}

\maketitle

{\abstract{
It is known that spectral networks naturally 
induce certain coordinate systems on moduli spaces 
of flat $SL(K)$-connections on surfaces, 
previously studied by Fock 
and Goncharov. We give a self-contained account of this 
story in the case $K=2$,
and explain how it can be extended to incorporate
the complexified Fenchel-Nielsen coordinates.
As we review, the key ingredient in the story 
is a procedure for passing between moduli of flat
$SL(2)$-connections on $C$ (equipped with a little extra 
structure) and moduli of equivariant $GL(1)$-connections
over a covering $\Sigma \to C$; taking holonomies of the equivariant
$GL(1)$-connections then gives the desired 
coordinate systems on moduli of $SL(2)$-connections.
There are two special types of spectral network,
related to ideal triangulations and pants decompositions of $C$;
these two types of network
lead to Fock-Goncharov and complexified 
Fenchel-Nielsen coordinate systems 
respectively.\\
\vspace{0.5mm}
}}

\setcounter{page}{1}

\tableofcontents

\section{Introduction}

In the recent paper
\cite{Gaiotto2012} a process called \ti{nonabelianization} is described.  Nonabelianization
and its inverse, abelianization,
allow one to study flat connections in rank $K$ vector bundles over a (perhaps punctured) surface $C$,
by relating them to simpler objects, namely connections in rank $1$ bundles
over $K$-fold branched covers of $C$. 
 
The main new result of this paper is that abelianization gives
a new way of thinking about
a classical object, namely the complexified \ti{Fenchel-Nielsen coordinate system} on a 
certain covering of the moduli space
of $PSL(2)$-connections over $C$.\footnote{In this paper we always consider the \ti{complex} groups,
so $PSL(2)$ means $PSL(2,\C)$ and $SL(2)$ means $SL(2,\C)$.}

\medskip

Before describing the main points of this paper, let us briefly review what is meant by (non)abelianization.
Suppose given a flat $GL(1)$-connection $\nabla^\ab$ in a complex line bundle
over a branched cover $\Sigma$ of $C$.\footnote{More precisely, it will turn out that we should consider $GL(1)$-connections with an additional
equivariance property, as explained in \S\ref{ssec:propab} below, and as a consequence
of this equivariance they will not be flat but rather have holonomy $-1$
around the branch locus of the covering $\Sigma \to C$; 
for the rest of this introduction we suppress this subtlety.}
A naive way of relating $\nabla^\ab$ to a rank $K$ connection
over $C$ would be just to take the pushforward.  This indeed gives a connection $\nabla = \pi_* \nabla^\ab$ in a rank $K$
bundle over $C$, or more precisely over the complement of the branch locus in $C$.
However, the connection $\pi_* \nabla^\ab$ does \ti{not} extend smoothly over the branch locus.
The reason is simple:  when one goes around a branch point of the covering $\Sigma \to C$,
the sheets of $\Sigma$ are permuted, and this permutation shows up as a nontrivial monodromy
for $\pi_* \nabla^\ab$.  Thus $\pi_* \nabla^\ab$ cannot be extended to a flat
connection on the whole of $C$.

The nonabelianization operation is a recipe for ``fixing up'' this monodromy problem, by
cutting the surface $C$ into pieces, then regluing
with unipotent automorphisms of $\pi_* \nabla^\ab$ along the gluing lines.  The effect of this
cutting-and-gluing is to cancel the unwanted monodromy around branch points 
while not creating any monodromy anywhere else.

The cutting-and-gluing is done along a collection of paths on $C$,
carrying some extra labels and obeying some conditions, called a
\ti{spectral network} $\cW$.
Spectral networks were described in \cite{Gaiotto2012}, where it was shown that they are a fundamental ingredient in the story
of BPS states and wall-crossing for $\N=2$ supersymmetric quantum field theories of ``class $S$'' \cite{Gaiotto:2009we,Gaiotto:2009hg}.
For subsequent work on spectral networks in $\N=2$ theories see \cite{Gaiotto:2012db,Galakhov:2013oja,Hori:2013ewa,Maruyoshi:2013fwa}.
Essentially the same object has also been considered in the mathematical literature on Stokes phenomena, e.g. \cite{MR2132714}.

Fixing a spectral network $\cW$ and given a generic flat
$GL(1)$-connection $\nabla^\ab$, we can nonabelianize $\nabla^\ab$. We can also consider
an inverse operation, abelianization:  given a generic flat rank $K$ connection $\nabla$ over $C$, 
we attempt to realize it by cutting-and-gluing of $\pi_* \nabla^\ab$
along the spectral network $\cW$, for some flat $GL(1)$-connection $\nabla^\ab$. 

\medskip 

In this paper we formalize the (non)abelianization process,
focussing on $K=2$. We start with the definition of a spectral network in 
\S\ref{ssec:defspectralnetwork}.\footnote{For our purposes in this
  paper we do not need the  
most general notion of spectral network, in the case $K=2$ things are simpler.
}
In the remainder of
\S\ref{sec:spectralnetworks} and \S\ref{sec:WKB} we introduce various
examples;  in particular we introduce two classes of spectral networks, which we 
call \ti{Fock-Goncharov} (which had been considered 
earlier, in  \cite{Gaiotto:2012db}) and \ti{Fenchel-Nielsen}
networks (which are new), for a reason which will  become apparent momentarily.
Fock-Goncharov networks are dual to ideal triangulations, whereas Fenchel-Nielsen ones are related to pants
decompositions of $C$.

\S\ref{sec:abelianization} introduces the concept of $\cW$-abelianization.
We formulate this as an extension of a flat $SL(2)$-connection $\nabla$
into a so-called $\cW$-pair, which contains the
connection $\nabla$ in the ``abelianized gauge'' $\pi_*
\nabla^\ab$. We explicitly construct $\cW$-pairs for 
Fock-Goncharov and Fenchel-Nielsen networks in 
\S\ref{sec:constructingabelianizations}. 

For the particular networks we consider in this paper,
there is an almost unique way of $\cW$-abelianizing a given generic
connection $\nabla$.  ``Almost unique'' 
means that there are some discrete choices to be made, which in this paper we encapsulate in 
a notion of \ti{$\cW$-framed} connection; then, we find a canonical
1-1 correspondence between $\cW$-framed 
flat $SL(2)$-connections $\nabla$ over $C$ and $\cW$-abelianizations.

In the latter half of the paper we discuss $\cW$-nonabelianization. Similarly, we formulate this as
an extension of a flat $GL(1)$-connection $\nabla^\ab$ into a $\cW$-pair. For
convenience later on, we 
do this in the language of  parallel transport of the relevant
connections. We explain the dictionary in
\S\ref{sec:integratedversion}. In \S\ref{sec:nonabelianization} we
introduce the notion of 
$\cW$-nonabelianization and  show that there is a canonical 1-1
correspondence between flat $GL(1)$-connections $\nabla^\ab$ and
$\cW$-nonabelianizations when $\cW$ is either a
Fock-Goncharov or Fenchel-Nielsen type network. In 
\S\ref{sec:examples}  we discuss concrete examples.

We combine these two pieces in \S\ref{sec:modulispaces}. Given any
spectral network $\cW$ there is a diagram
\begin{center}
\includegraphics{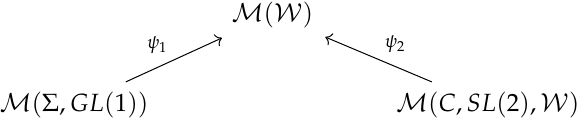}
\end{center}
where $\cM(\Sigma, GL(1))$ is the moduli space of flat
$GL(1)$-connections, $\cM(\cW)$ the moduli space of $\cW$-pairs and $
\cM(C, SL(2), \cW) $ the moduli space of $\cW$-framed
$SL(2)$-connections. In this diagram $\psi_1$ denotes
$\cW$-nonabelianization 
and $\psi_2$ denotes $\cW$-abelianization. 
Since we have shown that both maps $\psi_1$ and $\psi_2$ are bijections when $\cW$ is a
Fock-Goncharov or Fenchel-Nielsen type network, we obtain another
bijective map
\begin{equation}
\Psi_{\cW}: \cM(C, SL(2), \cW) \to \cM(\Sigma, GL(1)).
\end{equation}

\medskip

The relation to coordinate systems now comes about as follows.
Given a $\cW$-framed connection $\nabla$ and some $\gamma \in H_1(\Sigma,\Z)$,
let $\cX_\gamma(\nabla) \in \C^\times$ denote the holonomy around $\gamma$ of the corresponding $GL(1)$-connection 
$\nabla^\ab$.
For the networks $\cW$ we consider in this paper,
the collection of numbers $\cX_\gamma(\nabla)$ suffice to determine the connection $\nabla$.
In these cases, we thus get a coordinate system on the moduli space of $\cW$-framed flat connections $\nabla$; call this the \ti{spectral
coordinate system} associated to the network $\cW$.

 We discuss this in 
\S\ref{sec:modulispaces} and compute the spectral coordinates for both 
Fock-Goncharov and Fenchel-Nielsen type networks.

For the reader's convenience, we first
consider the class of Fock-Goncharov type network.  For these networks
it was shown in \cite{Gaiotto:2012db} 
(for all $K \ge 2$) that the spectral coordinates are the same as some
coordinates introduced earlier by Fock and Goncharov \cite{MR2233852}.
We give a self-contained re-derivation in the special case $K=2$
of the fact that Fock-Goncharov networks induce Fock-Goncharov coordinates.
We also clarify a few points which might be obscure in \cite{Gaiotto:2012db};
in particular we are more careful about how the story depends on whether
we consider connections for the group $SL(2)$ or $PSL(2)$.

We then move on to consider the new class
of Fenchel-Nielsen networks.
As we explain, the corresponding spectral
coordinates are a complexified version of Fenchel-Nielsen (length-twist)
coordinates.\footnote{Length-twist coordinates were originally defined
  by Fenchel and Nielsen, in a manuscript recently published
  as \cite{FN}. See also \cite{Wolpert1,Wolpert2}. A complex
  version of Fenchel-Nielsen coordinates was found in
  \cite{Kourouniotis, Tan,Nekrasov:2011bc}.}
This is one of our main points: 
\ti{spectral networks provide a framework
which naturally unifies Fock-Goncharov
and Fenchel-Nielsen coordinate systems.}

We finish this paper with some physical consequences in
\S\ref{sec:physcons}. As a preview, let us already make some remarks
here.

\subsection*{Spectral networks and quadratic differentials}

One suggestive way of understanding the results described above 
comes from the connection between spectral networks
and quadratic differentials.  Every meromorphic quadratic differential $\varphi_2$ on $C$ induces
a corresponding spectral network $\cW(\varphi_2)$, essentially the critical graph of $\varphi_2$.
When we consider $\varphi_2$ with at least one pole of order $\ge 2$,
and choose $\varphi_2$ generically within this class, the resulting
$\cW(\varphi_2)$ is a Fock-Goncharov network;
this is the basic reason why the previous work \cite{Gaiotto:2009hg} focused on
Fock-Goncharov networks and the associated Fock-Goncharov coordinate systems.

This motivates the question:  what about \ti{non}-generic $\varphi_2$?  If
$\varphi_2$ is a Strebel differential --- this case is maximally non-generic in 
a sense --- then the corresponding network
$\cW(\varphi_2)$ is a Fenchel-Nielsen network.
(If $\varphi_2$ is somewhere intermediate between being fully generic  
and being a Strebel differential,
then we obtain a more general kind of network, called \ti{mixed} in this paper;
we do not explore these much, although we expect they can be studied by
the methods of this paper.)

\subsection*{Motivations from physics}

Another motivation for being interested in Fenchel-Nielsen networks comes from the link between
spectral networks and $\N=2$ supersymmetric field theories.
Recall that given a punctured Riemann surface $C$
there is a corresponding $\N=2$ theory $S[A_1,C]$ as discussed in \cite{Witten:1997sc,Gaiotto:2009hg,Gaiotto:2009we}.
The quadratic differentials $\varphi_2$ just discussed
correspond to points of the Coulomb branch of this theory.

In some approaches to understanding the AGT correspondence \cite{Alday:2009aq},
complexified Fenchel-Nielsen coordinates on the moduli of $SL(2)$-connections play
an important role.
For example, in \cite{Nekrasov:2011bc}
a semiclassical limit of the Nekrasov partition function of the theory $S[A_1, C]$
gets identified with the generating function, in complexified Fenchel-Nielsen coordinates,
of the variety of opers inside the moduli space of flat $SL(2)$-connections on $C$.
Fenchel-Nielsen coordinates are also important in the interpretation of the AGT
correspondence given in \cite{Teschner2011}.

The result of this paper gives a new way of understanding the role of Fenchel-Nielsen coordinates
in the theories $S[A_1,C]$, and thus may shed further light on these approaches.
We limit ourselves here to a few optimistic comments, as follows.

\begin{itemize}

\item First, in this paper the complexified Fenchel-Nielsen coordinates get linked to
a particular locus of the Coulomb branch, corresponding to the set of Strebel differentials.
This locus is a
``real'' subspace inside the Coulomb branch:  indeed, it is the locus where all
of the vector multiplet scalars $a^I$ are real, with respect to a particular choice of
electromagnetic frame corresponding to the pants decomposition.
It is interesting that this real subspace has played a
distinguished role in applications to $\N=2$ theory in the past.
For example, in
the computation of the $S^4$ partition function via localization \cite{Pestun:2007rz},
one naturally considers integrals
not over the whole Coulomb branch but rather over this real subspace.
Also, it appears that the spectrum
of framed BPS states at this real subspace is especially simple:  as we describe in \S\ref{sec:defects} below,
at this locus the UV line defects corresponding to Wilson loops carry just two framed BPS
states, exactly matching the naive classical expectation.

 \item One would like to extend the approaches of \cite{Nekrasov:2011bc,Teschner2011}
to other Lie algebras, e.g. to the case $A_{K-1}$.  The question then arises:  \ti{what is the higher rank analogue of Fenchel-Nielsen coordinates?}

In this paper we find that Fenchel-Nielsen coordinates correspond to the spectral networks which appear at a
real locus of the Coulomb branch for $K=2$.  This suggests that one might be able to find higher rank Fenchel-Nielsen coordinates in a similar way,
by studying spectral networks appearing at a real locus of the Coulomb branch for $K>2$.

\item The constructions of \cite{Teschner2011} make use of a correspondence between
Fenchel-Nielsen length coordinates $\ell^I$ and local coordinate functions $a^I$ on the Coulomb branch
of the $\N=2$ theory, i.e. on the Hitchin base.
The invariant meaning of this correspondence is a bit mysterious at first glance.
After all, $\ell^I$ and $a^I$ are not the same.  Both of them can be viewed as functions on
the moduli space of flat connections, but they are different functions, indeed 
even holomorphic in different complex
structures on that \hk space --- the $\ell^I$ are holomorphic in the usual complex
structure, while the $a^I$ are holomorphic after rotating
to the Higgs bundle complex structure.
The correct relation between the two is an \ti{asymptotic} one, of the schematic form 
\begin{equation}
	\ell^I \sim const \times \exp(const \times a^I / \zeta) \quad \text{ as } \quad \zeta \to 0,
\end{equation}
obtained by applying the WKB approximation to certain 1-parameter
families of flat connections.  Of course, this relation is noted in \cite{Teschner2011} too.
What we add here is just the observation that
this asymptotic relation --- and an understanding of its precise range of validity ---
arises naturally from the point of view
of the Fenchel-Nielsen spectral networks; see \S\ref{sec:asymptotics} below.

\end{itemize}

\subsection*{Acknowledgements}

We thank Anna Wienhard for a very helpful discussion, and the anonymous referee
for many clarifying suggestions, comments and corrections. The research of
LH is supported by a Royal Society Dorothy Hodgkin fellowship.
The research of AN is supported by NSF grant 1151693. This work was
also supported in part by the National Science Foundation under Grant
No. PHYS-1066293 and we thank the hospitality of the Aspen Center for Physics.

\section{Spectral networks}\label{sec:spectralnetworks}

In this section we start with the definition of spectral network, and
then introduce the particular classes of spectral networks we will
consider in this paper.

We first describe the class of ``Fock-Goncharov networks.''
Such networks are associated to ideal triangulations of the
surface $C$.  They played the starring role in \cite{Gaiotto:2009hg}.
Second, we introduce a new class of spectral networks which we call ``Fenchel-Nielsen
networks.''  These networks are associated to pants decompositions of the
surface $C$ rather than triangulations.  They have not been considered explicitly before.
Third, we briefly consider two more general kinds of network:  one obtained by
collapsing some double walls in a Fenchel-Nielsen network, which we
call ``contracted Fenchel-Nielsen'',
and another which is part
Fock-Goncharov and part Fenchel-Nielsen, which we call ``mixed''.

\subsection{Spectral networks}\label{ssec:defspectralnetwork}

Let $C$ denote an oriented real surface of negative Euler characteristic.
$C$ may (but need not) have punctures and/or boundary components.

Let $\Sigma$ denote a branched double cover of $C$, 
which restricts to a trivial cover in a neighborhood of each puncture or
boundary component.
We often find it convenient to choose a system of branch cuts on $C$ and trivialize
the cover $\Sigma$ on the complement of the cuts, labeling its sheets by an index $i \in \{1,2\}$.

A \ti{spectral network $\cW$ subordinate to $\Sigma$}
is a collection of oriented open paths $w$ on $C$, which we call
\ti{walls}, with the following properties and structure:
\begin{itemize}
\item Each wall $w$ begins at a branch point of $\pi: \Sigma \to C$, and
ends either at a puncture or at a branch point. A wall ends
at a branch point if and only if it overlaps with another wall which is oriented
in the opposite direction.  We say that the locus where two walls
overlap is a ``double wall.''

\item Walls never intersect one another or self-intersect, except that it is possible 
for two oppositely oriented walls to overlap, as we have just remarked.

\item Each wall $w$ carries an ordering of the 2 sheets of the covering $\Sigma$ over $w$.
If we fix a system of branch cuts and  
trivialization of $\Sigma$, this amounts to saying $w$ locally carries either the label $12$
or $21$, and the label changes when $w$ crosses a branch cut.

\insfigscaled{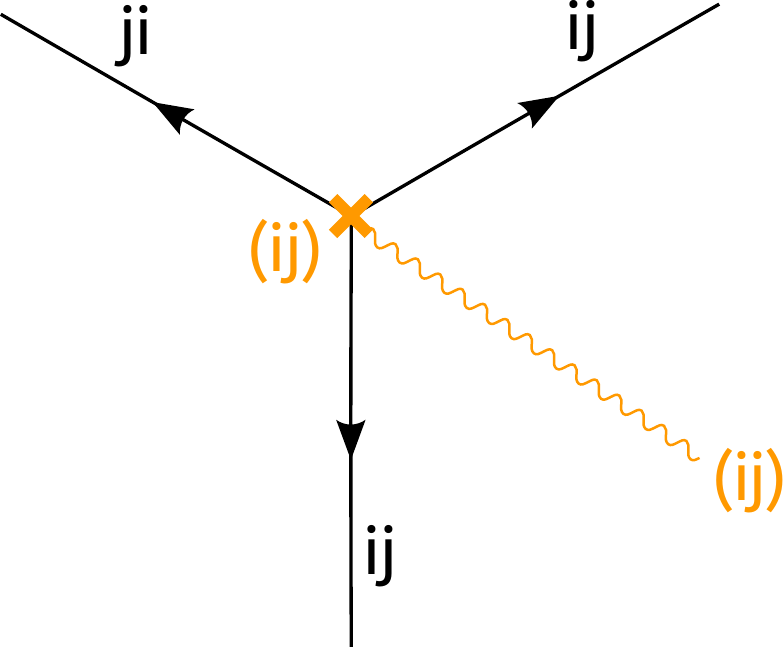}{0.3}{Local configuration of
walls near a branch point.}

\item Exactly three walls $w$ begin at each branch point.  If we fix
 a system of branch cuts and trivialization of $\Sigma$, then
consecutive walls carry opposite sheet orderings, unless we cross
a branch cut in moving from one to the next (or more generally, an odd number of branch cuts); see Figure
\ref{fig:branch-point-trajectories}.

\item The two constituents in a double wall carry opposite sheet orderings:  thus,
if we have fixed a labeling of the sheets,
then locally one is labeled $12$ and the other $21$.

\item If there is at least one double wall, then the network comes with an additional discrete datum, the choice of a \ti{resolution}, which is
one of the two words ``British'' or ``American''.  We think of the resolution as telling us how
the two constituents of each double wall are infinitesimally displaced from one another; see Figure \ref{fig:two-resolutions} for the rules.

\insfigscaled{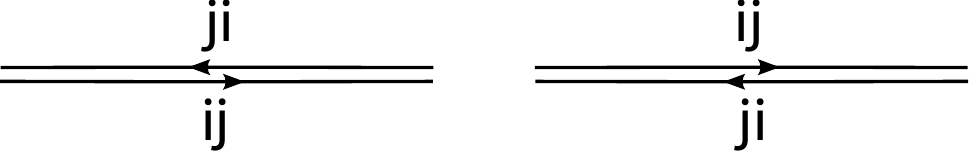}{0.4}{Left: a double wall with American resolution. Right: a double wall with British resolution.}

\item The closure of
each connected component of $C \setminus \cW$ must either
have one of the topologies shown in Figure~\ref{fig:topology-components}, or be obtained
from one of those topologies by some identifications along the
boundary. 
(In particular, this implies that each puncture has at least one wall ending on
it, while boundary components have no walls ending on them.)

\insfigscaled{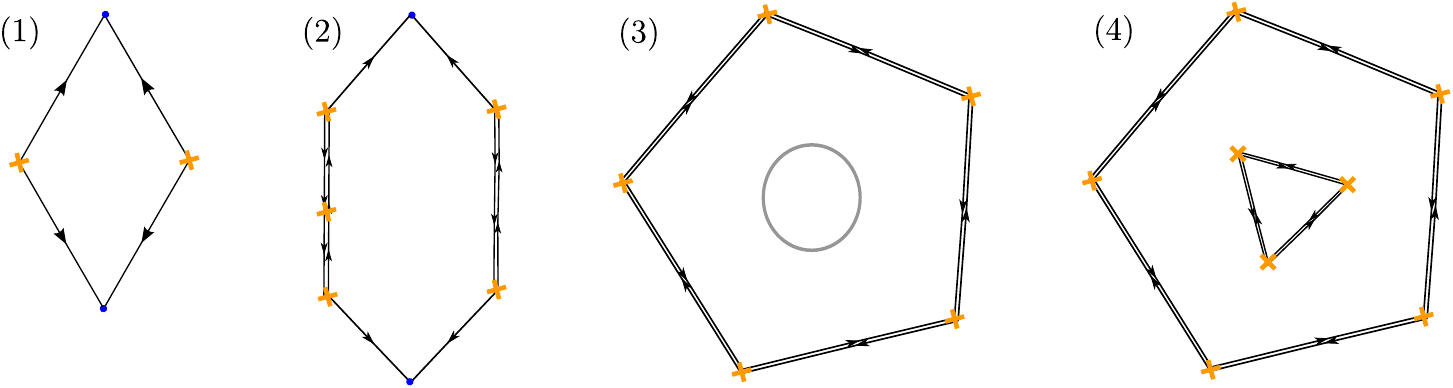}{0.65}{Possible topologies
  for the connected components of $C \setminus \cW$: (1) A generic
  cell, with two punctures and two branch points as vertices. (2)
  A degenerate cell, with two punctures and at least three 
  branch points as vertices; each pair of neighboring branch
  points is connected by a double wall. (3) An annulus bounded
  by a boundary component on one side and a polygon of double walls on the
  other. (4) An interior annulus, bounded by a polygon of double walls on each side,
  with the same resolution on each double wall.}

\item The spectral network carries a \ti{decoration}, which 
means a collection of decorations, one for each puncture on $C$ and
one for each annulus of $C \setminus \cW$:
\begin{itemize}
\item  A decoration of a puncture $z_l \in C$ is
an ordering of the sheets of $\Sigma$ over $z_l$, such that the
labels of all walls which end on $z_l$ match the decoration at $z_l$.
\item 
Call a cell of type (3) or (4) from Figure \ref{fig:topology-components} an 
\ti{annulus in $C \setminus \cW$}.
Each such annulus $A$ retracts to a circle $c(A) \subset A$.
A decoration of $A$ is an assignment of an ordering of the sheets of $\Sigma$ over $A$ 
to each orientation $o$ of $c(A)$, such that reversing $o$ also reverses the ordering.
The orientation of any wall $w$ on the 
boundary of $A$ induces an orientation $o_w$ on $c(A)$
(to see this, isotope $c(A)$ to the boundary of $A$, so that $w$ becomes
a subset of $c(A)$); the ordering which the decoration assigns to the orientation $o_w$
must match the label of $w$.
\end{itemize}

\end{itemize}

Let us emphasize that while the punctures and the covering $\Sigma \to C$
are part of the definition of a spectral network, the choice of
branch cuts and trivialization of $\Sigma$ are not:  they just give
an explicit presentation of $\Sigma$, which comes
in handy when we want to draw figures and describe examples explicitly.

For later use we remark that any spectral network as defined above
has $-2 \chi(C)$ branch points.
Indeed, we can compute $\chi(C)$ by summing up contributions from
the cells. We compute the contribution from a single cell
as follows: each edge contributes $-1/2$ (since it will appear 
on the boundary of two cells), each branch point contributes
$+1/3$ (since it will appear on the boundary of three cells),
and the interior face contributes $+1$.
For example, for a cell of type (2) in Figure \ref{fig:topology-components},
we have $b$ branch points and $b+2$ edges, so we
get in this way $b/3 - (b-2)/2 + 1 = -b/6$.
Similarly for the other cell topologies we always find $-b/6$.
Thus summing up over all cells 
we get $\chi(C) = - \sum b_i / 6$. Since $\sum b_i$ is 3 times
the total number of branch points on $C$ this gives
the desired result.

Finally we define the notion of \ti{equivalence}.
Roughly speaking, two spectral networks are equivalent if one
can be isotoped into the other.
More precisely, given two spectral
networks $\cW$ and $\cW'$ subordinate to covers $\Sigma$ and $\Sigma'$,
an \ti{equivalence} between $\cW$ and $\cW'$ 
is a 1-parameter family of spectral networks $\cW_t$
and covers $\Sigma_t$, together with an identification of 
$(\cW_0,\Sigma_0)$ with $(\cW,\Sigma)$
and $(\cW_1,\Sigma_1)$ with $(\cW',\Sigma')$.
Most of our constructions depend only on a spectral network
up to equivalence.

\subsection{Fock-Goncharov networks}

Suppose that $C$ has at least one puncture and
no boundary components.
Fix an
\emph{ideal triangulation}~$\cT$~of~$C$, i.e. a triangulation for
which the set of vertices is
precisely the set of punctures on $C$.\footnote{We could allow $\cT$ to include some ``degenerate'' triangles,
obtained e.g. by gluing together two vertices or two edges of an ordinary triangle; see \cite{MR2233852,Gaiotto:2009hg} for the
precise class of triangulations one can allow.}
We now construct a network $\cW(\cT)$, the \ti{Fock-Goncharov network} corresponding to $\cT$:
\begin{enumerate}
\item First, choose one point $b$ inside each triangle of $\cT$.
We will construct a double cover $\pi: \Sigma \to C$ which is branched exactly at these points.
$\Sigma$ is given explicitly by gluing together two copies of $C$ along a system of branch cuts,
with $3$ cuts emanating from each branch point $b$ and exiting the $3$ sides of the triangle.
\item Next we describe the walls of the network $\cW(\cT)$:
for each branch point $b$, there are three walls, running
from $b$ to the three vertices of the triangle in which $b$ sits.
\item Finally, according to the definition of spectral network, 
we must provide each wall $w$ with an ordering of the two sheets of
$\Sigma$ over $w$, and we must also decorate each puncture with an ordering of the two sheets
near the puncture.  Since our construction of $\Sigma$
equips it with a natural trivialization away from the branch cuts, 
this amounts concretely to putting labels $ij$ on the walls and punctures.
We simply assign the label $12$ to every wall and to every puncture.\footnote{Choosing instead 21 everywhere would lead to an equivalent network $\cW'(T)$. (Indeed, we define the equivalence between $\cW(T)$ and $\cW'(T)$ by the identification of $\cW(T)$ with $\cW'(T)$ that changes all labels.)}  
\end{enumerate}

The conditions in 
\S\ref{ssec:defspectralnetwork} are then easy to verify. See Figure \ref{fig:triangle-cuts} for the local picture in each triangle.
\insfigscaled{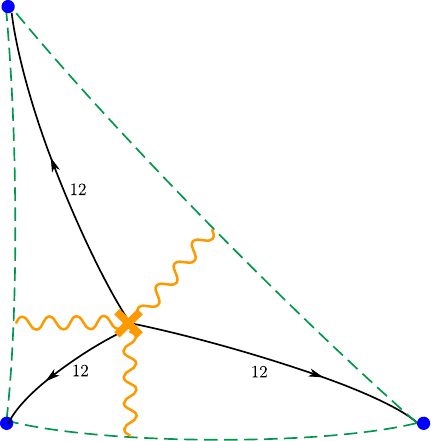}{0.55}{A Fock-Goncharov network $\cW(\cT)$, restricted to a single triangle of $\cT$.  Walls are
drawn as oriented black paths; the orange cross is
a branch point; wavy orange lines are branch cuts; dashed green paths are edges 
of $\cT$.
}

As just described, we build the covering $\Sigma$ by gluing along a 
complicated system of branch cuts; this presentation has the advantage of being canonical,
but in any particular example, $\Sigma$ can often be presented by gluing in a simpler way.
We will generally take advantage of that freedom when we show explicit examples.
Two complete networks are shown in
Figure~\ref{fig:FG-spectral-network-SU2Nf4}.
\insfigscaled{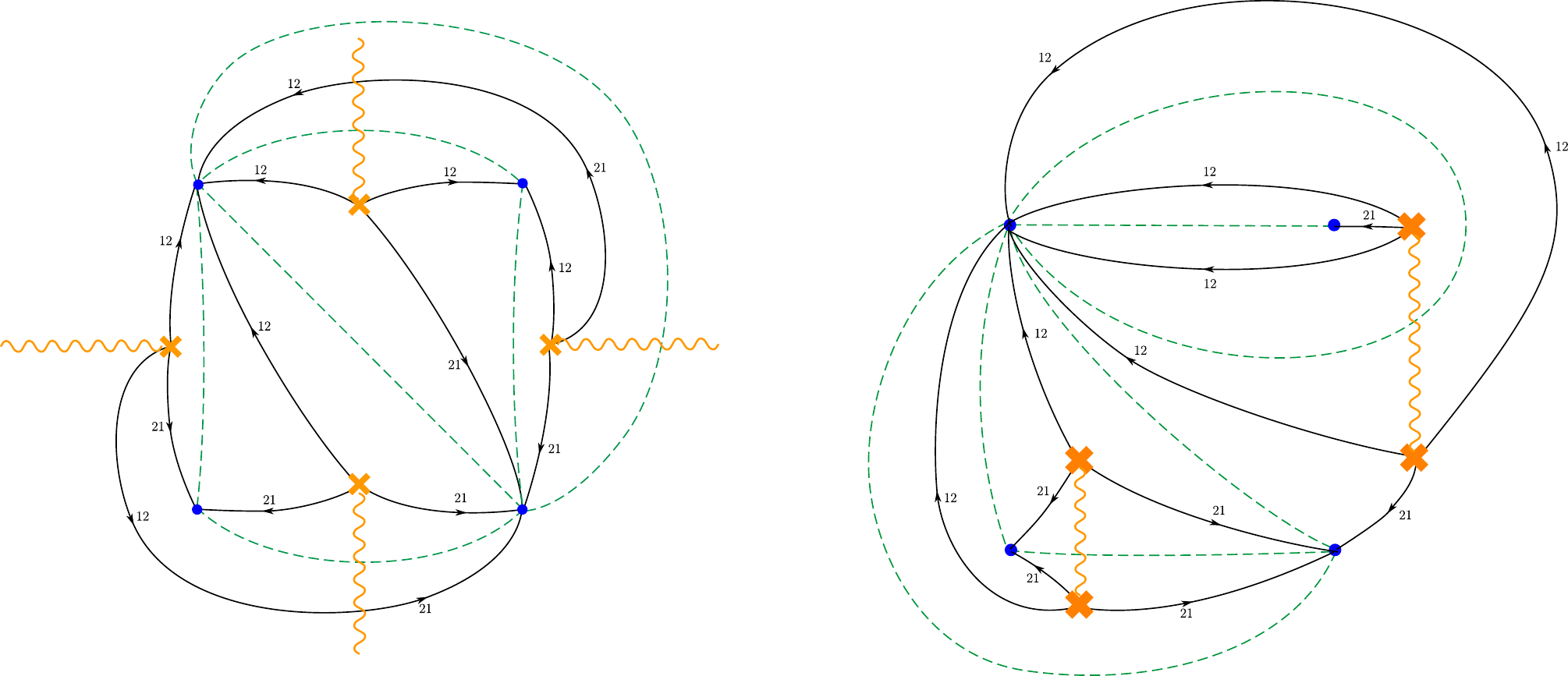}{0.45}{Two examples of
Fock-Goncharov networks on the four-punctured sphere (shown here as the plane, with the point 
at infinity omitted.)  The triangulation on the right includes two degenerate 
triangles.  To simplify the figures, here 
we chose to present the coverings with the minimal possible number of branch cuts;
nevertheless one can check directly that the networks shown here are equivalent
to the ones we described in the construction of the Fock-Goncharov network.
}

\subsection{Fenchel-Nielsen networks}

Now consider a surface $C$ with no punctures, and possibly with
boundary components.  Fix a pants decomposition $\cP$ of $C$.
Given this decomposition we now explain how to construct several different networks,
all of which we call \ti{Fenchel-Nielsen networks} associated to $\cP$.
On each pair of pants 
we fix the branched cover $\Sigma$ and spectral network $\cW$ to be one of the 
two shown in Figure \ref{fig:molecules}.
We then glue together along the boundaries of the pairs of pants,
inserting a circular branch cut around each pants curve.
Finally we must specify the decoration:  on each annulus in Figure \ref{fig:molecules}
it assigns the sheet ordering $12$ to the 
counterclockwise orientation.\footnote{Note that this would not have worked 
without the extra branch cuts around the pants curves.}
This completes the construction of the network.

  For instance,
Figure \ref{fig:four-holed-sphere} depicts a Fenchel-Nielsen network on the four-holed sphere (equivalent to the one obtained from the rules of gluing pants as just formulated). 

\insfigscaled{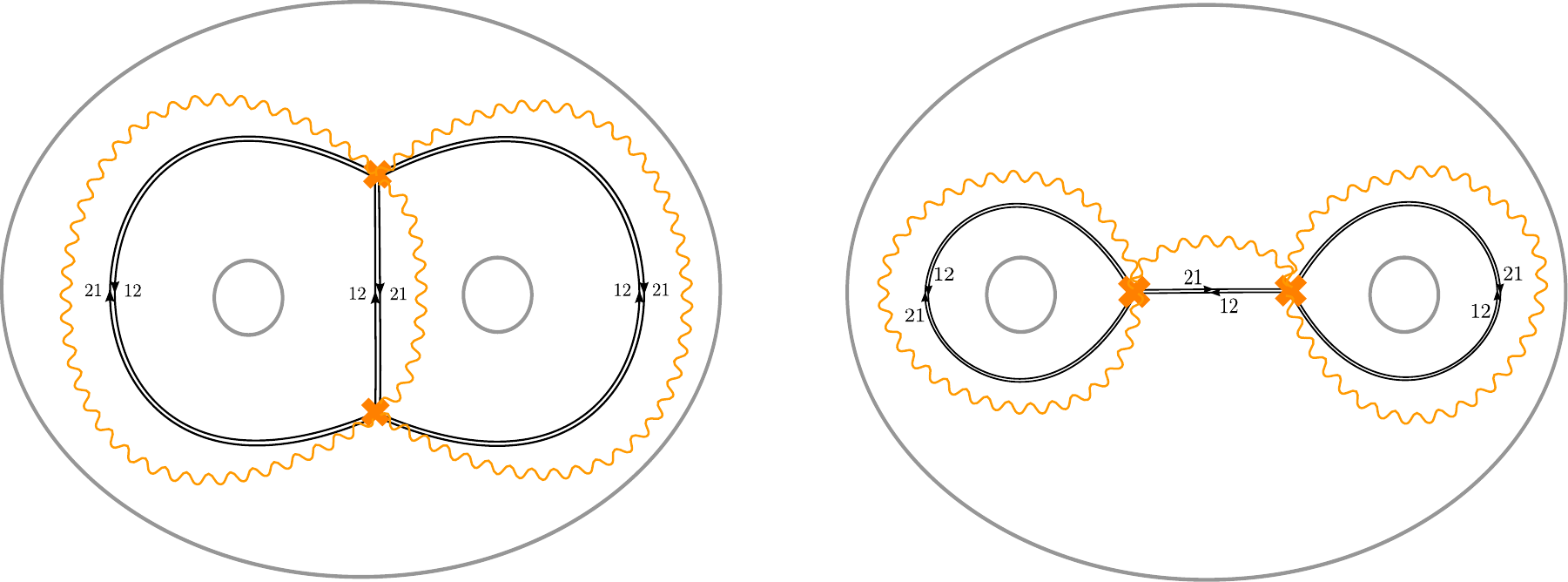}{0.55}{Two examples of
Fenchel-Nielsen spectral networks on pairs of pants.
The network on the left is ``molecule I'' and on the right is ``molecule II.''
Each network carries the ``British resolution''; to get the ``American
resolution'' one would reverse the orientation of every wall in the network.
}

\insfigscaled{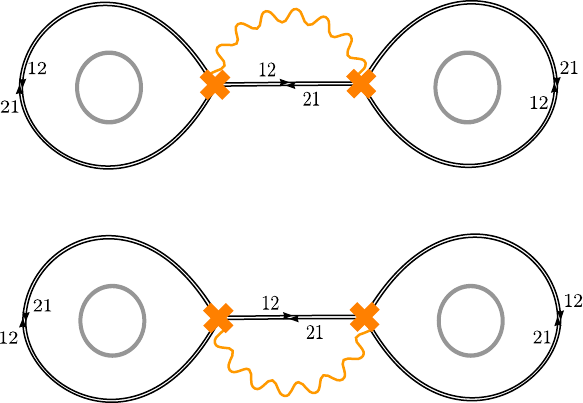}{0.8}{A Fenchel-Nielsen
network on the four-holed sphere.  This network is built
by gluing together two pairs of pants, each containing a copy of
molecule II from Figure \ref{fig:molecules}.}

There is a sense in which a
Fenchel-Nielsen network can be considered as a limit of Fock-Goncharov
networks.  See Figure~\ref{fig:molecule-limits} for an illustration of what we mean.  
This limiting procedure has a particularly transparent meaning if we consider the spectral
networks as coming from quadratic differentials on $C$. We explain
this in \S\ref{sec:WKB} below.

\insfigscaled{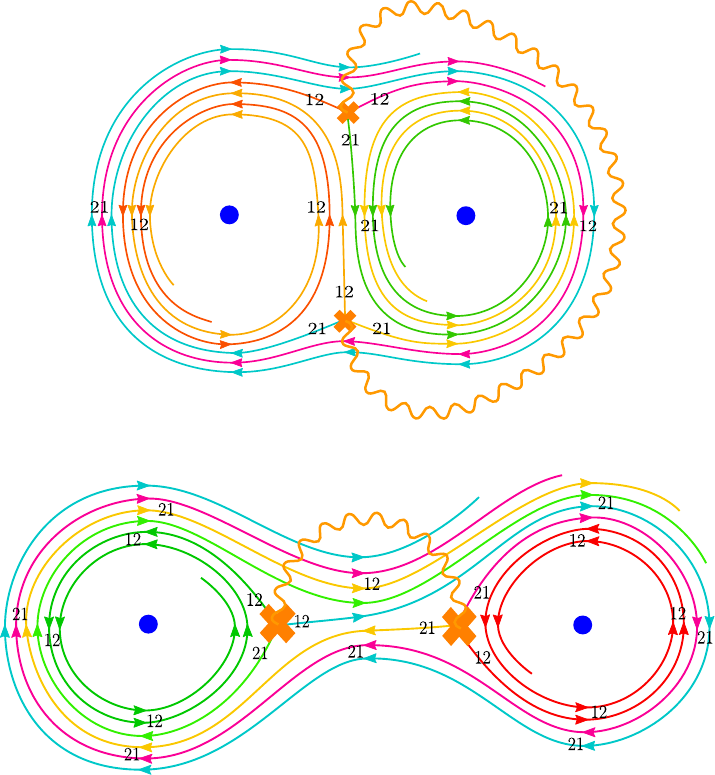}{0.8}{Two spectral
networks on the three-punctured sphere (with one puncture at infinity).
Each wall spirals many times around and asymptotically approaches some puncture; in the top network,
two walls approach each puncture, while in the bottom network, four walls approach the
puncture at infinity, and one wall approaches each of the other two punctures. 
To reduce clutter we show only truncated versions of the walls, and to help distinguish
the different walls, they are shown in different colors.
These networks are Fock-Goncharov networks, but
in the limit where the separation between walls goes to zero, they approach the
British resolution of Fenchel-Nielsen molecule~I (top) and molecule~II (bottom).}

\subsection{Contracted Fenchel-Nielsen networks}

As we have explained, in Fenchel-Nielsen networks all
walls are double walls; however, there are also spectral networks which are not
Fenchel-Nielsen networks in which all walls are double walls.   See e.g.
Figure~\ref{fig:deg-FN-network}. We refer to
such a network as a \emph{contracted} Fenchel-Nielsen network; 
the reason for the name is that any contracted Fenchel-Nielsen network
can be obtained by a degeneration of an ordinary Fenchel-Nielsen network where some annulus
shrinks to zero size,
in the sense indicated in Figure~\ref{fig:deg-FN-network-res}.\footnote{To see that \ti{every} contracted Fenchel-Nielsen network can be obtained in this way,
consider a connected component of the network containing more than two branch points, then look for 
two which are connected by exactly one double wall; one can define an operation which detaches these
two branch points from the rest, creating a new annulus; iterating this operation one eventually
reaches an ordinary Fenchel-Nielsen network. \label{footnote:FNcontracted}}
\insfigscaled{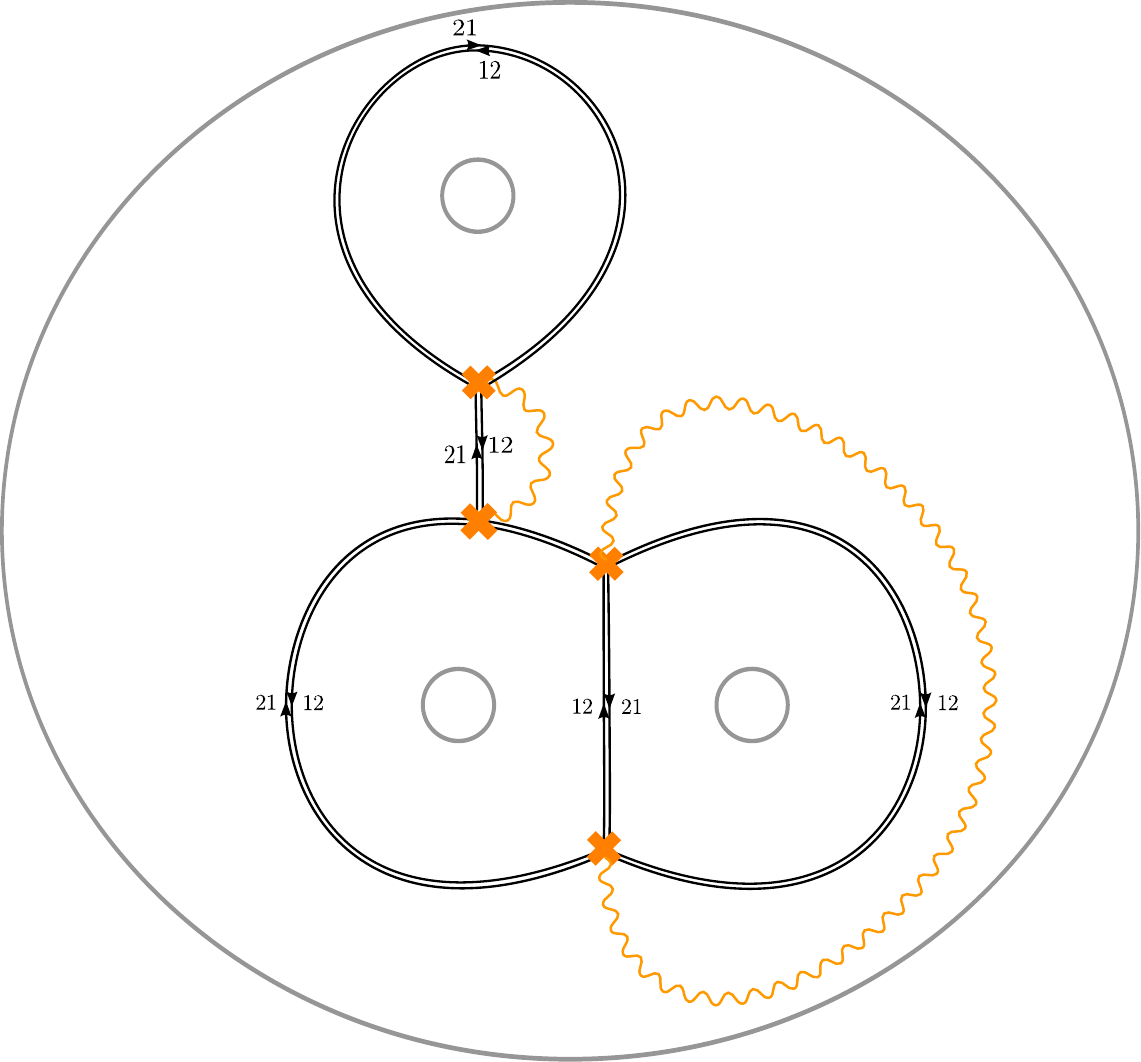}{0.4}{A contracted
Fenchel-Nielsen network on the four-holed sphere.}
In fact, as we see in that figure, one may obtain a single contracted Fenchel-Nielsen network
by degeneration from several distinct Fenchel-Nielsen networks, associated to different
pants decompositions.

\insfigscaled{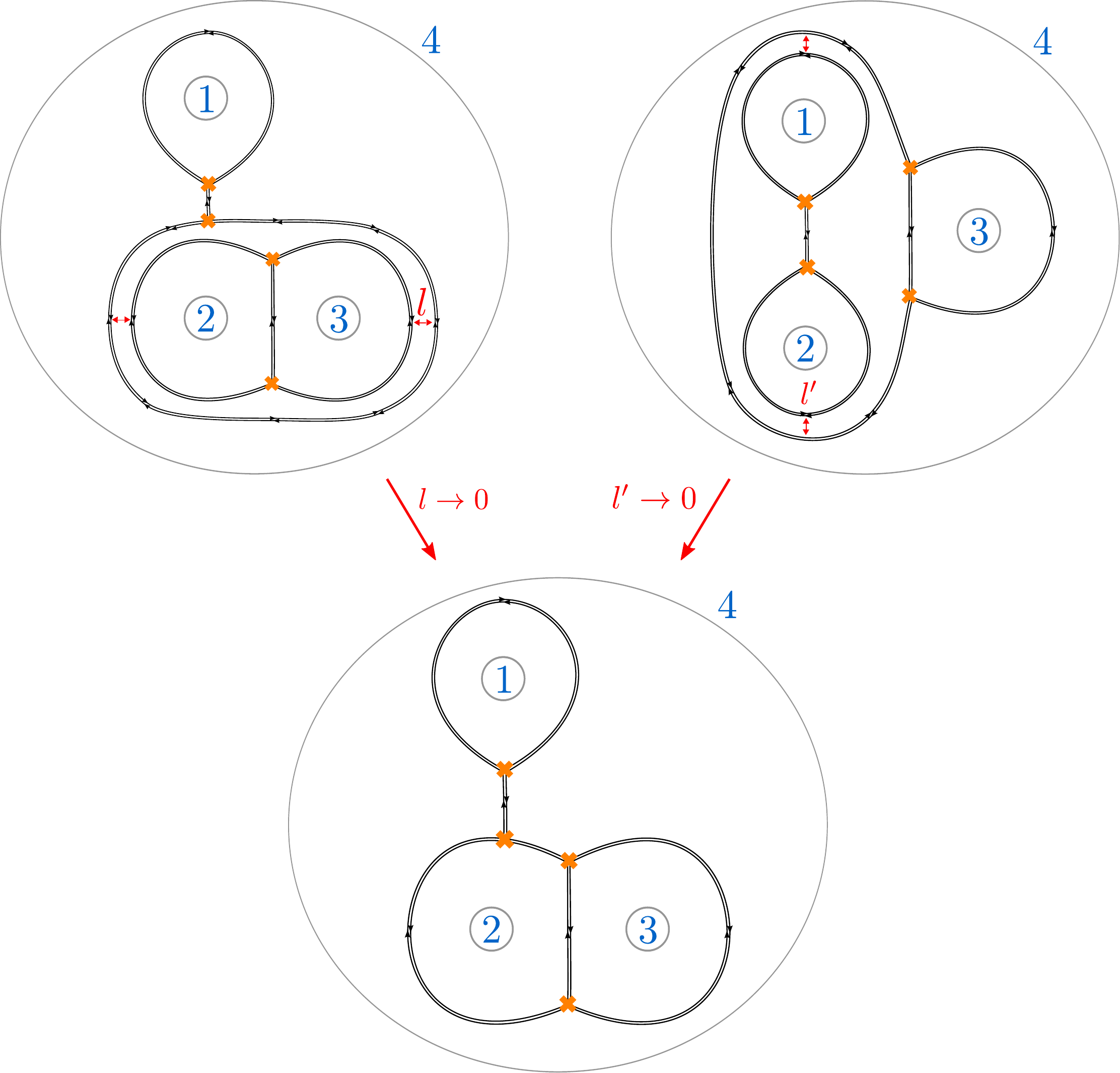}{0.3}{A contracted
  Fenchel-Nielsen network on the four-holed sphere, which can be
  obtained in a degeneration limit from two distinct Fenchel-Nielsen
  networks. These Fenchel-Nielsen networks are associated to two
  different choices of pants decomposition of the four-holed sphere. (The two branch points referred to in footnote~\ref{footnote:FNcontracted} are the ones connecting regions 1 and 2 in the $l\to0$ limit, and the ones connecting regions 2 and 3 in the $l'\to0$ limit.)  }

\subsection{Mixed spectral networks}\label{ssec:mixedspectralnetwork}

We call any spectral network that contains both single and double walls
a mixed spectral network.  An example is given in
Figure~\ref{fig:mixedsp1}.

\insfigscaled{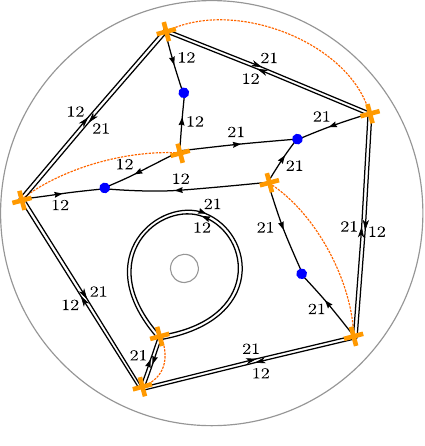}{0.95}{A mixed spectral network
on a sphere with four punctures and two holes.}

\section{WKB spectral networks}\label{sec:WKB}

Let $\bar C$ be a closed surface, $\{z_1, \dots, z_n\} \in \bar C$ some points,
and $C = \bar C \setminus \{z_1, \dots, z_n\}$.
In this section we review from \cite{Gaiotto:2009hg} how to construct a canonical spectral network
$\cW(\varphi_2, \vartheta)$ on $C$,
determined by the data of a complex structure on $\bar C$, 
a meromorphic quadratic differential $\varphi_2$ on $\bar C$ with double poles 
at the $z_l$, and a phase 
$\vartheta \in \R / 2 \pi \Z$.  Spectral networks
constructed in this way are known as \ti{WKB spectral networks}.

The class of WKB spectral networks overlaps with all classes
of spectral networks introduced above.
As we will describe, if $(\varphi_2, \vartheta)$ is chosen generically 
and $\varphi_2$ has at least one double pole, then $\cW(\varphi_2, \vartheta)$ is a
Fock-Goncharov network.  In contrast, if we choose 
$(\varphi_2, \vartheta)$ in a special way, so that 
$e^{- 2 \I \vartheta} \varphi_2$ is a \ti{Strebel differential} (i.e.
there is a pants decomposition of $C$ such that
the period of $\sqrt{\varphi_2}$ around each leg has phase
$\vartheta$),
then $\cW(\varphi_2, \vartheta)$ is a Fenchel-Nielsen network. (To be consistent with our previous definitions, we want to replace each puncture that is surrounded by a polygon consisting of double walls with a boundary.)

The reader who is only interested in understanding how
abelianization and nonabelianization work can safely skip this section.
However, the construction of WKB spectral networks
might help to motivate the
definition of spectral network, and will be relevant again in
the applications in \S\ref{sec:physcons}.

In the application to quantum field theories of class $S[A_1]$, the space of meromorphic 
quadratic differentials $\varphi_2$ on $C$ is the
Coulomb branch. Fock-Goncharov networks thus appear at generic
points of the Coulomb branch, and Fenchel-Nielsen networks
appear at distinguished ``real'' points.

\subsection{Foliations}

Let $\varphi_2$ be a
meromorphic quadratic differential on $\bar{C}$, holomorphic
away from the $z_l$. Locally such a differential is of the form
\begin{equation}
\varphi_2 = u(z) (\de z)^2.
\end{equation}
We suppose that $\varphi_2$ has only simple zeroes, and
second-order poles at the punctures $z_l$, with residues $-m_l^2 \neq 0$. Such a $\varphi_2$ has $-2 \chi(C)$ zeroes.

Fixing an angle $\vartheta \in \R / 2 \pi \Z$, define a \ti{$(\varphi_2,\vartheta)$-trajectory}
to be a real curve on $C$ such that, if $v$ denotes a nonzero tangent vector to the curve,
\begin{equation}
e^{-2 \I \vartheta} \varphi_2 (v^2) \in \R_+.
\end{equation}
Locally parametrizing such a curve as $z = \gamma(t)$, we can write this condition
more concretely as
\begin{equation}\label{eq:leaftheta}
e^{-2 \I \vartheta}\,  u(\gamma(t)) \left( \frac{\de \gamma(t)}{\de t}
\right)^2 \in \R_+.
\end{equation}
The $(\varphi_2, \vartheta)$-trajectories are leaves of a
singular foliation $\cF(\varphi_2, \vartheta)$ on $C$.
When $\vartheta = 0$ this is sometimes called the ``horizontal foliation'' attached
to $\varphi_2$.

The foliation $\cF(\varphi_2, \vartheta)$ has singularities at
the zeroes and poles of $\varphi_2$. Around a pole $z_l$ the foliation
depends on the value of $e^{-\I \vartheta} m_l$, as illustrated in
Figure~\ref{fig:foliation-pole}.
\insfigscaled{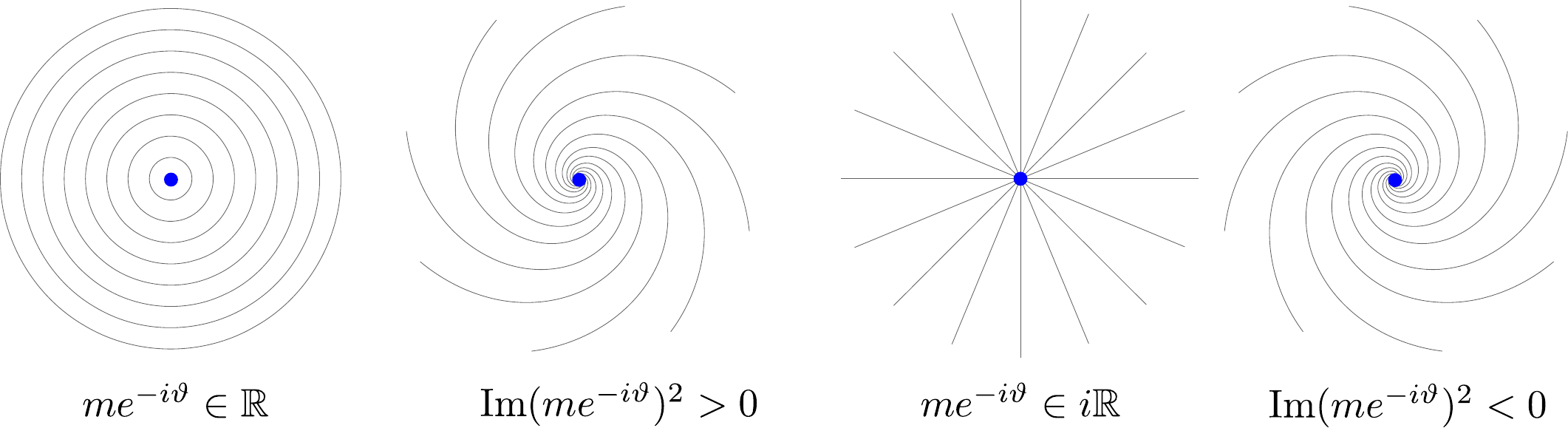}{0.4}{The behavior of the foliation $\cF(\varphi_2, \vartheta)$
around a pole of $\varphi_2$.}
In a neighborhood of a zero of $\varphi_2$, the foliation $\cF(\varphi_2, \vartheta)$ looks
like Figure~\ref{fig:foliation-zero}.
\insfigscaled{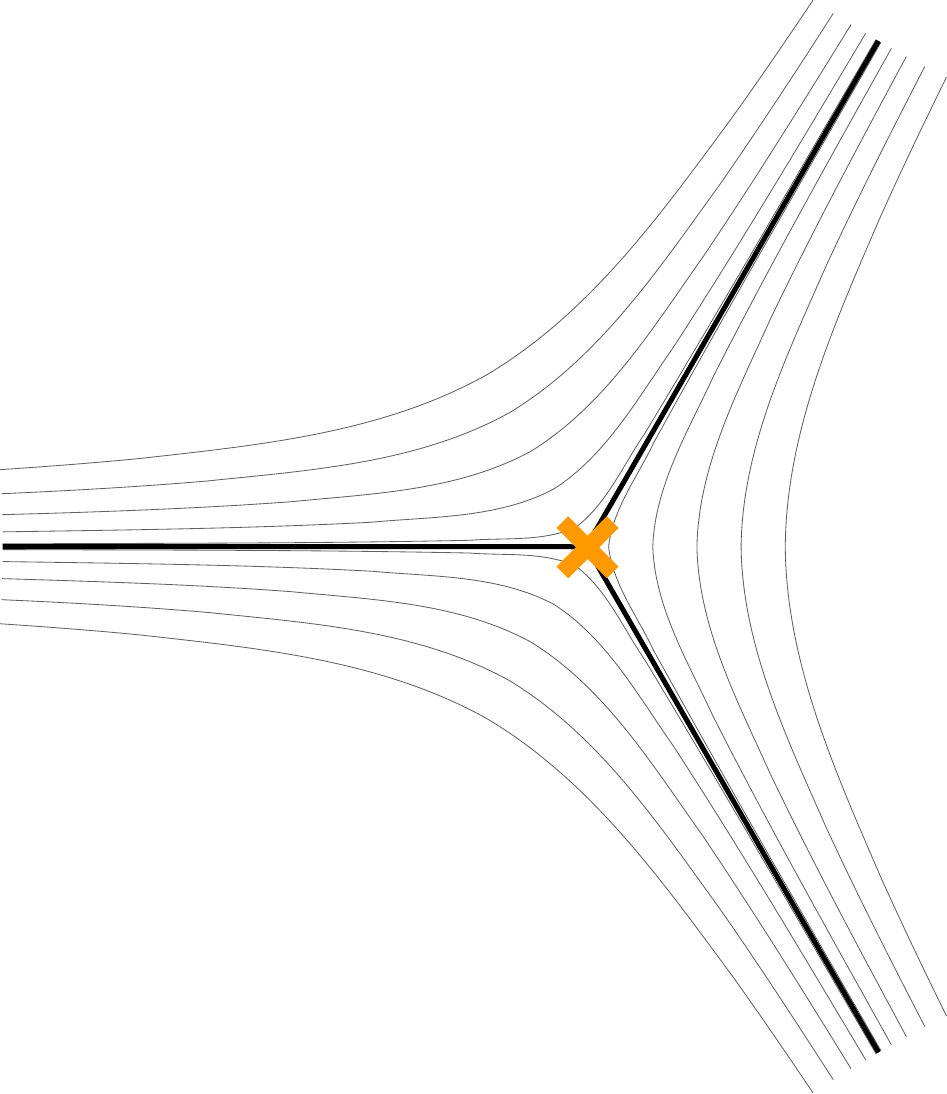}{0.3}{The behavior of the foliation $\cF(\varphi_2, \vartheta)$
around a zero of $\varphi_2$.  Three critical leaves emanate from the zero.}
The leaves which terminate at one end on a zero of $\varphi_2$ are called ``critical.''
The union of the critical leaves is known as the \ti{critical graph}
$CG(\varphi_2,\vartheta)$ of
the foliation.

Now we can define a spectral network $\cW(\varphi_2, \vartheta)$, as follows.
The double covering $\pi: \Sigma \to C$ will be the one given by the square roots of $\varphi_2$,
i.e.
\begin{equation}
 \Sigma = \{ (z \in C, \lambda \in T^*_z C): \lambda^2 = \varphi_2(z) \},
\end{equation}
with $\pi(z,\lambda) = z$.
In local coordinates, writing $\lambda = y\,\de z$, $\Sigma$ is given by the equation
\begin{equation}
y^2 = u(z). 
\end{equation}
Thus the branch points of the covering $\pi: \Sigma \to C$ are the zeroes of $\varphi_2$.

The walls of $\cW(\varphi_2, \vartheta)$ are roughly the leaves of $CG(\varphi_2, \vartheta)$,
but we have to take a bit of care concerning the orientations, as follows.
There are two possible behaviors for a leaf of $CG(\varphi_2, \vartheta)$.
One possibility is that the leaf has one end on a branch point and the other end 
on a puncture.  In this case the leaf gives a single wall of $\cW(\varphi_2, \vartheta)$,
oriented \ti{away} from the branch point.  The other possibility is that the leaf
has both ends on branch points (saddle connection.)
In this case the leaf gives a double wall of $\cW(\varphi_2, \vartheta)$.\footnote{By perturbing $\vartheta$ slightly to $\vartheta + \eps$ 
this wall would be resolved into two single walls which are close together; 
in the limit as $\eps \to 0^\pm$, we may think of these two walls as becoming 
infinitesimally separated, with the $\pm$ sign determining whether we get the American or the British resolution.}

To be consistent with our previous definitions of the various kinds of spectral networks, we replace a puncture with a boundary (defined by a non-critical leaf of the corresponding singular foliation) if the puncture is surrounded by a polygon consisting of double walls. 

Next we have to explain how the walls are labeled.  Fix a wall $w$ corresponding to a leaf $\gamma(t)$.  
The two sheets of $\Sigma$ over $w$, labeled by the index $i$, 
correspond to the two square roots $\lambda_i$ of $\varphi_2$ over $w$.  
Given a positively oriented tangent vector $v$ to $w$,
we consider the quantity $q_i = e^{-i \vartheta} \lambda_i(v)$.
From \eqref{eq:leaftheta} it follows that $q_i$ is real; it will be positive for one choice of 
square root (say $i = i_+$) and negative for the other (say $i = i_-$).  
We label the wall $w$ by the ordered pair $i_- i_+$.
\footnote{The reason for this choice of 
labeling is as follows \cite{Gaiotto:2009hg}: 
it guarantees that, when we study certain 1-parameter families of connections
$\nabla(\zeta)$ of the form \eqref{eq:wkb-family} below, the corrections $\cS_w$ attached to 
walls of the spectral network will be exponentially \ti{small} as $\zeta \to 0$ ---
thus these corrections will not affect the asymptotics in that limit.}

Finally we describe the decorations of the punctures and annuli.
First consider a puncture $z_l$.  In a neighborhood of 
$z_l$ we may label the sheets of $\Sigma$ by $i = 1, 2$;
each square root $\lambda_i$ of $\varphi_2$ has a simple pole near $z$, with residue 
$r_i = \pm \I m_l$.  
Generically, one of the two sheets will have $\re(e^{-\I \vartheta} r_i) > 0$; the decoration at 
the puncture $z_i$ is obtained by labeling the sheets with this sheet
first.
One then checks by a local analysis of \eqref{eq:leaftheta} \cite{Gaiotto:2009hg} 
that the only walls which fall into 
this puncture are ones whose labeling matches the decoration, as
required.
Next consider an annulus $A$ of $C \setminus \cW(\varphi_2, \vartheta)$.  Such an annulus
is foliated by closed leaves of $\cF(\varphi_2, \vartheta)$.  Fix one such leaf
and choose an orientation $o$ for it; having done so, we can again consider the quantity 
$q_i = \lambda_i(v)$ for $v$ a positively oriented tangent vector.
$q_i$ is positive for one choice of $i$, say $i_+$, and negative for the other,
say $i_-$.  Our decoration of the annulus $A$ assigns
the ordering $i_- i_+$ to the orientation $o$.  This definition is evidently
compatible with the labelings of the walls on the boundary.

\subsection{Generic differentials and Fock-Goncharov networks}

Suppose $C$ has at least one puncture and
$(\varphi_2, \vartheta)$ are chosen generically.
In this case
all leaves of $\cF(\varphi_2, \vartheta)$ end at punctures
\cite{Gaiotto:2009hg}.
The corresponding WKB network $\cW(\varphi_2, \vartheta)$ is a Fock-Goncharov network. 
See Figure~\ref{fig:Nf4-FG-network} for an example.

\insfigscaled{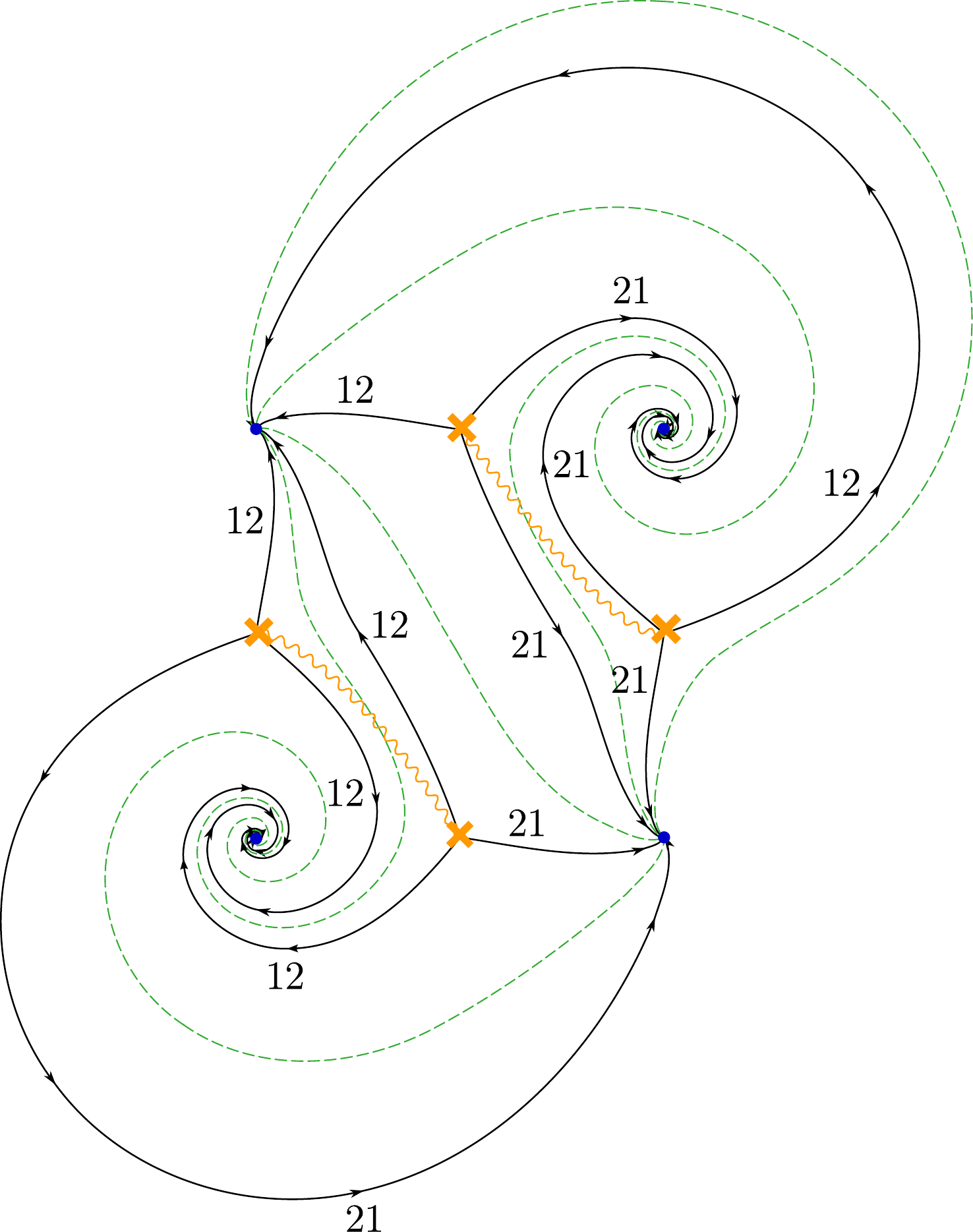}{0.28}{A generic WKB spectral
network $\cW(\varphi_2, \vartheta)$ on the four-punctured sphere, with $\varphi_2 =
   \frac{z^4 - 1}{(z^4+4)^2} \de z^2$
and $\vartheta = \frac{2 \pi}{3}$.  This network is a Fock-Goncharov network; 
the dashed green paths are the edges of the
corresponding (degenerate) ideal triangulation.
}

\subsection{Strebel differentials and Fenchel-Nielsen networks}

Now let us consider the opposite extreme:
suppose that the foliation $\cF(\varphi_2, \vartheta)$ has \ti{no} leaves ending on punctures,
and moreover that all leaves are compact.  Concretely this means that all leaves are either closed 
trajectories or saddle connections.
This is equivalent to saying that $e^{-2 \I \vartheta} \varphi_2$ is a \ti{Strebel differential}.

Under this condition, each component of the complement 
$C \setminus \cW(\varphi_2, \vartheta)$ is a punctured disc, swept out by closed $(\varphi_2, \vartheta)$-trajectories. 
The boundary of each such component is a polygon consisting of one or
more saddle connections. The WKB spectral network $\cW(\varphi_2,
\vartheta)$ contains only  double walls. After we replace all punctures with boundaries, each defined by a closed $(\varphi_2, \vartheta)$-trajectory, it is either a
Fenchel-Nielsen spectral network or a contracted Fenchel-Nielsen
spectral network.  

To each component of the complement $C \setminus \cW(\varphi_2,
\vartheta)$ we can associate a positive number. For any punctured disc $D_j$
this is the coefficient $m_j e^{- \I \vartheta}$ of the Strebel differential at 
the puncture, which is equivalent to the $(\varphi_2,\vartheta)$-length
of any closed horizontal trajectory $\gamma$ in $D_j$.  For any
annulus $R_k$ this is the $(\varphi_2,\vartheta)$-height $h_k$. 

Conversely, given any system of simple closed curves $\{ \alpha_k\}_{1
  \le k \le p}$ of a punctured Riemann surface $C$ and arbitrary $h_k>0$, $ 1 \le k \le p$, 
and $m_j >0$, $1 \le j \le n$, there is a unique 
Strebel differential $\varphi_2$ with the following properties. The
differential $  \varphi_2$ has $n$ punctured discs which are swept out
by closed trajectories around the punctures, that have
$\varphi_2$-length $m_j$. Moreover, it has $p$ characteristic annuli
with type $\alpha_k$ which have $\varphi_2$-height $h_k$. \cite{Liu}

For example, a Strebel differential on the three-punctured
sphere is of the form
\begin{equation}
\varphi_2 = - \frac{ m^2_{\infty} z^2 - (m^2_{\infty} + m^2_0 - m^2_1)
  z + m^2_0 } {z^2 (z-1)^2}~(\de z)^2,
\end{equation}
where all $m_i \in \R$.
This $\varphi_2$ has double poles at the punctures $z_i=0,1,\infty$, 
with residues $-m^2_0,-m^2_1,-m^2_{\infty}$,
respectively. The networks $\cW(\varphi_2,0)$ are depicted in 
Figure \ref{fig:molecules}: we get ``molecule I''
when $m_{\infty} < m_0 + m_1$, and ``molecule II'' when
$m_{\infty} > m_0 + m_1$.

Finally we remark that if we perturb $\vartheta$ slightly, each 
annulus of closed trajectories will be replaced
by a region of $C$ where the trajectories are spiraling around a finite number of times.
As $\vartheta$ is varied across the critical phase $\vartheta_c$ where the annulus appears,
clockwise spiraling leaves are replaced by anti-clockwise spiraling leaves, or vice versa.  
See Figure~\ref{fig:juggle} for an example.
\insfig{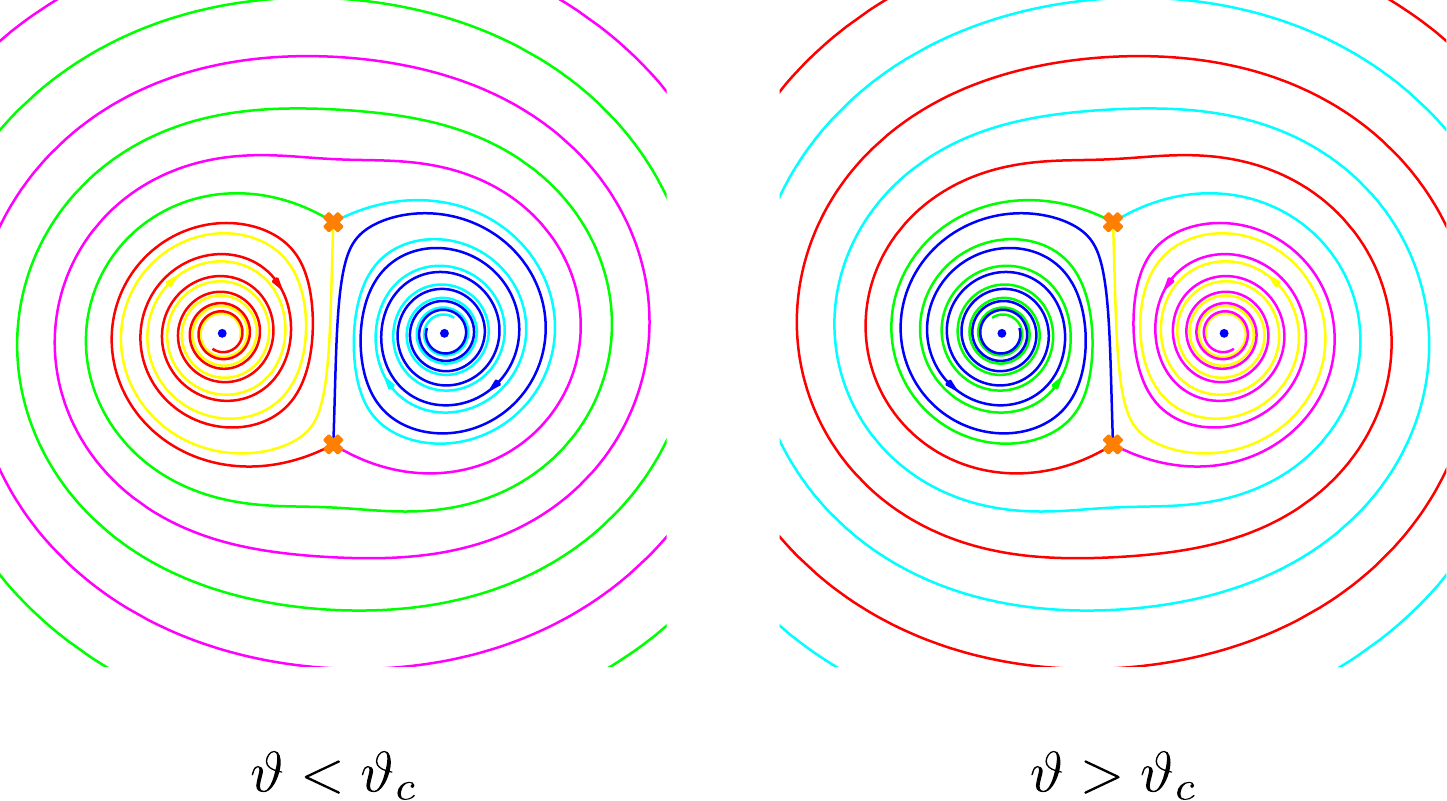}{An example of a 
foliation $\cF(\varphi_2, \vartheta)$ near a critical phase $\vartheta_c$ where
two families of compact leaves appear simultaneously. We label the walls by 
different colors to make it easier to tell them apart.}

\subsection{A counterexample}

We finish this section by noting that not all spectral
networks can arise as WKB spectral networks.
A counterexample is given in
Figure~\ref{fig:FNwithoutdifferential}. If this network 
were a WKB network $\cW(\varphi_2,0)$, then the 
$\varphi_2$-lengths of the critical leaves (measured by integrating 
$\sqrt{-\varphi_2}$) would obey
\begin{align*}
\ell_1 + \ell_6 &= \ell_1 + \ell_2 + \ell_3 + \ell_4, \\
\ell_4 + \ell_5 &= \ell_2 + \ell_3 + \ell_5 + \ell_6.
\end{align*}
Clearly, this is impossible for positive real parameters $\ell_1,
\ldots, \ell_6$.\footnote{This example appears in
  \cite{HubbardMasur}, as an example of a measured foliation
  which cannot be
  generated by a quadratic differential. Although the measured
  foliation is equivalent to one that can be generated by a quadratic
  differential, this is not the case for the associated spectral network.}

\insfigscaled{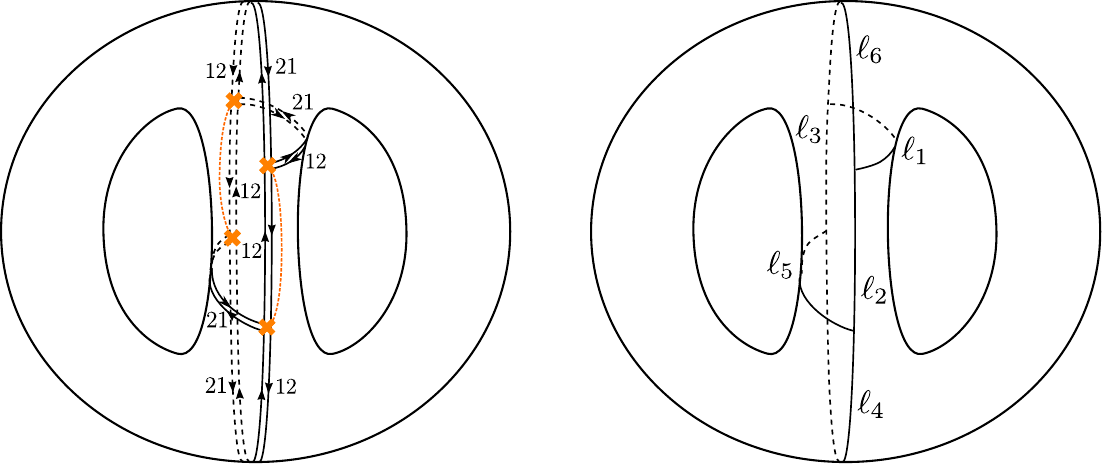}{0.7}{Left:  a spectral
  network on a genus $2$ surface that cannot arise as a WKB network $\cW(\varphi_2, 0)$
for any $\varphi_2$.  Right: critical leaves of the corresponding foliation.}

\section{Abelianization} \label{sec:abelianization}

What is a spectral network good for?
Here is one answer:  given a Fenchel-Nielsen or Fock-Goncharov spectral network $\cW$, and a generic
flat $SL(2)$-connection $\nabla$ in a complex rank 2 vector bundle $E$ over $C$, we can construct a \ti{\cW-abelianization} of $\nabla$.
This is the key to the construction of the ``spectral coordinate'' systems on the moduli
of flat $SL(2)$-connections.  In this section we explain this construction.

\subsection{$\cW$-pairs} \label{sec:what-is-abelianization}

A $\cW$-abelianization of $\nabla$ is a way of putting $\nabla$ in almost-diagonal form, by 
locally decomposing $E$ as a sum of two line bundles, $\cL_1 \oplus \cL_2$, 
which are preserved by $\nabla$.
It is impossible to do this on all of $C$, but we will come close to that ideal, 
in the following sense.

Let $\cW$ be a spectral network subordinate to the covering $\pi: \Sigma \to C$.
Let $\Sigma'$ denote $\Sigma$ with the branch
points removed.
Then:

\begin{itemize}
\item A \ti{$\cW$-pair} is a tuple $(E, \nabla, \cL, \nabla^\ab, \iota)$
consisting of:
\begin{itemize}
 \item a complex rank $2$ vector bundle $E$ over $C$,
 \item a flat $SL(2)$-connection $\nabla$ in $E$, 
 \item a complex line bundle $\cL$ over $\Sigma'$,
 \item an flat $GL(1)$-connection $\nabla^\ab$ in $\cL$,
 \item an isomorphism $\iota: E \simarrow \pi_* \cL$ defined over $C \setminus \cW$,
\end{itemize}
obeying the conditions that
\begin{itemize}
 \item the isomorphism $\iota$ takes $\nabla$ to $\pi_* \nabla^\ab$,
 \item $\iota$ jumps by a map $\cS_w = 1 + e_w \in \End(\pi_* \cL)$ at each single wall $w \subset \cW$,
where $e_w: \cL_i \to \cL_j$ if $w$ carries the label $ij$.  (Here by $\cL_i$ we mean the summand of $\pi_* \cL$ associated to sheet $i$.
Relative to diagonal local trivializations of $\pi_* \cL$, this condition says
$\cS_w$ is upper or lower triangular.) $\iota$ jumps by a map $\cS_{w'} \cS_{w}$ at each double wall $w' w$. 
\end{itemize}

\item Given a flat $SL(2)$-connection $\nabla$ in a complex rank $2$ bundle $E$ over $C$,
a \ti{$\cW$-abelianization} of $\nabla$ is any extension of 
$(E, \nabla)$ to a $\cW$-pair $(E, \nabla, \cL, \nabla^\ab, \iota)$.

\item Given a flat $GL(1)$-connection $\nabla^\ab$ in a complex line bundle $\cL$ over $\Sigma'$,
a \ti{$\cW$-nonabelianization} of $\nabla^\ab$ is any extension of 
$(\cL, \nabla^\ab)$ to a $\cW$-pair $(E, \nabla, \cL, \nabla^\ab, \iota)$.

\end{itemize}

We say $(E, \nabla, \cL, \nabla^\ab, \iota) \sim (E', \nabla', \cL', \nabla'^\ab, \iota')$,
and call these two $\cW$-pairs \ti{equivalent},
if there exist maps $\varphi: \cL \to \cL'$ and $\psi: E \to E'$ which take $\iota$ into $\iota'$,
$\nabla^\ab$ into $\nabla'^\ab$, and $\nabla$ into $\nabla'$.
In particular, in this case we have equivalences $(\cL, \nabla^\ab) \sim (\cL', \nabla'^\ab)$
and $(E, \nabla) \sim (E', \nabla')$.

Since the covering $\Sigma \to C$ has at least one branch point, every complex line bundle
over $\Sigma'$ is topologically trivial.  Thus, in particular, every $\cW$-pair is equivalent to 
one where $\cL$ is taken to be the trivial bundle, so we could have taken $\cL$ to be trivial
from the start; nevertheless this requirement is a bit unnatural from our point of view, and we prefer 
to emphasize that $\cL$ need not be \ti{canonically} trivial.

\subsection{Equivariant $GL(1)$-connections}\label{ssec:propab}

In this section we will show that the connections $\nabla^\ab$ which arise in $\cW$-pairs 
automatically carry some extra structure:  they are \ti{equivariant}, and as a consequence
of this, they are also \ti{almost-flat} (both terms to be defined below.)

Suppose given a  $\cW$-pair $(E, \nabla, \cL, \nabla^\ab, \iota)$.
The fact that $\nabla$ is an $SL(2)$-connection means that
the underlying bundle $E$ carries a $\nabla$-invariant, 
nonvanishing volume form $\varepsilon_E \in \wedge^2(E)$. 
Introduce the 
covering involution $\sigma: \Sigma \to \Sigma$. We have a
$\nabla^\ab$-invariant nondegenerate pairing 
\begin{equation} \label{eq:pairing}
\mu: \cL \otimes \sigma^* \cL \to \C
\end{equation}
given by\footnote{Here we abuse notation slightly: $\iota$ refers to the induced map $\iota: \pi^*E \to \pi^* \pi_* \cL \cong \cL \oplus \sigma^* \cL$ on $\Sigma\backslash \pi^{-1}(\cW)$.}
\begin{equation}
	\mu(s, s') = (\iota^{-1}(s) \wedge \iota^{-1}(s')) / \varepsilon_E. 
\end{equation}
The triangular property of the jumps $\cS_w$ shows that $\mu$
extends even over $\pi^{-1}(\cW)$.
In particular, $\mu$ gives a $\nabla^\ab$-invariant isomorphism
$\cL \simeq \sigma^* \cL^*$.

The antisymmetry of the wedge product means $\mu$ has 
a slightly tricky ``antisymmetry'' property:
if we consider
\begin{equation}
	\sigma^* \mu: \sigma^* \cL \otimes \cL \to \C
\end{equation}
using $\sigma^* \sigma^* \cL = \cL$,
and let $\tau: \sigma^* \cL \otimes \cL \to \cL \otimes \sigma^* \cL$
be the operation of reversing factors, we have
\begin{equation}
	\sigma^* \mu = - \mu \circ \tau.
\end{equation}
We call a connection $\nabla^\ab$ in a line bundle $\cL$
equipped with such a pairing \ti{equivariant}.

The equivariance of $\nabla^\ab$ has an important consequence:
it implies that the holonomy of $\nabla^\ab$
around a small loop on $\Sigma'$ encircling a branch point is always $-1$.
To see this, fix a basepoint $z$ near the branch point.
At this point the fiber of $\pi_* \cL$ is canonically decomposed
into two lines, labeled by the two preimages of $z$ in $\Sigma'$.
The holonomy of $\nabla^\ab$ around the loop exchanges these
two lines; so with respect to a basis $(s_1,s_2)$ respecting this decomposition,
it would be represented by an antidiagonal matrix.
Moreover the equivariance implies that this holonomy has determinant $1$.
Combining these facts, the holonomy of $\pi_* \nabla^\ab$ around the branch point 
can be represented by a matrix of the shape
\begin{equation}
\begin{pmatrix} 0 & D \\ -D^{-1} & 0 \end{pmatrix}. 
\end{equation}
The holonomy of $\pi_* \nabla^\ab$ around a path going twice around the branch point on $C$ is thus
\begin{equation}
\begin{pmatrix}
0 & D \\ -D^{-1} & 0
\end{pmatrix}
^2 = \begin{pmatrix} -1 & 0 \\ 0 & -1 \end{pmatrix}
\end{equation}
This says that the holonomy of $\nabla^\ab$ around each branch point on
$\Sigma$ is $-1$, as desired.
We refer to a connection
$\nabla^\ab$ with this property as an \ti{almost-flat} connection over $\Sigma$.

\subsection{$\cW$-pairs with boundary and gluing} \label{sec:w-pairs-with-boundary}

If $C$ has boundary, it is sometimes useful to consider connections 
and $\cW$-pairs with a bit 
of extra structure. We fix a marked point on each 
boundary component of $C$. Then, a \ti{$\cW$-pair with boundary}
consists of
\begin{itemize}
\item A $\cW$-pair $(E, \nabla, \cL, \nabla^\ab, \iota)$ as in \S\ref{sec:what-is-abelianization},
\item a trivialization of $E_z$ for each marked point $z$,
\item a trivialization of $\cL_{z_i}$ for each preimage $z_i \in \pi^{-1}(z)$ for a marked point $z$,
\item a trivialization of the covering $\Sigma$ over each marked
point $z$,
\end{itemize}
such that $\iota$ maps the trivialization of $E_z$ to
the trivialization of $\pi_* \cL_{z}$ induced from those 
of $\cL_{z_i}$ and $\Sigma$.

Given two surfaces $C$, $C'$ with boundary we can glue along a boundary 
component, in such a way that the marked points are identified.
Now suppose that we have a $\cW$-pair with boundary on each 
of $C$ and $C'$,
and that the monodromies around the glued component are the same (when
written relative to the given trivializations at the marked points).
Then using the trivializations we can glue the two $\cW$-pairs to 
obtain a $\cW$-pair over the glued surface.

\section{Constructing abelianizations}\label{sec:constructingabelianizations}

\subsection{A preliminary observation}

We begin with an observation which will be useful for the 
construction of $\cW$-pairs in the next few sections.
Suppose that in the definition of $\cW$-pair we require only that $\nabla^\ab$
is defined on $\Sigma' \setminus \pi^{-1}(\cW)$.  In fact this leads to an equivalent  definition: $\nabla^\ab$ automatically extends from $\Sigma' \setminus \pi^{-1}(\cW)$ 
to $\Sigma'$.

To see this, consider a single wall $w$ carrying the label (say) $21$.
Choosing a local trivialization of $\cL$ in some neighborhood of the interior of $\pi^{-1}(w)$,
$\pi_* \cL$ also becomes trivial in this neighborhood.  Then $\iota$ amounts to
two bases $(s_1, s_2)$ and $(s'_1, s'_2)$ for $E$ in a neighborhood of $w$, 
related by an upper-triangular transformation,
\begin{equation}
s'_1 = s_1, \qquad s'_2 = s_2 + \alpha s_1
\end{equation} 
for some function $\alpha$.
These bases diagonalize the connection $\nabla$, so that
we have $\nabla s_i = d_i s_i$ on one side of $w$ (for some closed $1$-forms $d_i$)
and on the other side $\nabla s'_i = d'_i s'_i$.
What we would like to show is that $d'_i = d_i$; this would say 
that $\pi_*\nabla^\ab$ extends over $w$, or in other words, 
$\nabla^\ab$ extends over $\pi^{-1}(w)$. But this is now a straightforward computation.
First, since $s'_1 = s_1$ it follows immediately that $d'_1 = d_1$.
Now 
\begin{equation}{\nabla s'_2 = \nabla (s_2 + \alpha s_1) = \nabla s_2 +
(\de \alpha) s_1 + \alpha \nabla s_1}
\end{equation}
and
\begin{equation}
\nabla s'_2 = d'_2 s'_2 = d'_2(s_2 + \alpha s_1).
\end{equation}
Comparing these two gives
$d'_2 = d_2$, as desired (and $\de \alpha = (d_2 - d_1) \alpha$.)

\subsection{Fock-Goncharov spectral networks}\label{ssec:FGcoord}

Suppose we are given a Fock-Goncharov network $\cW$ and a flat $SL(2)$-connection
$\nabla$ on $C$. How do we $\cW$-abelianize $\nabla$?

Given a flat $SL(2)$-connection $\nabla$ and a puncture $z_l$ on $C$, let a \ti{framing} of $\nabla$ at $z_l$ be a choice of a $\nabla$-invariant line subbundle $\ell_l$ of $E$, in a neighborhood of $z_l$.
We define a \ti{$\cW$-framed connection}
to be a flat $SL(2)$-connection $\nabla$ plus a choice of a framing at each 
puncture, subject to an additional condition:
if two punctures $z_i$, $z_j$ can be connected by a path $\cP$ which does not meet $\cW$, 
and we use the parallel transport along $\cP$ to transport $\ell_i$ and $\ell_j$
to a common point $z \in C$, then at this point we get $\ell_i(z) \neq \ell_j(z)$.
This condition is automatically satisfied for a generic $\nabla$.
Moreover, for generic $\nabla$, there are exactly $2^n$ possible $\cW$-framings,
where $n$ is the number of punctures; indeed
the monodromy around each puncture has two distinct eigenspaces, 
and we are just choosing one of them.

In the rest of this section we will show:
\ti{$\cW$-abelianizations of $\nabla$ up to equivalence
are canonically in 1-1 correspondence with $\cW$-framings of $\nabla$.}

Let $\cL$ be the trivial bundle over $\Sigma$.
Our $\cW$-abelianizations of $\nabla$ will be determined by giving
an isomorphism $\iota: E \simeq \pi_* \cL$ over $C \setminus \cW$.
$C \setminus \cW$ consists of ``cells'' with the topology of
component $(1)$ in Figure \ref{fig:topology-components}; let $c$ denote one
such cell.

Suppose we trivialize the covering $\Sigma$ over $c$.\footnote{We emphasize that this is only a \ti{local} trivialization; we can always do this since $c$ is simply connected. Thus
in what follows we will not have to worry about branch cuts.}
Giving $\iota$ over $c$ is then equivalent to giving $\iota^{-1}: \cO \oplus \cO \to E$,
i.e. a basis $(s_1, s_2)$ of $E$ over $c$. How will we get this basis?
On the boundary of $c$ we have two punctures, each carrying a decoration.
Moreover, our rules for spectral networks imply that the two decorations are opposite:
thus one puncture is labeled $21$ (call this puncture $z_1$) and the
other $12$ (call this puncture $z_2$).
Each puncture has a framing; call these $\ell_1$ and $\ell_2$.
By parallel transporting the framings $\ell_1$ and $\ell_2$ along a path in $c$ to a general $z \in c$,
we obtain lines $\ell_1(z)$ and $\ell_2(z)$ in $E_z$.
Our basis $(s_1, s_2)$ is obtained by choosing
$s_1(z) \in \ell_1(z)$ and $s_2(z) \in \ell_2(z)$. Crucially, for any such choice, $\nabla$ is diagonal
relative to the basis $(s_1, s_2)$; thus $\iota$ indeed identifies $\nabla$ with the pushforward of 
a connection $\nabla^\ab$ in the line bundle $\cL$.

Note that this construction would have failed if $\ell_1(z) = \ell_2(z)$,
since in that case $s_1(z)$ and $s_2(z)$ are not linearly independent, so they do not form a basis of $E_z$.
This is why we required that $\ell_1(z) \neq \ell_2(z)$ in the definition of $\cW$-framed connection.

\insfigscaled{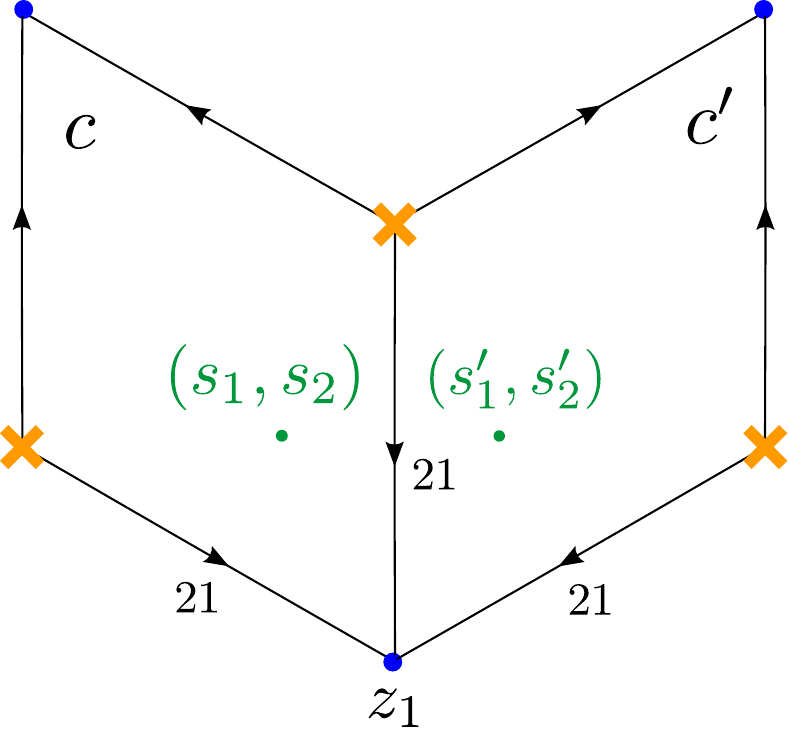}{0.4}{Bases $(s_1,s_2)$ and $(s'_1,s'_2)$ attached to the
  two sides of a wall that separates the cells $c$ and $c'$.}

Now let us consider the change-of-basis matrix $\cS_w$.
Consider the two bases $(s_1, s_2)$ and $(s'_1, s'_2)$
attached to the two sides of the wall in Figure~\ref{fig:FG-coord-sec}.  Because the two cells
$c$, $c'$ have the vertex $z_1$ in common, we will have $\ell_1 = \ell'_1$, so $s_1$
and $s'_1$ can differ at most by scalar multiple.  This implies that the transformation $\cS_w$ taking
$(s_1, s_2)$ to $(s'_1, s'_2)$ has the form $$\cS_w = \begin{pmatrix} * & * \\ 0 & * \end{pmatrix}$$ relative
to the basis $(s_1, s_2)$.
By a diagonal ``gauge transformation'' of the form
\begin{equation}
(s'_1(z), s'_2(z)) \mapsto (\lambda_1(z) s'_1(z), \lambda_2(z) s'_2(z)),
\end{equation}
where $\lambda_1$ and $\lambda_2$ are functions on $c'$,
we can thus arrange that
\begin{equation} \label{eq:S-form-uptri}
\cS_w = \begin{pmatrix} 1 & * \\ 0 & 1 \end{pmatrix}
\end{equation}
relative to the
basis $(s_1, s_2)$.
This is the desired form for the jump of $\iota$ according to the definition of $\cW$-abelianization.

\insfigscaled{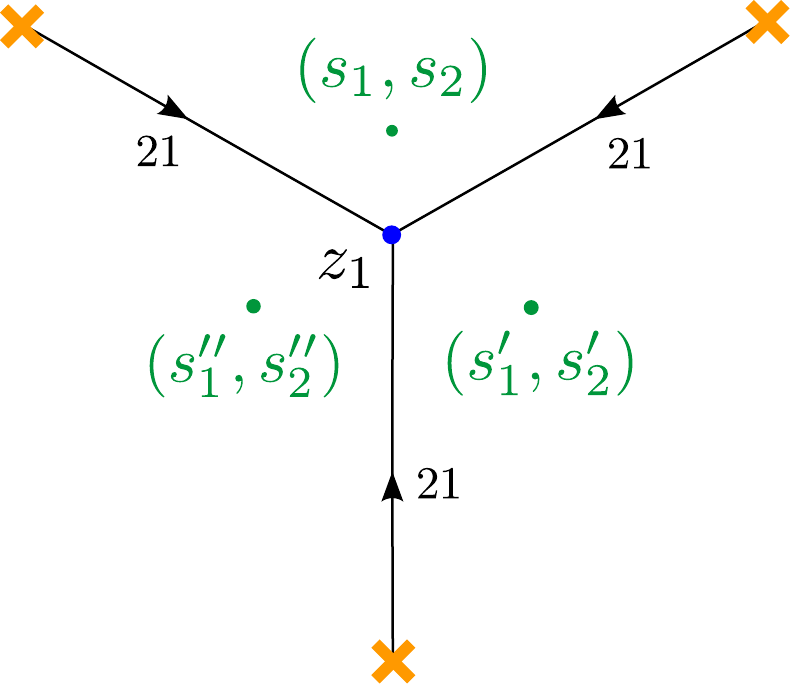}{0.4}{Neighborhood of a puncture with bases.}

We have now described a canonical $\cW$-abelianization for any
$\nabla$ equipped with a $\cW$-framing.
It only remains to show that in this way we obtain \ti{all}
$\cW$-abelianizations of $\nabla$.
So, suppose we have \ti{any} $\cW$-abelianization
of the connection $\nabla$.
Consider a neighborhood of a puncture $z_1$ as shown in Figure~\ref{fig:FG-coord-punct}.  In each cell, the $\cW$-abelianization amounts to
a pair of sections $s_1(z)$, $s_2(z)$, with respect to which $\nabla$ is diagonal.
The constraint that $\cS_w$ is of the form \eqref{eq:S-form-uptri}
says that $\cS_w$ preserves $s_1$, i.e. along the wall $w$ we have $s_1(z) = s'_1(z)$.
Continuing around the wall we see that $s_1, s'_1, s''_1$ are the restrictions of a
single section ${\hat s}_1(z)$ defined in a whole punctured
neighborhood of the puncture $z_1$, such that $\nabla {\hat s}_1(z)$ is a multiple of ${\hat s}_1(z)$.
In particular, ${\hat s}_1(z)$ is an eigenvector of the monodromy of $\nabla$ around the puncture.

Applying this condition to every puncture is almost enough to determine the 
$\cW$-abelianization.
Let us consider the choices remaining. First, there is the 
choice of \ti{which} eigenvector to take at each puncture. This choice 
is equivalent to the choice of a $\cW$-framing.
Second, there is the freedom to make an overall rescaling
$(s_1(z), s_2(z)) \to (\lambda_1(z) s_1(z), \lambda_2(z) s_2(z))$ in each cell 
(matching up along the walls).
Such a change only changes the
$\cW$-abelianization by equivalence.
So altogether we have found that $\cW$-abelianizations are in 1-1 correspondence
with $\cW$-framings, as desired.

\subsection{Fenchel-Nielsen spectral networks}\label{ssec:FN-ab}

Next suppose we are given a Fenchel-Nielsen spectral network $\cW$ and a flat $SL(2)$-connection
$\nabla$ over $C$.

As for Fock-Goncharov networks, we will 
need to enrich $\nabla$ by some extra data in order to get a
canonical $\cW$-abelianization, as follows.
The complement $C' = C \setminus \cW$
consists of various connected components $A$, each with the topology of an annulus.
There are two ways of going around
the annulus; define a \ti{framing} of $\nabla$ on $A$ to be
a 1-dimensional eigenspace of the monodromy in each direction (so if we label the two directions
$\pm$, we give an eigenspace $\ell_+$
of the monodromy $M_+$ and an eigenspace $\ell_-$ of the monodromy
$M_- = M_+^{-1}$.)  We require that $\ell_+ \neq \ell_-$,
and also that each of $\ell_\pm$ for any annulus is distinct from each of $\ell_\pm$ for
any adjacent annulus.\footnote{This last constraint would follow from
the requirement that $\nabla$ is indecomposable when restricted to any pair of
pants. Indeed, suppose that the monodromies around two adjacent annuli
$A$ and $A'$ have a common eigenline $\ell$, and consider a pair of pants that has $A$, $A'$ 
as two of its boundary components; the monodromy around the third boundary of this pair 
of pants will also have eigenline $\ell$.}
Note that a framing of $\nabla$ on $A$ exists only if $M_\pm$ are diagonalizable.
Finally, we define a $\cW$-framed connection to be a flat $SL(2)$-connection $\nabla$
plus a framing of $\nabla$ on each annulus $A$.

We can now proceed much as we did for the Fock-Goncharov networks.  Fix an annulus $A$.
For any $z \in A$, the framing of $\nabla$ on $A$
gives two lines $\ell_+(z) \subset E_z$ and $\ell_-(z) \subset E_z$.
Choose a section $s_+(z)$ valued in $\ell_+(z)$ and a section $s_-(z)$ valued in $\ell_-(z)$.
We want to use this pair of sections to build an isomorphism $\iota: E \to \pi_* \cL$,
where $\cL$ denotes
the trivial bundle over $\Sigma$.  Concretely, the covering $\Sigma$ is trivializable
over $A$ (no branch cut crosses the annulus from edge to edge), so we may choose a trivialization; after so doing,
what we need to do is to choose which of
$s_\pm$ will be $s_1$ and which will be $s_2$.  This extra information is provided
by the decoration of $\cW$ on $A$:  if the labeling attached to the $+$ direction is $ij$,
then we take $s_i = s_+$ and $s_j = s_-$.

\insfigscaled{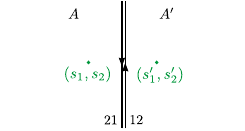}{1.5}{Bases of sections on either side of
  a double wall that separates two annuli $A$ and~$A'$. (Only a piece of each annulus is shown 
  here; for a sample global picture see Figure~\ref{fig:four-holed-sphere}.)
}

Now we need to arrange that the change-of-basis matrices $\cS_w$ attached to the walls
$w \subset \cW$ are of the required triangular form.
Consider the double wall shown in Figure \ref{fig:FN-coord-sec-1},
separating annuli $A$, $A'$.
What we require is that the matrix taking $(s_1, s_2)$ to $(s'_1, s'_2)$ is of the form
\begin{equation} \label{eq:S-form-double}
\begin{pmatrix} 1 & * \\ 0 & 1 \end{pmatrix} \begin{pmatrix} 1 & 0 \\ * & 1 \end{pmatrix}.
\end{equation}
This just says that this matrix has determinant $1$ and lower left entry $1$.
Since $s'_1$ and $s_2$ are not proportional, and $s'_2$ and $s_1$ are not proportional, we know
that neither of the diagonal entries vanishes; thus we can arrange
the desired form by a rescaling
of the form $(s'_1,s'_2) \to (\lambda_1(z) s'_1, \lambda_2(z) s'_2)$.

Finally, we observe that in this way we obtain \ti{all} $\cW$-abelianizations of
$\nabla$. Indeed, if we have an abelianization of $\nabla$ on an annulus,
then the two lines $\iota^{-1}(\cL_1)$ and $\iota^{-1}(\cL_2)$ must be
eigenlines of the monodromy of $\nabla$ around the annulus.
The only freedom is the choice of which of these lines is which
eigenline, i.e. the choice of framing. Thus $\cW$-abelianizations
of $\nabla$ are in 1-1 correspondence with $\cW$-framings of $\nabla$.

\subsection{Fenchel-Nielsen with boundary}

Finally we consider a small extension of the above.
Suppose that $C$ is a surface with boundary, carrying a
Fenchel-Nielsen network $\cW$, and a flat $SL(2)$-connection
$\nabla$ in a bundle $E$.
Suppose moreover that we fix a
marked point $z_a$ on each boundary component, along with a
trivialization of $E_{z_a}$, with respect to which the monodromy
is diagonal.
In this case our construction of a $\cW$-pair in \S\ref{ssec:FN-ab} 
is easily extended to construct a $\cW$-pair with boundary, 
in the sense of \S\ref{sec:w-pairs-with-boundary}.

Suppose we are given two surfaces $C$, $C'$ with boundary and a flat
$SL(2)$-connection on each of $C$ and $C'$ with trivialization at a
marked 
  point. If the monodromies around the glued component are the same
  (when written relative to the given trivializations at the marked
  points), then we can either first $\cW$-abelianize and then glue the
  resulting $\cW$-pairs (as in \S\ref{sec:w-pairs-with-boundary}), or
  first glue the two flat $SL(2)$-connections and then
  $\cW$-abelianize. Since the monodromy of the 
  $SL(2)$-connection is diagonal along the boundary component, it is
  clear that the two actions commute.

\subsection{Mixed spectral networks}\label{sec:abelianizationmixed}

In the last two sections we have described the process of $\cW$-abelianization when $\cW$
is a Fock-Goncharov or a Fenchel-Nielsen network.  More generally we could take $\cW$ to be a mixed
spectral network.  Such networks lie somewhere between the two extremes just considered.
We expect that all of our constructions have analogues when $\cW$ is mixed; indeed, all of our
constructions were \ti{local}, involving just a particular cell or a particular (single or double) wall.
Each cell of $C \setminus \cW$ has one of the two topologies we have just dealt with above, so we can
determine the abelianization on each cell using the same recipes we used above.
Moreover, by appropriate rescalings of the sections $(s_1, s_2)$ on each cell, we will be able to arrange
the desired form for the change-of-basis matrices $\cS_w$ at each wall, just as we did above.

\section{The integrated version}\label{sec:integratedversion}

Instead of working directly with flat connections,
it will be convenient in the rest of the paper to 
work with the corresponding integrated objects,
namely the parallel transport maps. 
While this change of emphasis is not strictly necessary, the
integrated point of view is
very convenient for practical computations,
as we will see in the examples of \S\ref{sec:examples} below.

In this section we spell out the details of this
integrated version of the story.
Thus we will 
replace flat $SL(2)$-connections on $C$ by $SL(2)$-representations 
of a certain groupoid $\cG_C$ of paths on $C$, defined below.
For the $GL(1)$-connections over $\Sigma$ we have to do something 
slightly more complicated, to capture the extra condition that 
the $GL(1)$-connections we use are 
equivariant (as we saw in \S\ref{ssec:propab}).
Thus we will define below a notion of 
\ti{equivariant $GL(1)$-representation} of a groupoid
$\cG_{\Sigma'}$ of paths on $\Sigma'$.
Finally, we will define a notion of \ti{representation $\cW$-pair}
which is the integrated counterpart of an ordinary $\cW$-pair.

\subsubsection*{Path groupoids}

First we define the path groupoid $\cG_C$.
We fix a set
$\cP_C$ of \emph{basepoints} on $C$, as follows.
For each single wall $w$ we fix two basepoints, which lie
very close to one another, one on each side of $w$.
For a double wall we fix three basepoints, one on each side 
of the wall and one between the two lanes. 
Finally, we fix one basepoint along each boundary component of $C$.
See Figure~\ref{fig:basepoints} for an example.
\insfigscaled{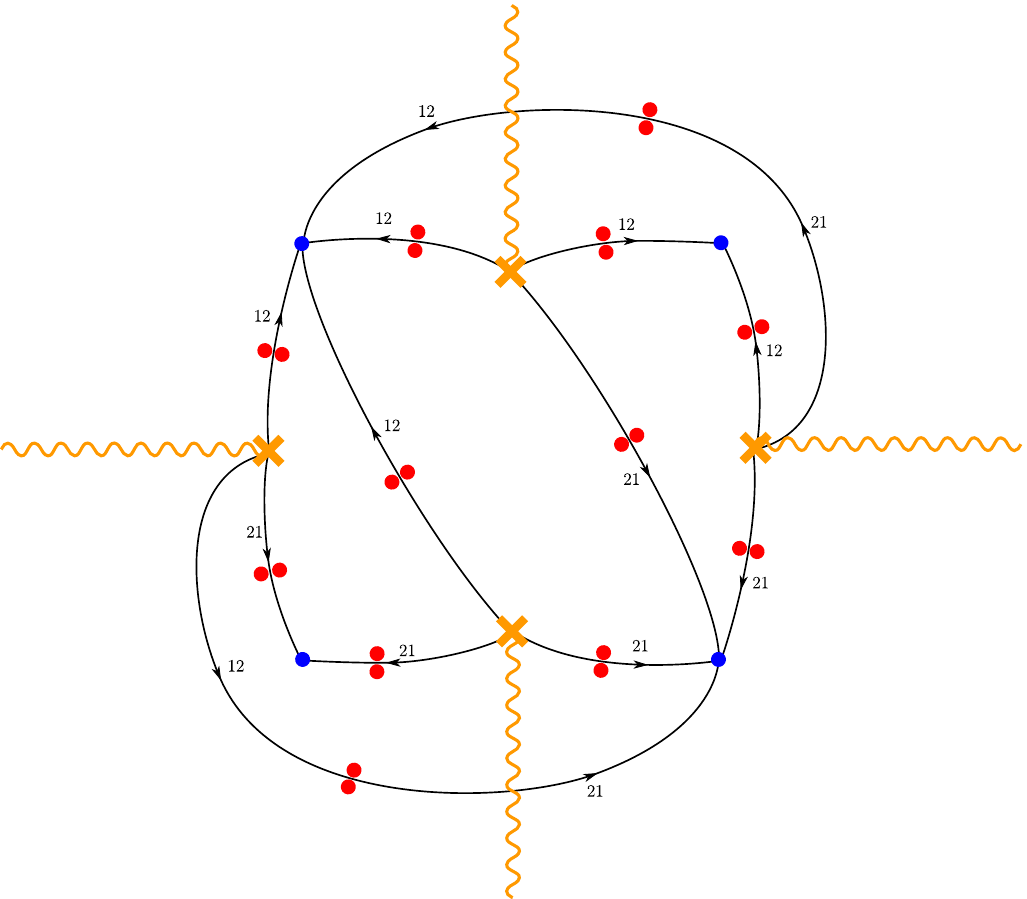}{0.6}{Example of a spectral network with basepoints.}

$\cG_C$ is the groupoid whose objects are the points of $\cP_C$,
and whose morphisms are homotopy classes of
paths $\wp$ which begin and end at points of $\cP_C$.

Next we define the path groupoid $\cG_{\Sigma'}$ on the double cover.
Let $\cP_{\Sigma'} = \pi^{-1}(\cP_C)$: thus for each $z
\in \cP_C$ we have two elements $z^{(i)} \in \cP_{\Sigma'}$.
The objects in $\cG_{\Sigma'}$ are just the points of $\cP_{\Sigma'}$.
The morphisms in $\cG_{\Sigma'}$ are homotopy classes of paths on $\Sigma'$
which begin and end at points of $\cP_{\Sigma'}$.

\subsubsection*{$SL(2)$-representations}

By an $SL(2)$-representation $\rho$ of $\cG_C$, we mean
\begin{itemize}
\item a $2$-dimensional vector space $E_z$ for each
$z \in \cP_C$, together with a nonzero element
$\varepsilon_z \in \wedge^2(E_z)$,
\item an element $\rho(\wp): E_{i(\wp)} \to E_{f(\wp)}$
for each morphism $\wp$ of $\cG_C$,
with initial point $i(\wp)$ and final point $f(\wp)$;
$\rho$ should be compatible with composition, and each
$\rho(\wp)$ should be compatible with $\varepsilon$ in the
obvious sense, i.e. have ``determinant $1$''.
\end{itemize}

An \ti{equivalence} between 
$SL(2)$-representations $\rho$, $\rho'$ is a
collection of isomorphisms
\begin{equation}
g_z: E_z \to E'_z
\end{equation}
such that for any $\wp \in \cG_C$
\begin{equation}
\rho'(\wp) = g_{f(\wp)} \rho(\wp) g_{i(\wp)}^{-1}
\end{equation}
and such that $g_z(\varepsilon_z) = \varepsilon'_z$.
Thus $SL(2)$-representations of $\cG_C$ form a category.

There is an equivalence between
the category of $SL(2)$-representations of $\cG_C$ and the category
of flat $SL(2)$-connections over $C$.
For completeness we briefly explain how this goes.
One direction of the equivalence is easy to describe: given a flat
$SL(2)$-connection $\nabla$ in a bundle $E$ over $C$, define an $SL(2)$-representation
$\rho$ of $\cG_c$
by the rule that $E_z$ is the fiber of $E$ over $z$, 
and $\rho(\cP)$ is the parallel transport of $\nabla$ along $\cP$.
For the other direction, we consider the based path space $\tilde{C}$,
consisting of homotopy classes of paths on $C$ which begin at points
of $\cP_C$ and end anywhere. $\tilde{C}$ has one connected
component $\tilde{C}_z$ for each $z \in \cP_C$. One can recover $C$ from $\tilde{C}$
by dividing out an equivalence relation: two paths $p$, $p'$ are equivalent
if they differ by concatenation with some path $\wp \in \cG_C$, $p = p' \wp $.\footnote{Notice that our notation for composing paths is
  opposite from usual: if the end point of $\wp_1$ agrees with the
  initial point of $\wp_2$, we write their composition as $\wp_2
  \wp_1$. \label{footnote:composepaths}}  
Now suppose given an $SL(2)$-representation $\rho$ of $\cG_C$.
We consider an
$SL(2)$-bundle $\tilde{E}$ over $\tilde{C}$, with flat connection, defined as follows: 
over the component $\tilde{C}_z$, $\tilde{E}$ is
the trivial bundle $E_z \times \tilde{C}_z$, with its trivial connection.
The desired $SL(2)$-bundle $E$ over $C$ with flat connection 
is then obtained by dividing out 
an equivalence relation: if $p = p' \wp$, we use $\rho(\wp)$ to identify
the fibers of $\tilde{E}$ over $p$ and $p'$.

\subsubsection*{Equivariant $GL(1)$-representations}\label{ssec:goodrep}

By an equivariant $GL(1)$-representation $\rho^\ab$ of $\cG_{\Sigma'}$ we mean
\begin{itemize}
\item a $1$-dimensional vector space $\cL_z$ for each
$z \in \cP_{\Sigma'}$,
\item an isomorphism $\rho^\ab(\wp): \cL_{i(\wp)} \to \cL_{f(\wp)}$
for each morphism $\wp$ of $\cG_{\Sigma'}$,
with initial point $i(\wp)$ and final point $f(\wp)$;
$\rho^\ab$ should be compatible with composition,
\item a nondegenerate pairing $\mu_z: \cL_z \otimes \cL_{\sigma(z)} \to \C$ for each
$z \in \cP_\Sigma$, with the antisymmetry property
$\mu(s,s') = - \mu(\sigma^*s', \sigma^* s)$, and
compatible with $\rho^\ab$ in the obvious sense.
\end{itemize}

In parallel to the previous section we can also define
\ti{equivalence} of equivariant $GL(1)$-representations:
the only new point is that an equivalence between $\rho^\ab$ 
and $\rho'^\ab$ is required to 
intertwine $\mu$ and $\mu'$ in the obvious sense.
Then, again in parallel to the previous section,
there is an equivalence between the category of
equivariant $GL(1)$-representations of $\cG_{\Sigma'}$ 
and the category
of equivariant $GL(1)$-connections over $\Sigma'$. (The construction
is just as above, now using the 
path space $\tilde \Sigma'$ instead of $\tilde C$; 
the equivariant structure goes through
straightforwardly.) 

Just as for equivariant $GL(1)$-connections
over $\Sigma'$,
equivariant $GL(1)$-representations are automatically
almost-flat, in the sense that they assign holonomy
$-1$ to a loop around a branch point.

\subsubsection*{Representation $\cW$-pairs}

Now we are ready to define the integrated analogue of a $\cW$-pair.
These definitions will be parallel to those of \S\ref{sec:what-is-abelianization}.

As before, suppose we have fixed a spectral network $\cW$ subordinate
to a covering $\Sigma$. Then:
\begin{itemize}
\item A \ti{representation $\cW$-pair} is a tuple $(\rho,
\rho^\ab, \iota)$, consisting of: 
\begin{itemize}
\item an $SL(2)$-representation $\rho$ of $\cG_C$,
\item an equivariant $GL(1)$-representation $\rho^\ab$ of $\cG_{\Sigma'}$,
\item an isomorphism $\iota_z: E_z \to (\pi_* \cL)_z$ for each 
$z \in \cP_C$,
\end{itemize}
such that
\begin{itemize}
\item if $\wp$ does not cross any walls of $\cW$, then $\iota$
intertwines $\pi_* \rho^\ab(\wp)$ with $\rho(\wp)$:
\begin{equation} \label{eq:intertwiner-1}
\rho(\wp) = \iota_{f(\wp)}^{-1} \circ \pi_* \rho^\ab(\wp) \circ \iota_{i(\wp)},
\end{equation}
\item if $\wp$ is a short path connecting the two basepoints
attached to a wall $w$ of $\cW$, then
$\iota$ intertwines $\pi_* \rho^\ab(\wp)$ with $\rho(\wp)$
up to a unipotent correction:
\begin{equation} \label{eq:intertwiner-2}
\rho(\wp) = \iota_{f(\wp)}^{-1} \circ \cS_w \circ \pi_* \rho^\ab(\wp) \circ \iota_{i(\wp)} 	
\end{equation}
where $\cS_w$ is a unipotent endomorphism of $(\pi_* \cL)_z$: if $w$ carries the 
label $ij$ then we have $\cS_w = 1 + e_w$, where $e_w: \cL_{z^{(i)}} \to \cL_{z^{(j)}}$,
\item the pairing $\mu_z$ corresponds to the volume form $\varepsilon_z$ under 
the isomorphism $\iota_z$.

\end{itemize}
\item Given a flat $SL(2)$-representation $\rho$ of the path groupoid
  $\cG_C$ on $C$, a \emph{$\cW$-abelianization} of $\rho$ is any
  extension of $\rho$ to a representation $\cW$-pair $(\rho,
\rho^\ab, \iota)$. 
\item Given an equivariant $GL(1)$-representation $\rho^\ab$ of $\cG_{\Sigma'}$ on
  $\Sigma$, a \emph{$\cW$-nonabelianization} of $\rho^\ab$ is any
  extension of $\rho^\ab$ to a representation $\cW$-pair $(\rho,
\rho^\ab, \iota)$. 
\end{itemize}

An \ti{equivalence} between representation $\cW$-pairs
$(\rho, \rho^\ab, \iota)$ and $(\rho',
\rho'^\ab, \iota')$ 
is a pair $(g,g^\ab)$ where 
$g$ is an equivalence between $\rho$ and $\rho'$ as defined above, and $g^\ab$ 
an equivalence between $\rho^\ab$ and
$\rho'^\ab$ as defined above, such that $g$ and $g^\ab$ intertwine $\iota$ with $\iota'$. 

We have now defined a category of representation $\cW$-pairs.
Naturally, the point of the definition is that this category
is equivalent to the category of ordinary $\cW$-pairs which we 
defined in \S\ref{sec:what-is-abelianization} above.
The construction of this equivalence is completely parallel
to the constructions for the individual categories of representations
which we considered above.

\subsubsection*{Representations with boundary and gluing}

If $C$ has boundary, it is sometimes useful to consider
representations and representation $\cW$-pairs with a bit of extra
structure. The definitions in \S\ref{sec:w-pairs-with-boundary}
translate to the following.  

\begin{itemize}
\item
An $SL(2)$-representation with boundary consists of an
$SL(2)$-representation together with a trivialization of $E_z$ for
each marked 
point $z$. 
\item
An equivariant $GL(1)$-representation with boundary consists of an
equivariant $GL(1)$-representation together with a trivialization of
$\cL_z$ for each marked point $z$. 
\item 
A representation $\cW$-pair with boundary is a tuple $(\rho,
\rho^{\ab}, \iota)$ consisting of an $SL(2)$-representation with
boundary $\rho$, an equivariant $GL(1)$-representation 
with boundary $\rho^{\ab}$, and a trivialization of the covering $\Sigma$ over each
marked point $z$
such that the isomorphism $\iota$ maps the trivialization of $E_z$ to
the trivialization of $\pi_* \cL_{z_i}$ induced from those of
$\cL_{z_i}$ and $\Sigma$.
\end{itemize}

Given two surfaces $C$, $C'$ with boundary, we can then define
the gluing of two $SL(2)$-representations or two representation
$\cW$-pairs with boundary, parallel to the discussion in
\S\ref{sec:w-pairs-with-boundary}.

\section{Constructing nonabelianizations}\label{sec:nonabelianization}

In this section, given a Fock-Goncharov or Fenchel-Nielsen
spectral network $\cW$ subordinate to a cover $\Sigma$,
and an equivariant $GL(1)$-representation $\rho^\ab$ of $\cG_{\Sigma'}$,
we construct a $\cW$-nonabelianization of $\rho^\ab$.

First, for $z \in \cP_C$ we define
\begin{equation}
	E_z = (\pi_* \cL)_z,
\end{equation}
and let $\iota_z$ be the identity map.
(Up to equivalence, this choice does not involve any loss of generality: 
any $\cW$-nonabelianization of $\rho^\ab$
is trivially equivalent to one for which $E_z = (\pi_* \cL)_z$
and $\iota_z$ is the identity.) The pairing $\mu_z$ induces a
volume element $\varepsilon_z \in \wedge^2(E)$.

To complete the $\cW$-nonabelianization, it only remains to determine 
the parallel transports
$\rho(\wp)$. These are constrained by the intertwining
relations \eqref{eq:intertwiner-1}, \eqref{eq:intertwiner-2},
which \ti{almost} determine them in terms of the given $\rho^\ab$:
the only remaining ambiguity is in the choice of the
$e_w$ attached to the walls. Each $e_w$ lies in a $1$-dimensional
vector space $Hom(\cL_i, \cL_j)$.

If we choose the $e_w$ arbitrarily, the resulting 
elements $\rho(\wp)$ will not necessarily define an honest
representation of $\cG_C$: they will not satisfy the condition
of homotopy invariance. Concretely, this condition amounts to 
requiring that $\rho(\wp) = 1$ whenever $\wp$ is a contractible
loop on $C$. 
A general contractible loop $\wp$ on $C$ can be written as 
a composition of 
loops which do not cross any walls and loops which encircle a single branch point.
If $\wp$ does
not cross any walls, then $\rho(\wp) = 1$ 
follows simply from \eqref{eq:intertwiner-1}
and the fact that $\rho^\ab$ is an honest representation of $\cG_{\Sigma'}$.
So it only remains to consider loops $\wp$ which encircle a single branch point:
for each such loop, the condition $\rho(\wp) = 1$ gives a nontrivial constraint 
on the $e_w$, which we call a \ti{branch point constraint}.

We will now show that,
for all Fock-Goncharov or Fenchel-Nielsen networks $\cW$, the
branch point constraints admit a unique solution,
and thus completely determine the $e_w$.
For a Fock-Goncharov network, we consider the path 
$\wp$ in Figure \ref{fig:around_branchpoint_symm}.
\insfigscaled{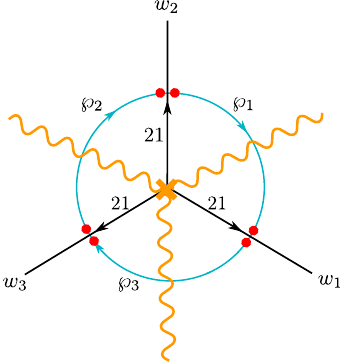}{1.1}{A small loop around a branch point in 
a Fock-Goncharov network.}
Since $\wp$ is contractible, we must have 
\begin{equation} \label{eq:rho-contract}
\rho(\wp) = 1.
\end{equation}
On the other hand, using \eqref{eq:intertwiner-1}, \eqref{eq:intertwiner-2}
we have 
\begin{align}
\rho(\wp) &= \cS_{w_1} \circ \pi_* \rho^\ab(\wp_1) \circ \cS_{w_2} \circ \pi_* \rho^\ab(\wp_2) \circ \cS_{w_3} \circ \pi_* \rho^\ab(\wp_3) \\
&= \begin{pmatrix} 1 & e_{w_1} \\ 0 & 1 \end{pmatrix}  \begin{pmatrix} 0 & D_1 \\ D'_1 & 0 \end{pmatrix} \begin{pmatrix} 1 & e_{w_2} \\ 0 & 1 \end{pmatrix}  \begin{pmatrix} 0 & D_2 \\ D'_2 & 0 \end{pmatrix} \begin{pmatrix} 1 & e_{w_3} \\ 0 & 1 \end{pmatrix}  \begin{pmatrix} 0 & D_3 \\ D'_3 & 0 \end{pmatrix} \label{eq:product-out}
\end{align}
where each matrix is written with respect to the decomposition 
of $\pi_* \cL$ into the two lines $\cL_i$, and
$D_n$, $D'_n$ are obtained by applying $\rho^\ab$ to the
two lifts of $\wp_n$.
Combining \eqref{eq:rho-contract} and \eqref{eq:product-out},
and using the fact that $D_1 D'_2 D_3 D'_1 D_2 D'_3 = -1$,
gives the unique solution
\begin{equation}
e_{w_1} = -D_1 D'_2 D_3, \qquad e_{w_2} = -D_2 D'_3 D_1, \qquad e_{w_3} = -D_2 D'_1 D_3.
\end{equation}
Incidentally, this solution has a nice interpretation: it says that 
$-e_w$ is obtained by applying $\rho^\ab$ to a ``detour'' path on $\Sigma$, 
shown in Figure \ref{fig:detour} for the wall $w = w_1$. This was 
the point of view taken in 
\cite{Gaiotto2012}.
\insfigscaled{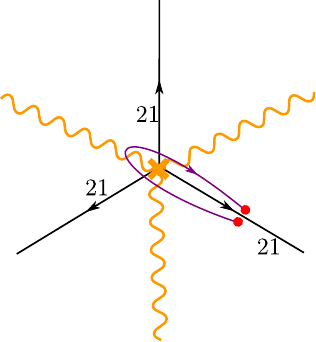}{1.1}{A ``detour'' path on $\Sigma$.}

Next suppose we have a Fenchel-Nielsen network. In this
case the walls are grouped into molecules of the shapes 
shown above in Figure \ref{fig:molecules}. Let us consider
molecule II 
as shown in Figure \ref{fig:FN-local-computation}.
\insfigscaled{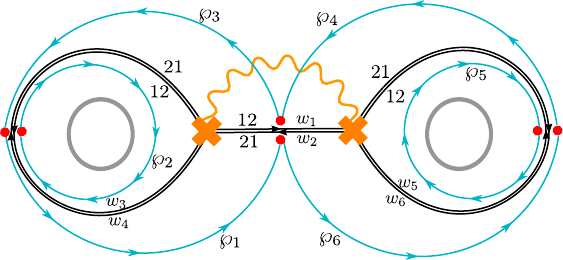}{1.1}{Local configuration around
  branch-points in a Fenchel-Nielsen molecule.}
Here we have $6$ undetermined $e_{w_1}, \dots, e_{w_6}$,
$2$ constraint equations (see footnote~\ref{footnote:composepaths} for our convention on composing paths)
\begin{equation} \label{eq:constr-FN}
\rho(\wp_1 \wp_2 \wp_3) = 1, \qquad \rho(\wp_4 \wp_5 \wp_6) = 1,
\end{equation}
and relations $D_1 D'_1 D_2 D'_2 D_3 D'_3 = -1$ and 
$D_4 D'_4 D_5 D'_5 D_6 D'_6 = -1$, just in parallel to the above.
The equations \eqref{eq:constr-FN} determine the $e_{w_i}$ uniquely;
for example one gets
\begin{equation}
e_{w_3} = D_1 D'_3 \frac{1}{D_2 - D'_2} 
\end{equation}
where in writing the last factor we use the fact that $D_2$ and $D'_2$
are just \ti{numbers}, obtained by applying $\rho^\ab$ to \ti{closed} paths.
The other $e_{w_i}$ are given by similar (slightly more complicated)
expressions.
 
Thus, for any Fock-Goncharov or Fenchel-Nielsen network $\cW$, 
there exists a (unique up to equivalence)
$\cW$-nonabelianization for
any $\rho^\ab$. Moreover, all our constructions were canonical;
from this it easily follows that given an equivalence between
equivariant $GL(1)$-connections $\rho^\ab$ and $\rho'^\ab$
we get an equivalence of their $\cW$-nonabelianizations.

\subsubsection*{Fenchel-Nielsen with boundary}

As a small extension of the above, suppose that $\cW$ is a
Fenchel-Nielsen network subordinate to a cover $\Sigma$ with
boundary. Suppose moreover that we are given an equivariant
$GL(1)$-representation with boundary. Then we can easily extend 
the nonabelianization construction to construct a $\cW$-pair with
boundary. 

Suppose that we are given two covers $\Sigma$, $\Sigma'$ with boundary
and an equivariant $GL(1)$-representation on each of $\Sigma$ and
$\Sigma'$ with trivialization at the two preimages of a marked point
$z$. If the monodromies around the glued components are the same (when
written relative to the trivializations), then we can either first
$\cW$-nonabelianize and then glue the resulting $\cW$-pairs, or first
glue the two abelian $GL(1)$-representations and then
$\cW$-nonabelianize.

\section{Moduli spaces and spectral coordinates}\label{sec:modulispaces}

\subsection{Moduli spaces}

Given the surface $C$ and a spectral network $\cW$ of Fock-Goncharov or 
Fenchel-Nielsen type, subordinate to a covering $\Sigma$,
we have considered three closely related categories.
Now we consider the corresponding moduli spaces: 
\begin{itemize}
\item let $\cM(C, SL(2), \cW)$ be the moduli space parametrizing flat
$\cW$-framed $SL(2)$-connections over $C$, up to equivalence,
\item let $\cM(\Sigma, GL(1))$ be the moduli space parametrizing
equivariant $GL(1)$-connections over $\Sigma$ up to equivalence,
\item let $\cM(\cW)$ be the moduli space parameterizing $\cW$-pairs
up to equivalence.
\end{itemize}
The constructions of the last few sections lead to a diagram
relating these spaces, as follows:
\begin{center}
\includegraphics{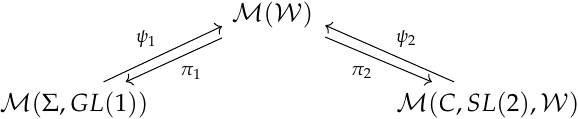}
\end{center}
Here,
\begin{itemize}
\item  $\pi_1$ and $\pi_2$ are the \ti{forgetful} maps
which map a $\cW$-pair to the underlying equivariant $GL(1)$-connection
or $\cW$-framed flat $SL(2)$-connection respectively,
\item $\psi_1$ is the \ti{nonabelianization} map constructed in
\S\ref{sec:nonabelianization}, which extends an equivariant $GL(1)$-connection 
to a $\cW$-pair,
\item $\psi_2$ is the \ti{abelianization} map constructed 
in \S\ref{sec:constructingabelianizations}, which extends a $\cW$-framed flat $SL(2)$-connection 
to a $\cW$-pair.
\end{itemize}
From this description it is evident that $\pi_1 \circ \psi_1$
and $\pi_2 \circ \psi_2$ are the identity maps.
Moreover, the uniqueness we have proven for abelianization
and nonabelianization say that $\pi_1$ and $\pi_2$ are both injective.
It follows that all of the maps are bijections.
In particular, $\Psi = \pi_1 \circ \psi_2$ is a bijection
(in fact diffeomorphism) 
\begin{equation} \label{eq:coord-diffeo}
\Psi: \cM(C, SL(2), \cW) \to \cM(\Sigma, GL(1)).
\end{equation}

Let us verify that the dimensions match. The dimension of the moduli space
$\cM(\Sigma, GL(1))$ is
$2g' - 2g + n$, where $g'$ and $g$ are the genera of $\bar \Sigma$ and $\bar C$
respectively, and $n$ the number of punctures or boundary components.
Using the fact that the number of
branch points is $-2 \chi(C)$ (proven in \S\ref{ssec:defspectralnetwork}), 
we get $\chi(\bar \Sigma) = 4 \chi(\bar C)$, which implies
$2 g' = -6 + 8g + 2n$. Thus we get
\begin{equation}
	\dim \cM(\Sigma, GL(1)) = 6g-6+3n
\end{equation}
indeed matching $\dim \cM(C, SL(2), \cW) = -3 \chi(C)$, as desired.
If we fix the conjugacy classes of the monodromies around the 
punctures on $C$ then we get a similar matching with both
dimensions reduced by $n$.

\subsection{Spectral coordinates}\label{ssec:spectralcoordinates}

Given an equivariant $GL(1)$-connection $\nabla^\ab$ we can construct
some interesting numbers, as follows.
Given any class $\gamma \in H_1(\Sigma', \Z)$ we can
consider the holonomy
\begin{equation}
\cX_\gamma = \Hol_{\gamma} \nabla^\ab \in \C^\times.
\end{equation}
From their definition it immediately follows that
they are multiplicative:
\begin{equation} \label{eq:x-mult}
\cX_\gamma \cX_{\gamma'} = \cX_{\gamma + \gamma'}.
\end{equation}
Moreover, if $\gamma_b$ denotes a small loop around a branch point $b$, then we have
\begin{equation} \label{eq:x-branch}
\cX_{\gamma_b} = -1.
\end{equation}
Finally, the equivariance of $\nabla^\ab$ implies the relation
\begin{equation} \label{eq:x-equivariant}
\cX_{\gamma + \sigma_{*} \gamma} = 1.
\end{equation}
It follows from \eqref{eq:x-mult}, \eqref{eq:x-branch}, \eqref{eq:x-equivariant} that if we fix
a collection $\{\gamma_i\} \subset H_1(\Sigma',\Z)$ which generates
$H_1(\Sigma',\Z) / \IP{\gamma + \sigma_* \gamma}$, then the $\cX_{\gamma_i}$
are enough to determine all of the $\cX_\gamma$.
They thus give a coordinate system on $\cM(\Sigma, GL(1))$.

Alternatively, via the diffeomorphism
$\Psi$ in \eqref{eq:coord-diffeo}, the
$\cX_\gamma$ can be thought of as functions on $\cM(C, SL(2), \cW)$,
and the $\cX_{\gamma_i}$ as above give a coordinate system.
With this in mind, we call the $\cX_{\gamma_i}$ \ti{spectral coordinates}.

Spectral coordinate systems have some nice properties which distinguish them from arbitrary
coordinate systems.  In particular, it was argued in
\cite{Gaiotto2012} that they are \ti{Darboux} coordinates, in the following sense.
The moduli space $\cM(SL(2), C, \cW)$ has a natural holomorphic Poisson structure, described
e.g. in \cite{prr} (generalizing the symplectic structure which one gets when $C$ is 
closed \cite{MR85k:14006}).
The functions $\cX_\gamma$ have simple Poisson brackets with respect
to this structure:
\begin{equation}
\{ \cX_\gamma, \cX_{\gamma'} \} = \IP{\gamma, \gamma'} \cX_{\gamma+\gamma'},
\end{equation}
where $\IP{\cdot,\cdot}$ denotes the intersection pairing on $H_1(\Sigma', \Z)$.
In particular, if we fix our basis $\{\gamma_i\}$ so that $\IP{\cdot,\cdot}$ is of
the block form $$\begin{pmatrix} 0 & 1_{r \times r} & 0 \\ -1_{r \times r} & 0 & 0 \\ 0 & 0 & 0_{m \times m} \end{pmatrix},$$ the spectral coordinates $\cX_{\gamma_i}$ for $i = 1, \dots, r$ are
Darboux coordinates for the symplectic leaves of $\cM(SL(2), C, \cW)$.

So far we have been rather general. In the next two sections we want to 
show how these general considerations recover concretely known Darboux 
coordinate systems on $\cM(C, SL(2), \cW)$,
namely the Fock-Goncharov and complexified Fenchel-Nielsen coordinates.

\subsection{Fock-Goncharov coordinates from Fock-Goncharov networks}

Suppose $\cW$ is a Fock-Goncharov network. Then, given a $\cW$-framed
$SL(2)$-connection $\nabla$, we have explained in \S\ref{ssec:FGcoord}
the construction of the corresponding equivariant $GL(1)$-connection
$\nabla^\ab$. Now we will use this construction to 
compute some of the spectral coordinates.

Consider the cycle $\gamma$ on $\Sigma'$ shown in Figure \ref{fig:gamma-one-quadrilateral}.
We would like to compute the holonomy $\cX_\gamma$ of $\nabla^\ab$ around $\gamma$.
\insfigscaled{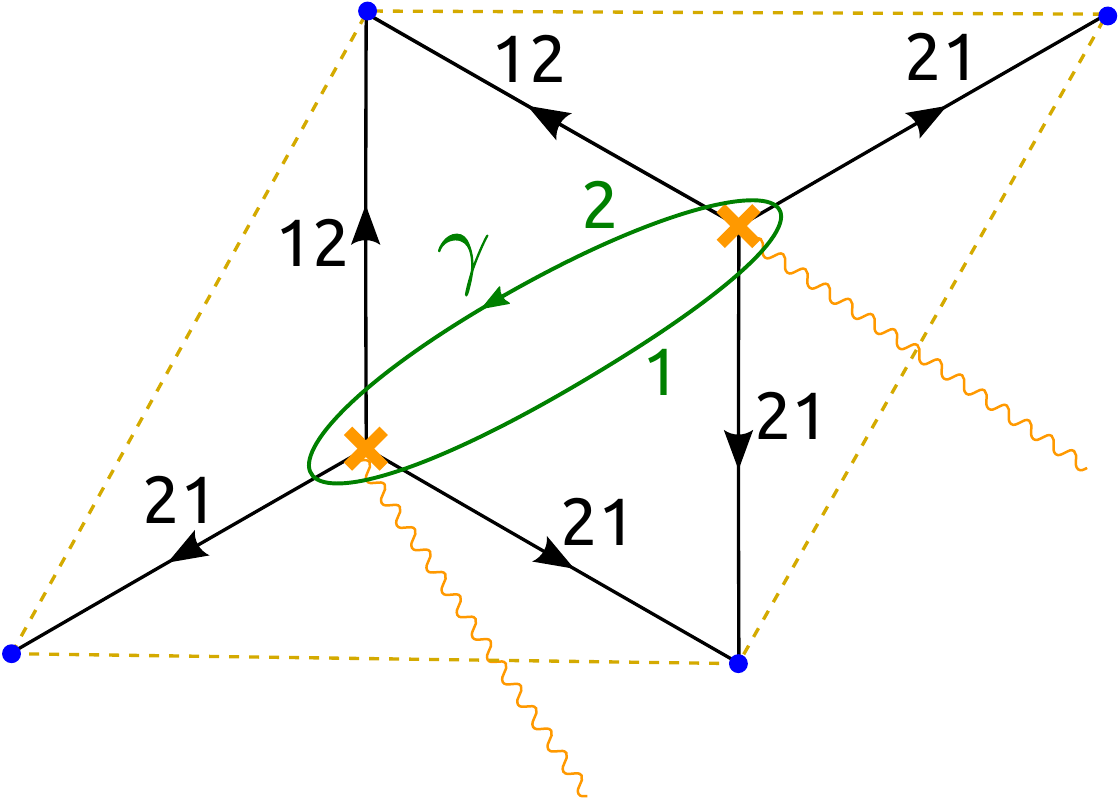}{0.4}{A canonical 1-cycle $\gamma$ on $\Sigma'$, lying 
above the pictured quadrilateral $Q$ in $C$.}
To compute $\cX_\gamma$ it is convenient to split $\gamma$ into two open paths
and compute their parallel transports separately.
Thus, consider a single branch point and an open path $a$ on $\Sigma$,
as shown in Figure \ref{fig:branch-point}.
\insfigscaled{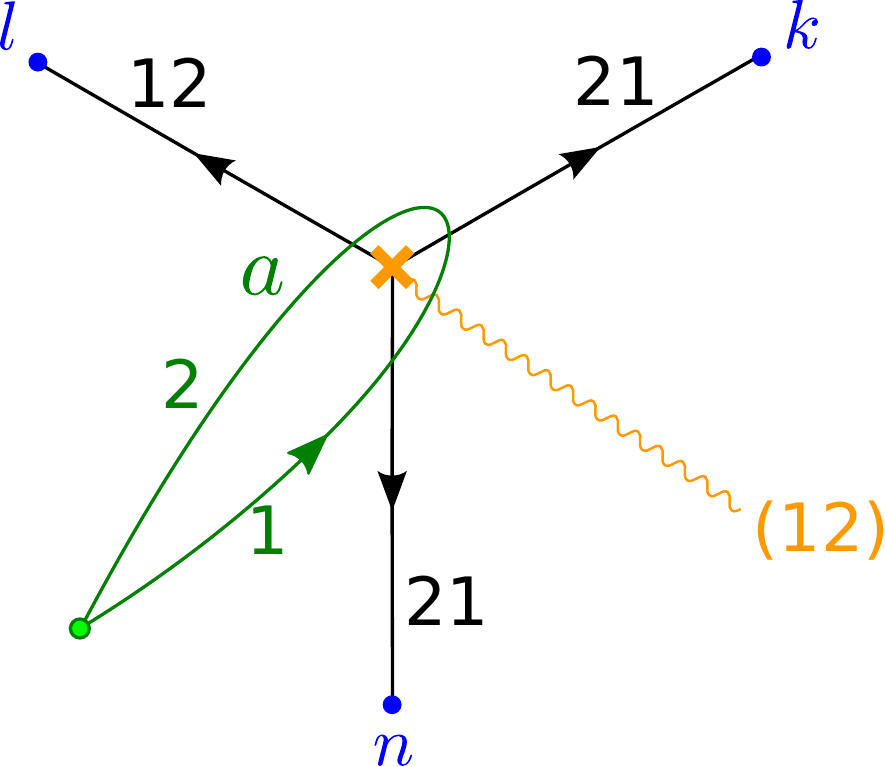}{0.35}{Path $a$ around a branch point. The green labels $1$ and $2$ along the path refer to the sheets at which the path $a$ starts and ends respectively. }
We aim to describe the parallel transport $\cX_a$ of $\nabla^\ab$ along $a$.

Let $s_n$ denote some $\nabla$-flat section of $\ell_n$ 
over the quadrilateral $Q$
(recall that $\ell_n$ denotes the framing at the puncture $z_n$.)
Similarly define $s_k, s_l, s_m$ associated to the 
framings at the other punctures.
Since $a$ is a path beginning on sheet $1$ of $\Sigma$ and ending
on sheet $2$, we have
\begin{equation}\label{eq:cXlambda}
 \cX_a (s_n) = \lambda s_l
\end{equation}
for some scalar $\lambda$.  Our job is to determine $\lambda$.
For this, divide $a$
into three segments $a = a_3 a_2 a_1$,
each crossing one of the walls. Using $\cX_{a}=
\cX_{a_3}\cX_{a_2}\cX_{a_1}$, \eqref{eq:cXlambda} becomes
\begin{equation}\label{eq:cXlambda2}
 \cX_{a_2} \cX_{a_1} (s_n) = \lambda \cX_{a_3^{-1}} s_l.
\end{equation}
The triangular structure
of the $\cS_w$ says that we have
$\cX_{a_1}(s_n) = s_n$, $\cX_{a_3}(s_l) = s_l$,
so that we can simplify
\eqref{eq:cXlambda2} to
\begin{equation}
 \cX_{a_2}  (s_n) = \lambda s_l.
\end{equation}
The triangular structure of the $\cS_w$ also says
$(\cX_{a_2}(s_n) - s_n) \sim s_k$, which means
$(\cX_{a_2}(s_n) - s_n) \wedge s_k = 0$,
i.e.
\begin{equation}
 \lambda = \frac{s_n \wedge s_k}{s_l \wedge s_k}.
\end{equation}

Now consider two neighboring branch points, as in Figure \ref{fig:two-branch-points}.
\insfig{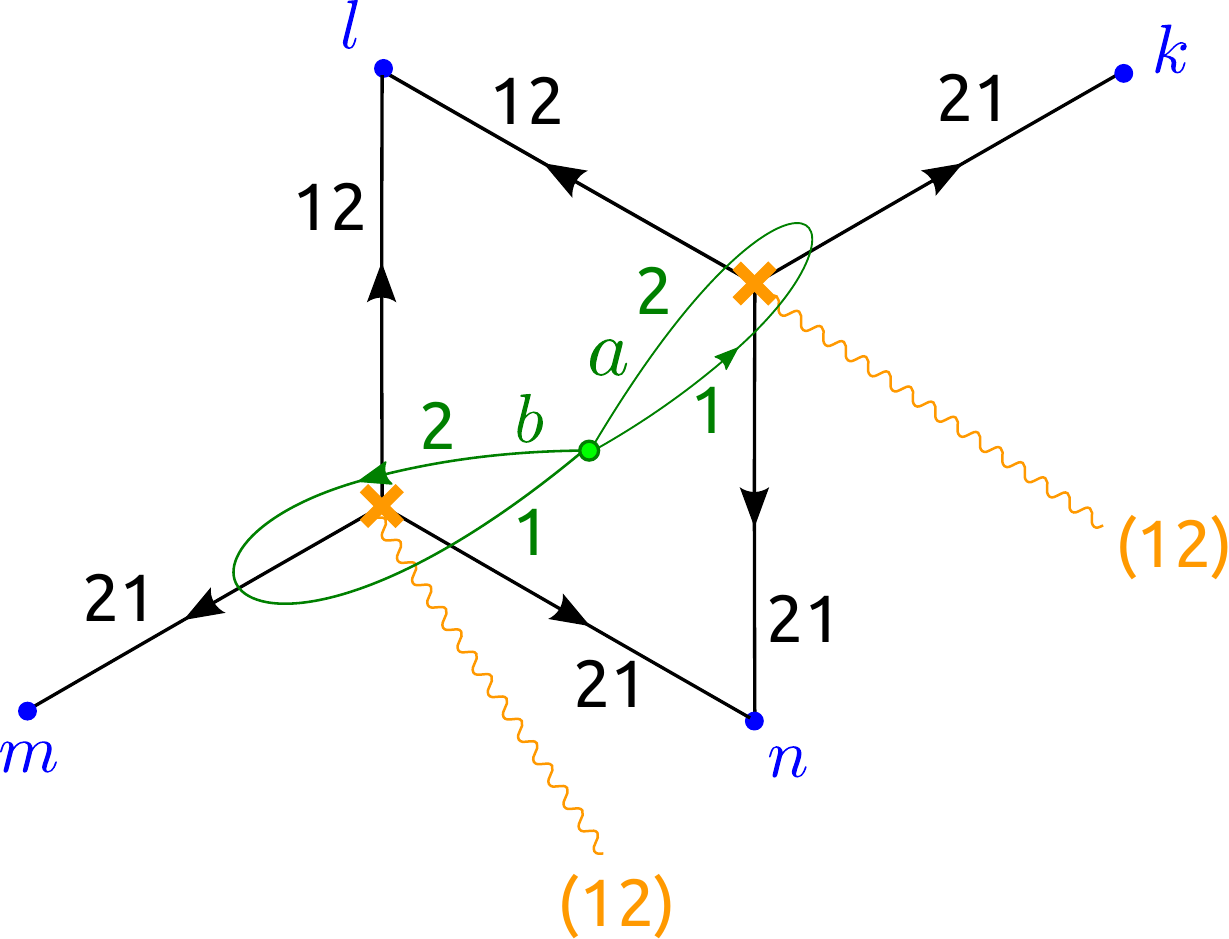}{Cycle $\gamma=a+b$ in the
  quadrilateral $Q$ with vertices $k$, $l$, $m$ and $n$. The green labels $1$ and $2$ refer to sheets of the covering.}
The previous analysis shows that
\begin{equation}
 \cX_a(s_n) = \frac{s_n \wedge s_k}{s_l \wedge s_k} s_l, \qquad \cX_b(s_l) = \frac{s_l \wedge s_m}{s_n \wedge s_m} s_n.
\end{equation}
Combining these gives
\begin{equation}
 \cX_a \cX_b (s_l) = \frac{s_n \wedge s_k}{s_l \wedge s_k} \frac{s_l \wedge s_m}{s_n \wedge s_m} s_l,
\end{equation}
but $a + b = \gamma$, so finally we can write the holonomy along the closed cycle $\gamma$:
\begin{equation} \label{eq:fg-coordinate}
 \cX_{\gamma} = \frac{(s_n \wedge s_k)(s_l \wedge s_m)}{(s_l \wedge s_k)(s_n \wedge s_m)}.
\end{equation}
We emphasize that the expression~\eqref{eq:fg-coordinate} is $SL(2)$
invariant, and can be evaluated at any internal point of the
quadrilateral $Q$.
This is the main result of this section:  it gives a concrete formula for the spectral coordinate
$\cX_\gamma$ attached to the cycle $\gamma$ of Figure \ref{fig:gamma-one-quadrilateral},
in terms of data attached to the original
$\cW$-framed $SL(2)$-connection $\nabla$.

This formula is a familiar one.
Indeed, consider the ideal triangulation $\cT$ which corresponds to the Fock-Goncharov network $\cW$.
The quadrilateral $Q$ contains a unique edge of $\cT$.
On the other hand, for each edge of $\cT$ there is a corresponding
Fock-Goncharov coordinate on the space of $\cW$-framed connections \cite{MR2233852}.
$\cX_\gamma$ given by \eqref{eq:fg-coordinate} is this Fock-Goncharov coordinate multiplied by $-1$.
The sign discrepancy can be eliminated
by shifting $\gamma$ by a loop around a branch point, i.e. replacing $\gamma$ by $\gamma'$ shown in Figure \ref{fig:cycle-sign-shift}.
\insfig{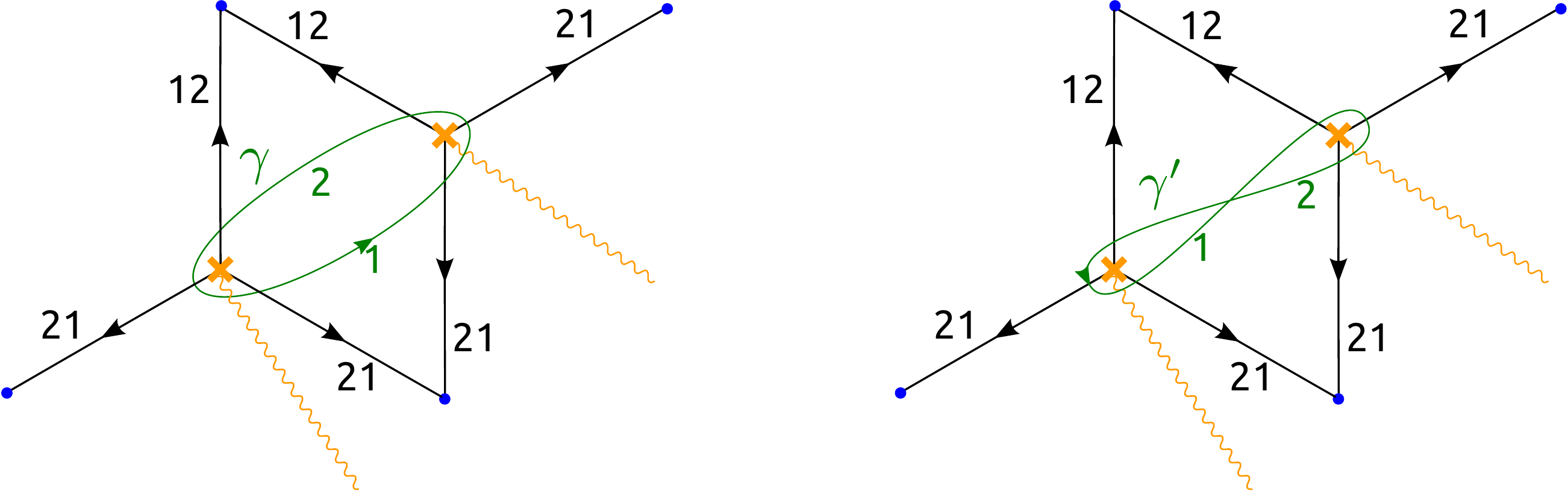}{The homology classes $\gamma$ and $\gamma'$ differ only by the addition of
a loop around one of the branch points on $\Sigma$.  Thus $\cX_{\gamma'} = - \cX_\gamma$.}

Applying this result to each cell in turn, we obtain a collection of cycles 
$\gamma_i \in H_1(\Sigma,\Z)$
such that the spectral coordinates 
$\cX_{\gamma_i}$ are the Fock-Goncharov coordinates attached
to the ideal triangulation $\cT$.
We have thus rederived the result of \cite{Gaiotto:2012db} in the 
special case $K=2$ (see also \cite{Gaiotto:2009hg} where essentially
the same result appeared, albeit phrased in a different
language.)

Knowing the $\cX_{\gamma_i}$ does 
not quite determine the $GL(1)$-connection $\nabla^\ab$.
Indeed, the cycles $\gamma_i$
descend to a $\Q$-basis for $H_1(\Sigma, \Q) / \IP{\gamma + \sigma_* \gamma}$,
but not necessarily to a $\Z$-basis for $H_1(\Sigma, \Z) / \IP{\gamma + \sigma_* \gamma}$.
Thus the full collection of all spectral coordinates $\cX_\gamma$ 
contains slightly more information than the Fock-Goncharov coordinates $\cX_{\gamma_i}$.
This extra information is enough to determine the $SL(2)$-connection $\nabla$ up to equivalence 
(whereas the Fock-Goncharov coordinates would be enough
to determine only its projection to $PSL(2)$.) We will see how this 
works explicitly in some examples in \S\ref{sec:examples}.

\subsection{Fenchel-Nielsen coordinates from Fenchel-Nielsen networks}

Next, suppose $\cW$ is a Fenchel-Nielsen network. Then, given a $\cW$-framed
$SL(2)$-connection $\nabla$, we have explained in \S\ref{ssec:FN-ab}
the construction of the corresponding equivariant $GL(1)$-connection
$\nabla^\ab$. Now we will use this construction to 
compute some of the spectral coordinates.
We will see that they are essentially complexified Fenchel-Nielsen 
coordinates.

One new subtlety here is that, unlike the Fock-Goncharov case just considered,
a Fenchel-Nielsen network does not quite give a
basis for $H_1(\Sigma,\Q) / \IP{\gamma + \sigma_*\gamma}$.
Instead it gives a basis up to some shifts: 
more precisely,
half of the basis elements are canonical (these
will be associated to the Fenchel-Nielsen length coordinates, which are 
likewise canonical) while the other half are ambiguous by certain shifts
(these will be associated to the Fenchel-Nielsen twist coordinates, which are likewise known to suffer from some ambiguities.)

\insfigscaled{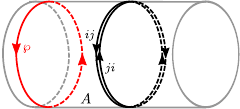}{1.5}{Cycle $\wp$ used to compute the
  spectral length coordinate associated to an annulus $A$. 
  (Only the region near one boundary component of $A$ is shown.)
  }

Consider the annulus $A$ pictured in Figure \ref{fig:FN-coord-path}, and the path $\wp$ going around in the same direction as the nearest walls on the boundary
(call this the $+$ direction).
We choose a labeling of the sheets
of $\Sigma$ over $A$, such that the nearest boundary walls are of type $ij$.
Let $\gamma \in H_1(\Sigma,\Z)$ be the lift of $\wp$ to sheet $j$.
The holonomy $\cX_\gamma$ of $\nabla^\ab$ along $\gamma$
is equal to the eigenvalue of $M_+$ acting on the section
$s_j$.  According to our rules above, $s_j = s_+$.
Thus $\cX_\gamma$ is the monodromy eigenvalue corresponding 
to the section $s_+$.
This monodromy eigenvalue is a complexified version of a
square-root of the exponentiated Fenchel-Nielsen length coordinate
of the connection $\nabla$.
Thus we have found that \ti{square roots of complexified 
exponentiated Fenchel-Nielsen 
length coordinates arise as $GL(1)$ holonomies of $\nabla^\ab$.}

The cycles $\gamma$ obtained in this way do not descend to a basis of 
$H_1(\Sigma,\Q) / \IP{\gamma + \sigma_* \gamma}$. This corresponds
to the fact that the exponentiated complexified 
Fenchel-Nielsen length coordinates do not give a complete coordinate system
on the moduli of flat $PSL(2)$-connections.
To get more cycles we now consider ones which cross the annuli, such as
$\gamma$ pictured in Figure \ref{fig:FN-coord-twist-1}.

\insfigscaled{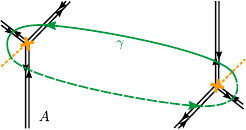}{1.5}{A local region of an annulus $A$,
and a cycle $\gamma \in H_1(\Sigma,\Z)$ which can be used to define the
spectral twist coordinate associated to $A$.}

How should we understand the $GL(1)$ holonomy $\cX_\gamma$?
Rather than trying to compute it directly, let us take a more indirect route:
we want to consider
how $\cX_\gamma$ is transformed under a modification of the connection $\nabla$,
as follows.  Suppose we cut the surface $C$ into two pieces along the annulus $A$.  We will obtain
two surfaces with boundary, say $C_1$ and $C_2$, carrying connections $\nabla_1$ and $\nabla_2$,
as well as an isomorphism $i: E_1 \simeq E_2$ along the boundary, with
$\nabla_1 = i^*\nabla_2$.  We could now glue
$C$ back together along the boundary using the isomorphism $i$ to recover the original $\nabla$.
However, we could instead glue with an isomorphism $i' = i \circ a$, where $a$ is any automorphism of $E_1$
which preserves $\nabla_1$.  Such an automorphism must preserve the
monodromy
eigenspaces:  thus, in terms of the sections $s_1$, $s_2$ of $E$ over $A$, the action of
$a$ can be written as $s_1 \to \lambda s_1$, $s_2 \to \lambda^{-1} s_2$.  Thus we obtain
a 1-parameter family of modified connections $\nabla(\lambda)$.  This operation is sometimes called
the ``twist flow.''

If the original connection $\nabla$ is abelianized by $\nabla^\ab$,
then $\nabla(\lambda)$ is also abelianized by a connection $\nabla^\ab(\lambda)$, constructed by
the $GL(1)$ analogue of the twist flow:  we cut $\Sigma$ apart along the two components of $\pi^{-1}(A)$,
then reglue using the automorphism $s \to \lambda s$ on sheet $1$ and $s \to \lambda^{-1} s$ on
sheet $2$.  From this description, it follows that the twist flow acts simply
on the coordinate $\cX_\gamma$:
\begin{equation} \label{eq:twist-flow}
\cX_\gamma \mapsto \lambda^2 \cX_\gamma.
\end{equation}
One sees this by computing the parallel transport of $\nabla^\ab(\lambda)$
along $\gamma$:  it is identical to that of $\nabla^\ab$, except for two extra factors
of $\lambda$ for the two times the path $\gamma$ crosses $\pi^{-1}(A)$.

The transformation law \eqref{eq:twist-flow} under twist flow
is the characteristic property of a (complexified, exponentiated) ``Fenchel-Nielsen
twist coordinate.''\footnote{We thank Anna Wienhard for emphasizing this 
point of view on the Fenchel-Nielsen coordinates.} We call it \ti{a} twist coordinate rather than 
\ti{the} twist coordinate because of the ambiguity mentioned above. 
To fix this ambiguity one needs
some further choice beyond that of a pants decomposition. What we have seen here is that given a
particular Fenchel-Nielsen network $\cW$ and a particular choice of a cycle $\gamma$ on $\Sigma$ crossing the annulus (as in Figure \ref{fig:FN-coord-twist-1}), 
we get a particular way of fixing the ambiguity, and
\ti{the spectral coordinate
$\cX_\gamma$ is a complexified exponentiated Fenchel-Nielsen twist coordinate.}
Changing the choice of $\gamma$ while keeping the network $\cW$ fixed 
multiplies the twist coordinate $\cX_\gamma$ by a power of the 
exponentiated length coordinate.

The $GL(1)$ holonomies of $\nabla^\ab$ contain
slightly more information than the Fenchel-Nielsen
coordinates. 
Indeed, the spectral coordinates determine the full $SL(2)$-connection $\nabla$ up to equivalence, while
the complexified Fenchel-Nielsen coordinates determine only its
projection to $PSL(2)$. (One reflection of this is the fact that, 
as we have noted, the exponentiated complexified Fenchel-Nielsen length
coordinate is equal to the \ti{square} of the
spectral coordinate $\cX_{\gamma}$.) We will  
see how this works explicitly in the examples in \S\ref{sec:examples}.

\subsection{Mixed coordinates}

Finally we briefly discuss the mixed networks of 
\S\ref{ssec:mixedspectralnetwork} above.
As we discussed in \S\ref{sec:abelianizationmixed}, we can construct 
$\cW$-abelianizations in this case just as for
Fock-Goncharov or Fenchel-Nielsen networks.
We expect that the spectral coordinates 
in this situation will be a kind of hybrid
between Fock-Goncharov and Fenchel-Nielsen coordinate systems 
on the relevant $\cW$-framed moduli spaces;
these hybrid coordinate systems should
contain some coordinates of each kind, along with new 
types of coordinates. We have not worked this story out in detail.

\section{Monodromy representations}\label{sec:examples}

In this section we discuss a few instructive examples of
the nonabelianization construction. In particular, we verify that the
nonabelianization map is injective, and give an 
interpretation of the unipotent elements $\cS_w$. Furthermore,
these examples clarify how the construction depends on whether we consider
connections for the group $SL(2)$ or $PSL(2)$.   

We start this section by outlining how the computations
work in practice. That is, we explain how the somewhat formal description in    
\S\ref{sec:nonabelianization} can be turned into a very concrete
matrix calculation. The first result of this exercise is an explicit
description of the space of equivariant $GL(1)$-representations. The
second result is an explicit description of the nonabelianization
map, which turns equivariant $GL(1)$-representations into 
$\cW$-pairs. As a byproduct we find the monodromy representation of the 
corresponding $SL(2)$-representations in terms of
spectral coordinates.

\subsection{Strategy}

As before we fix a branched covering $\Sigma \to C$
and a spectral network $\cW$. But now we also make a choice of branch
cuts and trivialize $\Sigma$. 

We furthermore choose a set of generators for the paths in
$\cG_C$ such that they either do not cross any walls, or are short
paths connecting the two basepoints attached to a single wall. We
label the former generators as $\wp_n$, and call them \emph{simple}, and
label the latter generators as
$w_k$ (by slight abuse of notation).  
Their lifts to $\Sigma$ form a set of generators for the paths in
$\cG_{\Sigma'}$. We denote them as $\wp^{(ij)}_n$  and
$w_k^{(ii)}$, respectively.\footnote{Our conventions are such that
$\wp^{(ij)}_n$ is a path that starts at sheet $j$ and ends at sheet
$i$.} 

Our first aim is to give an explicit characterization of the space of 
equivariant $GL(1)$-representations. 

An equivariant representation $\rho^\ab$ can be ``trivialized''
 by choosing a basis vector $e_{z^{(i)}}$ in each line 
$\cL_{z^{(i)}}$. We require these choices to be
 compatible with the equivariant structure, i.e. such that over 
each point $z \in \cP_C$ the two elements  $e_{z^{(1)}}$,
$e_{z^{(2)}}$ obey the condition $\mu(e_{z^{(1)}},e_{z^{(2)}}) = 1$.   
In this trivialization $\rho^\ab$  is simply 
given by a collection of $SL(2)$-valued matrices  
$\cD_{\wp_n}$. One for each generating path $\wp_n$, where $\cD_{\wp_n}$
is (off-)diagonal if $\wp_n$ crosses an even (odd) 
number of branch cuts. 

Explicitly, for  a simple generating path $\wp_n^{(ij)}$ the
isomorphism $\rho^\ab(\wp_n^{(ij)})$ is just given by
multiplication with the element $$(\cD_{\wp_n})_{ij},$$ 
whereas for a short generating path $w^{(ii)}_k$ it is simply the identity. When $\wp$ is not
simple, for example $\wp = \wp_3 w_2 \wp_1$, and 
has a lift $\wp^{(ij)} = \wp_3^{(ik)} w_2^{(kk)} \wp_1^{(kj)}$, the
isomorphism $$\rho^\ab(\wp^{(ij)}) = \rho^\ab(\wp^{(ik)}_3)
\rho^\ab(w^{(kk)}_2) \rho^\ab(\wp^{(kj)}_1)$$ corresponds to
multiplication with 
\begin{align*}
                      (\cD_{\wp_3})_{ik} (\cD_{\wp_1})_{kj} 
                     &= (\cD_{\wp_3} \cD_{\wp_1})_{ij}.
\end{align*}
There is an obvious extension to longer paths.

Two different trivializations $\cD_{\wp_n}$  and $\cD'_{\wp_n}$ of
$\rho^\ab$ are
related by an \emph{abelian gauge transformation} $g: \cP_C \to
\C^\times$. If we label the initial point of the generating path $\wp_n$ 
by a number $i(n)$ and its final point by a number $f(n)$, this gauge
transformation acts as 
\begin{align}\label{abeliangaugefreedom}
d_n \to d_n \frac{g_{i(n)}}{g_{f(n)}}
\end{align}
on a diagonal matrix $\cD_{\wp_n} $, with entries $d_n$ and
$d_n^{-1}$, and as
\begin{align}
d_n \to d_n \, g_{i(n)} \, g_{f(n)}
\end{align}
on a strictly off-diagonal matrix $\cD_{\wp_n}$.

The other way around, a collection of matrices $\cD_{\wp_n}$
can only be interpreted as the trivialization of an equivariant
$GL(1)$-representation if there is no monodromy around contractible
cycles on 
$\Sigma'$. We say in that case that the \emph{bulk constraints} are
satisfied.  
Indeed, we can construct an equivariant $GL(1)$-representation $\rho^\ab$ from such a
collection by taking a copy of $\mathbb{C}$ for each line $\cL_z$ and by defining the
isomorphisms $\rho^\ab(\wp)$ through multiplication with
$\cD_{\wp_n}$. 

First trivializing an equivariant $GL(1)$-representation $\rho^\ab$
into a collection of matrices $\cD_{\wp_n}$ obeying the bulk
constraints, and then constructing an equivariant
$GL(1)$-representation from this collection, yields an equivariant
$GL(1)$-representation $\rho'^\ab$ that is equivalent to $\rho^\ab$. 
On the other hand, if we start with a collection of matrices
$\cD_{\wp_n}$ obeying the bulk constraints, turn this into a
representation $\rho^\ab$ and then trivialize again, we find a new
collection of matrices $\cD'_{\wp_n}$ related by an abelian gauge
transformation.

We have thus obtained an explicit characterization of the space of
equivariant $GL(1)$-representations in terms of a collection of
matrices $\cD_{\wp_n}$ obeying the bulk constraints.
Up to abelian gauge transformations the only parameters
that enter the matrices $\cD_{\wp_n}$ are spectral coordinates and 
holonomies around punctures on $\Sigma$. 

Now we are ready to explicitly construct representation $\cW$-pairs.

Given an equivariant $GL(1)$-representation $\rho^\ab$,
\S\ref{sec:nonabelianization} prescribes how to obtain a
$\cW$-pair, and in particular an $SL(2)$-representation $\rho$ (in
``nonabelianization gauge''). Similar as
before, we can trivialize this representation using the two basis vectors
$e_{z^{(1)}}$, $e_{z^{(2)}}$ at every $z \in \cP_C$. It follows that
 $\rho$ can be explicitly described in terms of the  
matrices $\cD_{\wp_n}$ as well, but with some additional unipotent
matrices ``spliced in''. 

More precisely, introduce a normal vector for each wall $w$, a
parameter $S_w \in \C$ and a matrix of the form  
\begin{equation}\label{eq:Smatrix}
\cS_w = \begin{cases} \begin{pmatrix} 1 & S_w \\ 0 & 1 \end{pmatrix} & \text{ for $w$ of type } 21, \\
                      \begin{pmatrix} 1 & 0 \\ S_w & 1 \end{pmatrix} & \text{ for $w$ of type } 12.
        \end{cases}
\end{equation}
Then for any path $\wp \in \cG_C$, the isomorphism $\rho(\wp)$
corresponds to multiplication by matrices $\cS$ and $\cD$. For
example, when $\wp = w_1 \wp_2$ with $w_1$ a short path
crossing the wall labeled by the same symbol and $\wp_2$ not crossing a wall,
$\rho(\wp) $ is given by multiplication with  
\begin{equation}
\cS_{w_1}  \cD_{\wp_2},
\end{equation}  
if the path $\wp$ runs in the direction of the normal vector attached
to the wall $w_1$ that is being crossed. The extension to longer paths
should be clear.

The fact that the matrices $\cD_{\wp}$ obey the bulk constraints on
$\Sigma$ implies that they obey the bulk constraints on $C$ as
well. But if the numbers $S_w$ are chosen arbitrarily, there is
no reason to expect that the branch point constraints are obeyed. The
arguments of \S\ref{sec:nonabelianization} show that the branch point
constraints can be solved uniquely for Fock-Goncharov and
Fenchel-Nielsen networks. This determines the numbers $S_w$ in terms
of the matrices 
$\cD_{\wp}$.  

Thus, for each choice of abelian $GL(1)$-representation
we obtain an explicit characterization of a unique representation
$\cW$-pair in terms of the matrices $\cD_{\wp}$.  In particular, we
can compute the monodromy representation of the resulting
$SL(2)$-representation in terms of purely abelian data.  

Let us illustrate this
with some examples.  

\subsection{Fock-Goncharov networks}

First we describe nonabelianization for two examples of Fock-Goncharov
networks, on the four-punctured sphere and on the one-punctured torus.

\subsubsection{Four-punctured sphere}

Consider the Fock-Goncharov
network on the four-punctured sphere shown in
Figure~\ref{fig:FG-spectral-network-SU2Nf4}.
We make a choice of basepoints $z_k \in \cP_C$ and simple paths
$\wp_m \in \cG_C$ as in Figure~\ref{fig:FG-spectral-network-SU2Nf4-paths-numb}.
Along each simple path $\wp_m$ we choose an $SL(2)$ matrix
$\cD_{\wp_m}$ that is diagonal if the path has an even number of 
intersections with branch cuts and 
strictly off-diagonal if this number is odd. For instance,
\begin{align}
\cD_{\wp_1} = \begin{pmatrix}  d_{1} & 0  \\ 0 &
  d_{1}^{-1}  \end{pmatrix}, \quad \cD_{\wp_6} = \begin{pmatrix}  0 & -d_{6}^{-1}  \\ 
  d_{6} & 0  \end{pmatrix}.
\end{align}
\insfigscaled{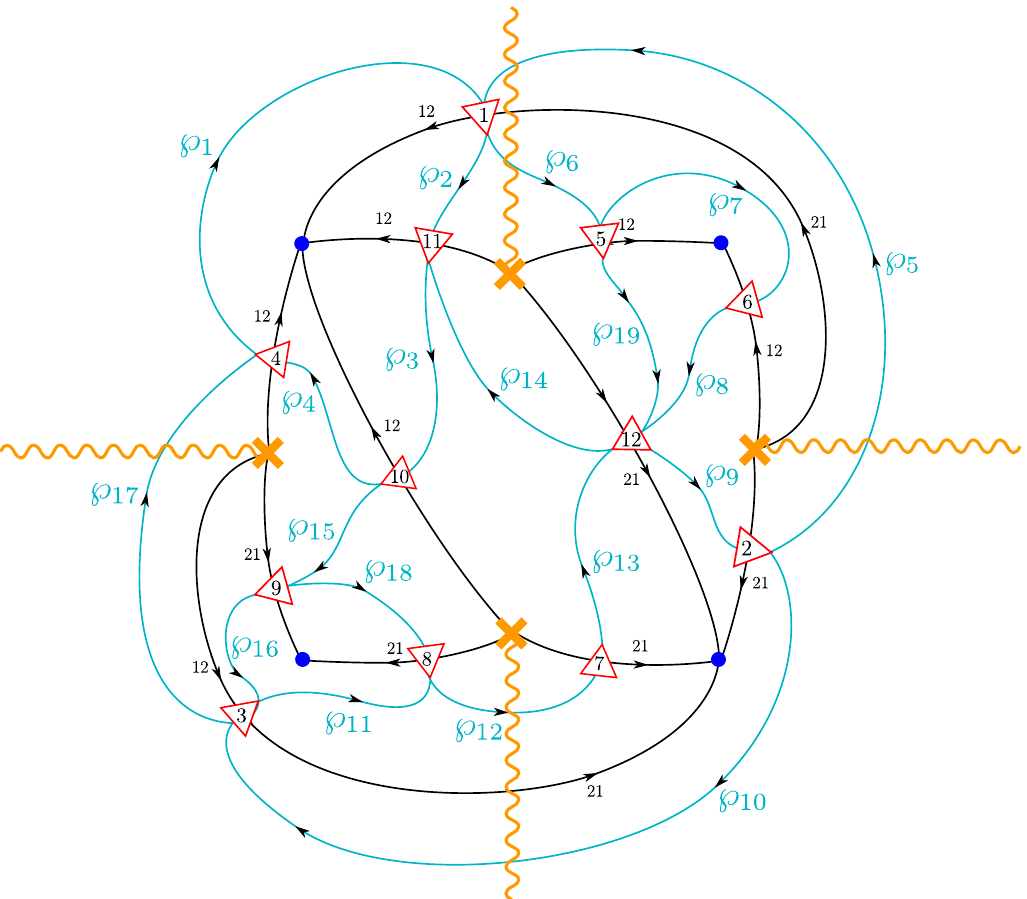}{0.65}{Spectral
  network on the four-punctured sphere with auxiliary data:
  blue lines are simple paths; the arrowhead and tail of the red arrows
  indicate the two base-points at either side of the wall; the direction of the red arrows specifies the normal
 vector at each wall.} 
The collection of matrices $\cD_{\wp_m}$ determines an equivariant abelian
representation $\rho^{\ab}$, provided that the bulk
constraints are satisfied. We choose holonomies
\begin{align}
M^{\mathrm{ab}}_{1,2} = \begin{pmatrix} \cX_{1,2}^{-1} & 0 \\ 0 & \cX_{1,2} \end{pmatrix}
\end{align}
in the clockwise direction around the lifts of the
punctures $z_{1,2}$ to the cover $\Sigma$, and similarly
holonomies
\begin{align}
M^{\mathrm{ab}}_{3,4} = \begin{pmatrix} \cX_{3,4} & 0 \\ 0 & \cX_{3,4}^{-1} \end{pmatrix}
\end{align}
in the clockwise direction around the lifts of the punctures $z_{3,4}$
to the cover.
Then we enforce the bulk equations, for instance
\begin{align}
\cD_{\wp_4} \cD_{\wp_3} \cD_{\wp_2} \cD_{\wp_1} = M^{\mathrm{ab}}_1.
\end{align}
There are four of these bulk equations, plus the equation at infinity
\begin{align}
\cD_{\wp_17} \cD_{\wp_{10}} \cD_{\wp_{5}}^{-1} \cD_{\wp_{1}} = \mathbf{1}.
\end{align}
Thus the nineteen unknowns $d_m$ are constrained by five equations in total, so
that we have fourteen leftover degrees of freedom. This includes the gauge
degrees of freedom, which we discussed around equation~(\ref{abeliangaugefreedom}).

Fixing the abelian gauge freedom, i.e. choosing convenient values for
the twelve $g_k$'s, only leaves two degrees of freedom in the
$d_k$'s. These two unknowns correspond to gauge-invariant  abelian
holonomies around an $A$-cycle and $B$-cycle on the cover $\Sigma$.
Let us call these abelian holonomies $\cX_A$ and $\cX_B$. In this
example we may choose
\begin{align}
\cX_A & = d_2 \,d_{14}^{-1} d_9 \,d_5 \quad \mathrm{and} \quad
\cX_B  = d_5^{-1} d_{10} \,d_{16}\, d_{15} \,d_3 \,d_{14} \,d_8 \,d_{7} \,d_{6}^{-1}.
\end{align}
Up to equivalence the equivariant $GL(1)$-representation
$\rho^{\ab}$ is determined completely in terms of the abelian holonomies $\cX_A, \cX_B,
\cX_1, \ldots, \cX_4$.

Having parametrized all equivariant $GL(1)$-representations on
$\Sigma$, we now pick one of them and apply the non-abelianization
construction.

We fix a normal
vector at each basepoint $z_k$ and write down a matrix $\cS_k$ as in
equation~(\ref{eq:Smatrix}). The matrix $\cS_k$ will determine the
parallel transport across the wall in the direction of the
normal vector.

\insfigscaled{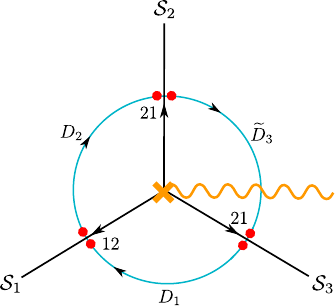}{0.9}{Local configuration around a
branch point in a Fock-Goncharov network.}

Requiring the $SL(2)$-representation $\rho$ to be flat is
equivalent to enforcing branch point constraints. There is
one constraint for each branch point. For any branch point in a
Fock-Goncharov network this constraint is of the form
\begin{align}
\cS_{3,ji} \widetilde{\cD}_{3} \cS_{2,ji}  \cD_{2} \cS_{1,ij}  \cD_{1} = \mathbf{1},
\end{align}
where $\widetilde{\cD}_3$ is a strictly off-diagonal matrix. See
Figure~\ref{fig:around_branchpoint_FG}. The $\cS$-matrices that appear in this equation are
fully determined by this constraint:  their off-diagonal elements are of the form
\begin{align}
S_{k,ij} = (\cD_{k} \cD_{k+2} \cD_{k+1})^{-1}_{ij},
\end{align}
where one of the matrices $\cD$ has a tilde on it. 
(As noted before, this formula has a nice geometric interpretation:  it represents
the $GL(1)$ parallel transport along a  path $a^{(ij)}_k$ that starts at the lift
$z_k^{(j)}$, runs backward along the wall $w^{ij}$ until it
hits a branch point, loops around the branch point, then follows the
wall $w^{ij}$ forward to the lift $z_k^{(i)}$.)

Now that we have determined the $\cS$-matrices, we have an explicit
description of the $SL(2)$-representation $\rho$ in terms of abelian data. In particular, we can
compute the monodromies of $\rho$ around any loop on $C$.

We fix a
presentation of the fundamental group of $C$, i.e. four cycles
$\cP_1, \ldots, \cP_4$
with a common basepoint, chosen such
that $\cP_p$ loops around the $p$th puncture, and with the relation
$\cP_4 \cP_3 \cP_2 \cP_1 = 1$.  Then we compute the $SL(2)$ monodromies $\cM_{\cW}(\cP_p)$
in terms of the abelian data. For instance, if we choose
\begin{equation}
 \cP_1 = w_{11} \wp_2 w_1 \wp_1 w_4 \wp_4 w_{10} \wp_3
\end{equation}
the corresponding monodromy matrix is obtained by splicing in the $\cS$-matrices,
\begin{align}
M_{\cW}(\cP_1) &= \cS_{11} \cD_{\wp_2} \cS_1 \cD_{\wp_1} \cS_4 \cD_{\wp_4} \cS_{10}  \cD_{\wp_3} 
\\ &= \left( \begin{array}{cc} \cX_1^{-1} & 0 \\ g
    & \cX_1 \end{array} \right), \notag
\end{align}
where $g$ is a somewhat complicated expression that changes under
abelian gauge transformations. We can easily check though that (in any gauge)
\begin{align}
M_{\cW}(\cP_4) \cdot M_{\cW}(\cP_3) \cdot M_{\cW}(\cP_2) \cdot
M_{\cW}(\cP_1) = 1,
\end{align}
and hence we have found the monodromy representation of the
$SL(2)$-representation $\rho$ on the four-punctured sphere, in terms of the
abelian holonomies
$\cX_A$, $\cX_B$, $\cX_1, \ldots, \cX_4$.

\subsubsection*{Complexified shear coordinates}

So far we have constructed the  $SL(2)$-representation $\rho$ for a Fock-Goncharov
network in terms of abelian data. Now let us discuss how the Fock-Goncharov coordinates
are explicitly realized in terms of the abelian data.

Recall the definition of the Fock-Goncharov coordinates reviewed in \S\ref{ssec:FGcoord},
which assigns one Fock-Goncharov coordinate to each edge of the ideal triangulation
corresponding to the network.  In our case this triangulation
has six edges; we label the Fock-Goncharov coordinates for these edges
$\cZ_1, \ldots, \cZ_6$ as illustrated in
Figure~\ref{fig:FG-spectral-network-SU2Nf4-shear}. For instance,
\begin{align}\label{eq:FGdef}
\cZ_2 =  \frac{(s_1 \wedge s'_3)(s_2 \wedge s''_3)}{(s_2 \wedge s'_3)(s_1
  \wedge s''_3)}.
\end{align}
Here $s'_3$ and $s''_3$ represent the parallel transport of the
section $s_3$ along two inequivalent paths in the two triangles
on either side of the edge labeled by the coordinate $\cZ_2$.

\insfigscaled{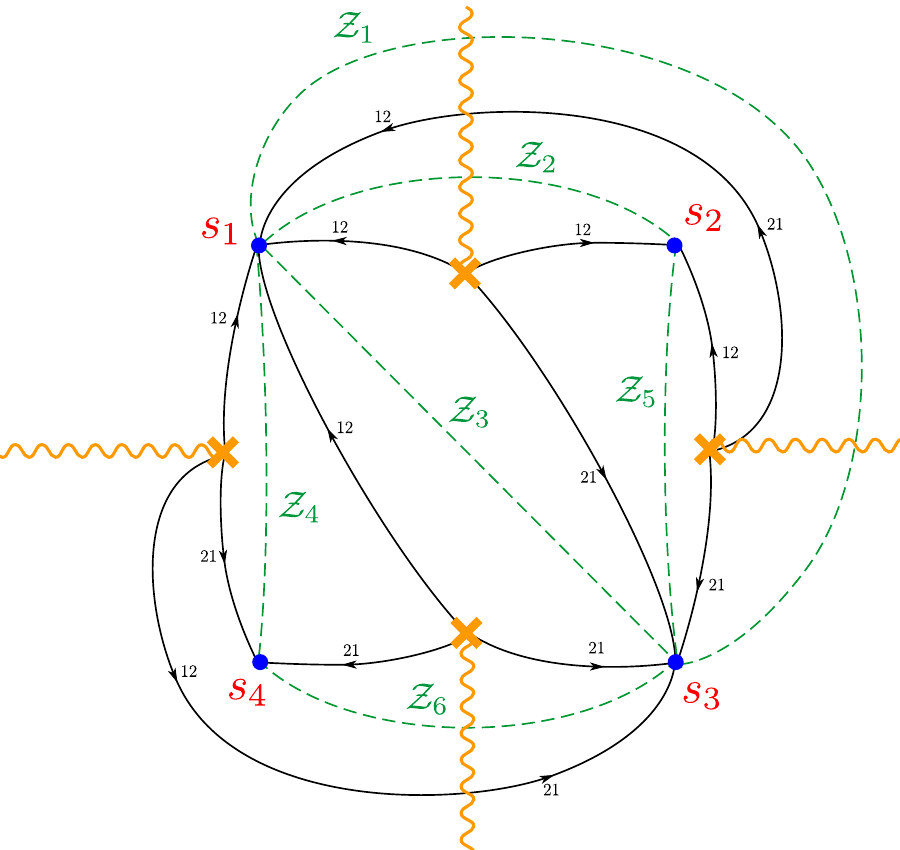}{0.5}{Ideal
  triangulation (in dark green) on the four-punctured sphere, and
  the corresponding shear coordinates
  $\cZ_1, \ldots, \cZ_6$. }

Now, from the discussion of  \S\ref{ssec:FGcoord} we expect that the coordinates
$\cZ_1,\ldots, \cZ_6$ can be expressed in terms of abelian holonomies along odd cycles on $\Sigma$.
Concretely, for example,
$$
\cZ_2 = (\cD_{\wp_9} \cD_{\wp_8} \cD_{\wp_{7}} \cD_{\wp_{19}}^{-1} \cD_{\wp_{14}}^{-1} \cD_{\wp_2} \cD_{\wp_5})_{22}.
$$
In turn these holonomies can be expressed in terms of
$\cX_A, \cX_B, \cX_1, \dots, \cX_4$.  Doing this we obtain
\begin{align}
\cZ_1 = \cX_1 \cX_A \cX_B, \quad
\cZ_2 &=  - \frac{\cX_2}{ \cX_A}, \quad
\cZ_3 =  \frac{\cX_3}{\cX_2 \cX_4 \cX_A \cX_B}, \quad \notag \\
\cZ_4 = - \frac{\cX_1 \cX_4 \cX_A}{ \cX_3}, \quad
\cZ_5 &= -\cX_2 \cX_A, \quad
\cZ_6 = - \frac{ \cX_3 \cX_4}{ \cX_1 \cX_A}. \label{eq:ZrelX4ptsphere}
\end{align}
As a check, the consecutive product of shear
coordinates on the edges ending at the $p$th puncture is equal to
$\cX_p^2$.

If we attempt to invert the relations~(\ref{eq:ZrelX4ptsphere}) to find the abelian
holonomies ($\cX$'s) in terms of the $\cZ$'s, we
find that the $\cX$'s involve square roots in the $\cZ$'s. This is consistent with the fact
that the Fock-Goncharov coordinates $\cZ$ do not completely determine
an  $SL(2)$-representation, only its projection to $PSL(2)$.  In contrast, the $\cX$'s
really do determine the $SL(2)$-representation.

Rather than comparing only the gauge invariant quantities ($\cX$'s and $\cZ$'s) we can also
go a bit further and compare the actual $SL(2)$ parallel transport matrices which we obtain
by nonabelianization (depending on $\cX$'s) to the ones given in the Fock prescription
(depending on $\cZ$'s).  This prescription was given in \cite{Fock97} (see also Appendix A of
\cite{Gaiotto:2009hg}). Let us briefly review it.

We compute the monodromy matrix $M_{\cF}(\cP_p)$ as follows.
Follow the loop $\cP_p$ and write down a matrix
\begin{equation}
E(\cZ_E) = \left( \begin{array}{cc} 0 &1 \\ \cZ_E &
    0 \end{array} \right),
\end{equation}
when we cross an edge $E$ labeled by the shear
coordinate $\cZ_E$. When turning
right (left) inside a triangle, we write down the matrix $V$ ($V^{-1}$) given by
\begin{equation}
V = \left( \begin{array}{cc} -1 & 1 \\ -1 &
    0 \end{array} \right)
\end{equation}
The monodromy matrix $M_{\cF}(\cP_p)$ is then
the ordered product of the matrices $E(\cZ)$ and $V$, multiplied by a
normalization factor which makes the determinant $1$.
For example, in our case
\begin{eqnarray}
M_{\cF}(\cP_1) &= \pm \cX_1^{-1} V E(\cZ_2) V E(\cZ_1) V E(\cZ41) V
E(\cZ_3).
\end{eqnarray}
The normalization factor is only determined up to a sign, so
the resulting $M_{\cF}(\cP_p)$ is really
valued in $PSL(2)$.  We want to compare the $M_{\cF}(\cP_p)$
with the matrices $M_\cW(\cP_p)$ which arise from nonabelianization,
which are  $SL(2)$-valued.  Thus we fix the sign
in $M_{\cF}(\cP_p)$ in such a way that the eigenvalues match.
Denote the resulting $SL(2)$-valued matrices by $M^{\epsilon}_{\cF}(\cP_p)$.
After so doing, we can indeed verify directly that the
monodromy matrices
$M^{\epsilon}_{\cF}(\cP_p)$ and $M_{\cW}(\cP_p)$ are
conjugate.

This comparison also gives a concrete way of seeing
how the reduction from $SL(2)$-representations to $PSL(2)$-representations
works on the level of the spectral coordinates.
Let us call the
$SL(2)$-representation
$\rho$, as before, and its reduction to a $PSL(2)$-representation
$\tilde{\rho}$. The normalization factors of the
Fock matrices $M_{\cF}(\cP_{1,2,3,4})$ are given by $1/\cX_1$, $\cX_2$, $
\cX_1/(\cX_2 \cX_4)$ and $ \cX_4$, respectively, and are dependent on
the choice of the lift of $\tilde{\rho}$ to $\rho$. In
contrast, the shear coordinates $\cZ_E$ should
only depend on $\tilde{\rho}$. This suggests that any two lifts of $\tilde{\rho}$
to $SL(2)$ are related to one another by a transformation
\begin{align}
\cX_1 &\mapsto \epsilon_1 \, \cX_1\,, \notag \\
\cX_2 &\mapsto \epsilon_2 \, \cX_2\,,  \\
\cX_4 &\mapsto \epsilon_3 \, \cX_4\,, \notag
\end{align}
for a choice of $\epsilon_1,\epsilon_2,\epsilon_3 \in \{\pm1\}$. Indeed, from the
relations~(\ref{eq:ZrelX4ptsphere}) it is easy to see that this
operation can be uniquely completed to a transformation
on all the $GL(1)$ holonomies $\cX_A,\, \cX_B$, $\cX_1, \ldots, \cX_4$, in such a way that the
shear coordinates $\cZ_E$ are invariant. 

We call two equivariant $GL(1)$-representations $PSL(2)$-equivalent if they
are related by a transformation of this kind. Then we find a
1-1 correspondence between 
equivariant $GL(1)$-representations, up to $PSL(2)$-equivalence, 
and $PSL(2)$-representations.

\subsubsection{Once-punctured torus}

Next, we work out the nonabelianization map for a Fock-Goncharov network
on the once-punctured torus. Figure~\ref{fig:FGnetwork1torus}
illustrates the chosen spectral network with auxiliary data such as
basepoints and simple paths.

\insfigscaled{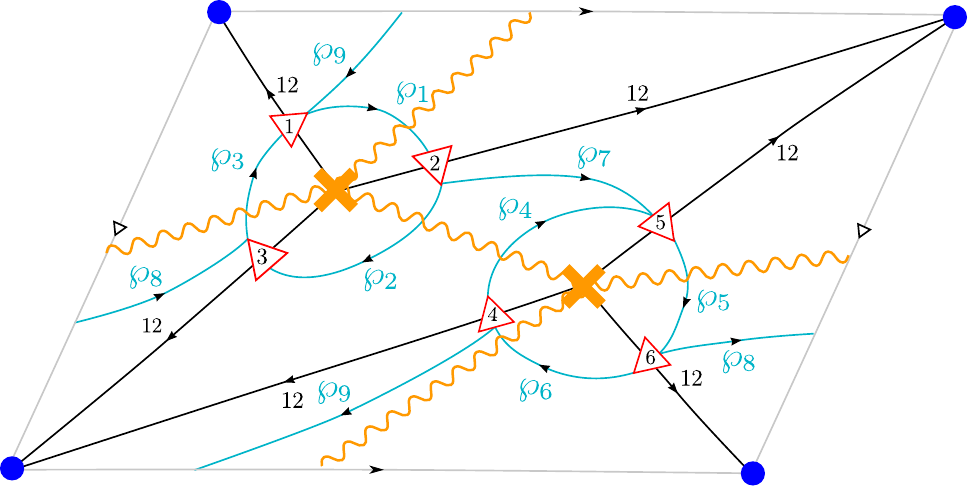}{0.65}{Spectral
  network on the once-punctured torus with auxiliary data.}

We choose an $SL(2)$ matrix $\cD_{\wp}$ for each simple path
$\wp$. This matrix 
is diagonal if the path has an even number of
intersections with branch cuts and 
strictly off-diagonal if this number is odd. Altogether these matrices contain 9 degrees of
freedom $d_{\wp}$. 

The collection of matrices $\cD_{\wp}$ determines an equivariant $GL(1)$-representation $\rho^{\ab}$ iff the matrices obey the bulk
constraints.
In this example there are three bulk equations:
\begin{align}
M^{\ab}_a  &= \cD_4^{-1} \cD_7 \cD_1 \cD_9, \\
M^{\ab}_b &= \cD_5 \cD_7 \cD_2^{-1} \cD_8, \notag \\
M^{\ab}_m &= \cD_5 \cD_7 \cD_1 \cD_9 \cD_6 \cD_8^{-1}  \cD_2 \cD_7^{-1}  \cD_4 \cD_9^{-1} 
\cD_3 \cD_8, \notag
\end{align}
where we define
\begin{align}
M^{\mathrm{ab}}_{a,b,m} &= \begin{pmatrix} \cX_{a,b,m} & 0 \\ 0 & \cX_{a,b,m}^{-1} \end{pmatrix}.
\end{align}
Each equation fixes the abelian holonomy around the lifts of a 1-cycle
on the once-punctured torus to the covering. 
This completely fixes the equivariant $GL(1)$-representations on the once-punctured torus in terms of the abelian holonomies $\cX_a$, $\cX_b$ and
$\cX_m$.

The next step is to apply the nonabelianization construction. This
works similar to the previous example. For
any abelian $GL(1)$-representation, we can determine the
$\cS$-matrices using the branch point equations and find the resulting
$SL(2)$-representation $\rho$ in terms of abelian data. 

Given a presentation of the fundamental group of the once-punctured
torus, we can then compute the corresponding $SL(2)$ monodromy matrices
$M_{\cW}(\cP_{a,b,m})$.  For instance, if we fix three cycles  $\cP_a$, $\cP_b$ and $\cP_m$ with
the relation 
\begin{align}
\cP_b^{-1} \cP_a^{-1} \cP_b \cP_a = \cP_m,
\end{align}
we indeed verify that 
\begin{align}
M_{\cW}(\cP_b)^{-1} M_{\cW}(\cP_{a})^{-1}  M_{\cW}(\cP_b) M_{\cW}(\cP_a) = M_{\cW}(\cP_m).
\end{align}
For example,
\begin{align}
M_{\cW}(\cP_a) &= \cS_2 \cD_1 \cD_9  \cS_{4}^{-1} \cD_{4}^{-1}  \cD_{7} 
\\ &= \left( \begin{array}{cc} \cX_a^{-1} & \cX_a \cX_b^2 \cX_m^{-1} g^{-1}
    \\ g \cX_a^{-1}
    & \cX_a (\cX_b^2 + \cX_m) \cX_m^{-1} \end{array} \right), \notag
\end{align}
where $g$ is a gauge-dependent quantity.

As in the previous example, we might compare these monodromy
matrices with those which arise from the Fock prescription. There are  three
shear coordinates $\cZ_{1,2,3}$ on the
once-punctured torus. Expressed in terms of the abelian
holonomies $\cX_{a,b,m}$ they read
\begin{align}
\cZ_1 =  \cX_a^2 \cX_m^{-1}, \quad \cZ_2 =  \cX_b^2 \cX_m^{-1},  \quad
\cZ_3 =   \cX_m \cX_a^{-2} \cX_b^{-2}.
\end{align}
Notice that even though $\cX_m$ can be expressed as a ratio of shear
coordinates, $\cX_a$ and $\cX_b$ cannot.

We compute the Fock matrices $M^{s}_{\cF}(\cP_{a,b,m})$ as
\begin{align*}
M_{\cF}^{+}(\cP_a) &= \cX_a V^{-1}  E(\cZ_2) V E (\cZ_3), \\
M_{\cF}^{+}(\cP_b) &= \cX_b V  E(\cZ_1) V^{-1} E (\cZ_2), \\
M_{\cF}^{+}(\cP_b) &= M_{\cF}^{+}(\cP_b)^{-1}  M_{\cF}^{+}(\cP_a)^{-1} M_{\cF}^{+}(\cP_b)
M_{\cF}^{+}(\cP_a) ,
\end{align*}
and compare them to the spectral matrices
$M_{\cW}(\cP_{a,b,m})$.  As before, we find that the two sets of matrices
are conjugate.

We can reduce from $SL(2)$ to $PSL(2)$-representations by considering
transformations
\begin{align}
\cX_a &\mapsto \epsilon_a \cX_a, \\
\cX_b &\mapsto \epsilon_b \cX_b, \notag \\
\cX_m &\mapsto \cX_m, \notag
\end{align}
for any choice of $\epsilon_a, \epsilon_b \in \{\pm 1\}$. Clearly this leaves the shear
coordinates $\cZ_E$ invariant (and hence leaves the $PSL(2)$-representation
invariant), while changing the $SL(2)$-representation $\rho$.

\subsection{Fenchel-Nielsen examples} \label{sec:fn-examples}

Next we study the nonabelianization for Fenchel-Nielsen
networks. We apply the construction to
molecule I and II on the three-holed sphere, and also to Fenchel-Nielsen
networks on the four-holed sphere and the genus two surface.

Let us note that, in principle, it would be enough to illustrate the
procedure for three-holed spheres. We can obtain $\cW$-pairs on any
other surface by the gluing construction.

\subsubsection{Pair of pants: molecule I}

We start with the spectral network ``molecule I'' on the three-holed
sphere. As always, we first specify the auxiliary data. Here we fix three
basepoints for each double wall and a basis of simple paths
$\wp_n$. Additionally, we  fix a basepoint $z_{1,2,3}$ on 
each boundary component and assign a sheet ordering to each
orientation of that boundary component. All choices
are illustrated in Figure~\ref{fig:molecule-II-Swall-base-paths}.

\insfigscaled{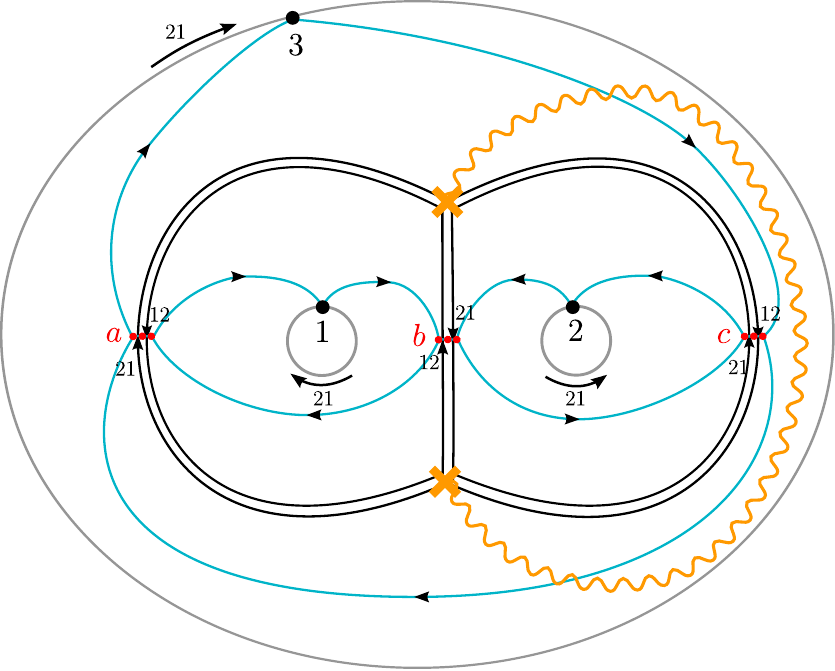}{0.5}{Molecule I
  with a marked point on each boundary (black dot), three
  basepoints for every double wall (red dots) and a collection of simple paths
  (blue lines). The boundary components are labeled as $1$, $2$ and $3$ (in
  black) and the double walls as $a$, $b$ and $c$ (in red). The arrow
  on each boundary component indicates the sheet ordering in that orientation.}

As before, we choose an $SL(2)$ matrix $\cD_{\wp}$ for each simple
path $\wp$.  This matrix is diagonal if the path has an even number of
intersections with branch cuts and 
strictly off-diagonal if this number is odd. In this example these matrices
have 9 degrees of freedom $d_{\wp}$.
The collection of matrices
$\cD_{\wp}$ forms an equivariant $GL(1)$-representation $\rho^{\rm ab}$ iff the
matrices obey the bulk constraints.

We choose a closed path $\gamma_{1,2,3}$ along each boundary
component (in the same orientation as that of the boundary component), and
enforce the $GL(1)$ holonomy along the lift $\gamma^{(ii)}_{1,2,3}$
to be equal to the $i$th diagonal entry of the matrix
\begin{align}
M^{\mathrm{ab}}_{1,2,3} = \begin{pmatrix} \cX_{1,2,3} & 0 \\ 0 & \cX_{1,2,3}^{-1} \end{pmatrix},
\end{align}
where $\cX_{1,2,3} \neq 0, \pm 1$.\footnote{To restrict to irreducible
representations we furthermore impose that $\cX_1 \neq \cX_ 2 \cX_3$,
$\cX_2 \neq \cX_ 3 \cX_1$, $\cX_3 \neq \cX_ 1 \cX_2$ and $\cX_1 \cX_ 2
\cX_3 \neq 1$. We find that these constraints translate to the
constraints that we imposed when abelianizing Fenchel-Nielsen networks
at the level of the image of the non-abelianization map
 (see for instance Proposition 3.2 in \cite{Kabaya}).}
This gives three constraints on the nine unknowns
$d_{\wp}$. 

Abelian gauge transformations at the basepoints
reduce the number of unknowns to three.
The three remaining unknowns can be expressed in terms of parallel
transport coefficients
$\cX_{12,23,31}$ along open paths that run from one boundary component to
another.
We can choose to additionally divide out by abelian gauge transformations at the
 marked points $z_{1,2,3}$, in which case the parallel transport
 coefficients drop out and the number of
unknowns reduces to zero. 

\insfigscaled{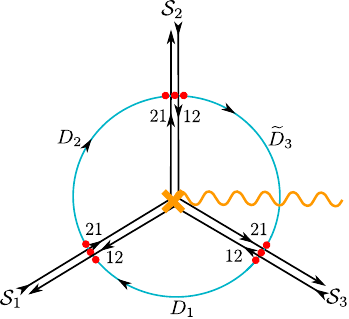}{0.9}{Local configuration around a
branch point in a Fenchel-Nielsen network.}

This gives a parametrization of all equivariant $GL(1)$-representations with
boundary $\rho_{\mathrm{ab}}$ on the cover $\Sigma$.  We now pick one and apply the
non-abelianization construction. 

We construct an $SL(2)$-representation from this 
$GL(1)$-representation by splicing in $\cS$-matrices
whenever the path $\wp$ crosses a wall. Requiring the resulting $SL(2)$-representation to be flat is equivalent to enforcing
branch point constraints. There is one constraint for each branch
point (see Figure~\ref{fig:around_branchpoint_FN}). 
The $\cS$-matrices are fully determined once we impose the
branch point equations.
 For instance, consider the
$\cS$-matrix attached to the double wall $(w_a^{21},w_a^{12})$ on the left in
Figure~\ref{fig:molecule-II-Swall-base-paths}. When crossing the
double wall from the left to the right, this matrix can be written in the form
\begin{align}
\cS_{a} &= \begin{pmatrix} 1 & 0  \\  S(w_a^{12}) &
    1 \end{pmatrix} \begin{pmatrix} 1 &S(w_a^{21}) \\
     0 & 1 \end{pmatrix},
\end{align}
with off-diagonal components
\begin{align}
S(w_a^{12}) & = \cX(a_a^{12})  \left( \frac{1 + \cX_b^2}{1 - \cX_a^2
  \cX_b^2} \right) \\
S(w_a^{21}) & = \cX(a_a^{21}) \left( \frac{1 + \cX_c^2}{1 - \cX_a^2
    \cX_c^2} \right),
\end{align}
where we redefined $\cX_1 = \cX_{a} \cX_{b}$, $\cX_2 =
\cX_{b} \cX_{c}$ and $\cX_3 = \cX_{a} \cX_{c}$. The novel factor in these
$\cS$-matrices is an infinite sum of
monomials in the coordinates $\cX_{a,b,c}$.

\insfigscaled{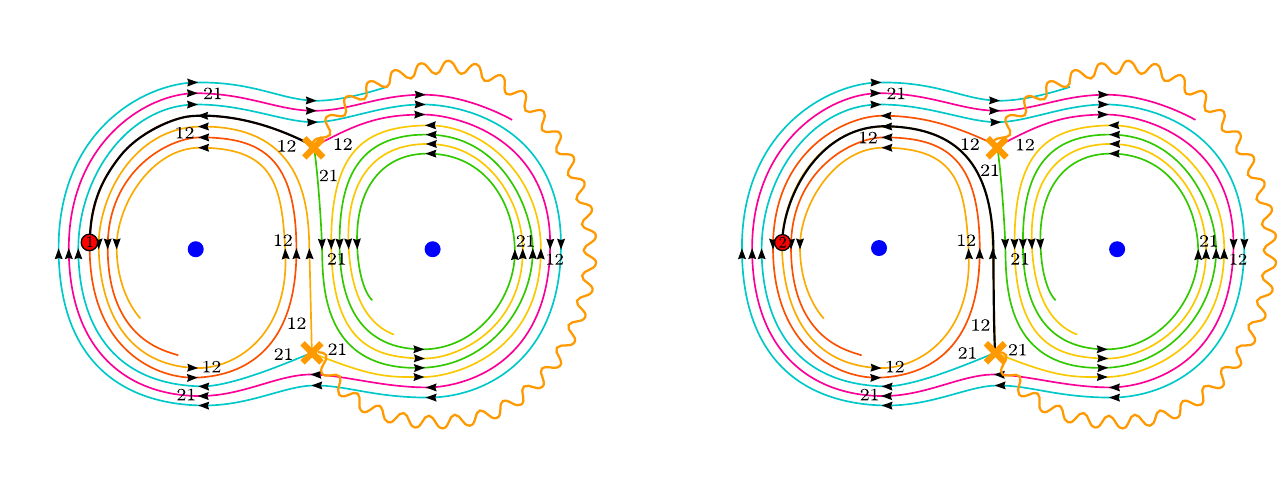}{0.7}{Left and right illustrate the
  interpretation of the terms $\cX(a_a^{12})$ and  $\cX_b^2 \cX(a_a^{12})
  $ in the $\cS$-matrices, respectively, in the limiting spectral
  network to molecule II. }

The interpretation for each
monomial in the sum is best understood from the limiting spectral network
illustrated in Figure~\ref{fig:molecule-II-british-aux}.  The term
$$ \cX(a^{12}_a)$$ computes the abelian parallel transport along a
path on $\Sigma$ that starts at the red dot labeled by 1 (lifted to
the second sheet of the cover), follows
the (black part of the) orange wall to the upper branch point, circles around it, and runs
back
 to the same red dot (lifted to the first sheet of the cover). The term
$$ \cX_b^2~\cX(a^{12}_a)$$
computes the abelian parallel transport along a path that starts at the
red dot labeled by 2 (lifted to the second sheet of the cover), follows
the (black part of the) yellow
wall all the way to the lower
branch point, circles around it, and runs back to
the same red dot (lifted to the first sheet of the cover). All other monomials
in the $\cS$-wall matrices can be
explained similarly.

The $\cS$-matrices attached to the two other double walls in the
spectral network are obtained from $\cS_a$ using the cyclic
symmetry of the spectral network.

\subsubsection*{Length-twist coordinates}

The resulting $SL(2)$-representation $\rho$ has a diagonal monodromy
matrix
\begin{equation}
 \rho(\gamma_{1,2,3}) = \begin{pmatrix} \cX_{1,2,3} &0 \\
     0 & \cX_{1,2,3}^{-1} \end{pmatrix}
\end{equation}
along the closed paths $\gamma_{1,2,3}$. The
$GL(1)$ holonomies $\cX_{1,2,3}$ are each equal to a square-root of an
exponentiated complex Fenchel-Nielsen length parameter. Thus we have
found that these length parameters occur as an example of spectral
parameters.

If we fix a trivialization of the equivariant $GL(1)$-representation
at the marked points, the resulting $SL(2)$-representation $\rho$ not only depends on the abelian
holonomies $\cX_{1,2,3}$, but also on the parallel transports
$\cX_{12,23,31}$ along open paths that begin and end on different
boundary components. If we divide out by
diagonal gauge transformations at the marked points the
dependence on the parallel transport coefficients drops out. 

We can furthermore reduce to $PSL(2)$-representations by considering
transformations that map
\begin{align}
(\cX_1,\cX_2,\cX_3) \mapsto􏰀 (\epsilon_1 \cX_1, \epsilon_2 \cX_2, \epsilon_3 \cX_3),
\end{align}
for any choice of $(\epsilon_1,\epsilon_2,\epsilon_3) \in (\Z/2\Z)^3$
satisfying $\epsilon_1 \epsilon_2 \epsilon_3 = 1$ (see also \cite{Kabaya}).
Let us call two equivariant $GL(1)$-representations $PSL(2)$-equivalent if they
are related by such a transformation. Then we find a
1-1 correspondence between 
equivariant $GL(1)$-representations, up to $PSL(2)$ 
equivalence, and $PSL(2)$-representations.

\subsubsection{Pair of pants: molecule II}

\insfigscaled{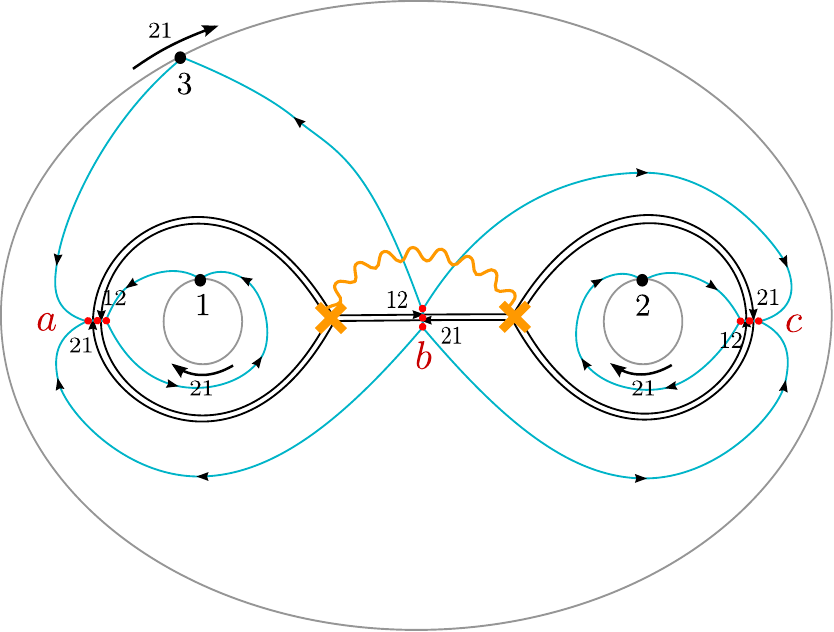}{0.5}{Molecule II
  with a marked point on each boundary (black dot), three
  basepoints for every double wall (red dots) and a collection of simple paths
  (blue lines). The boundary components are labeled as $1$, $2$ and $3$ (in
  black) and the double walls as $a$, $b$ and $c$ (in red). The arrow
  on each boundary component indicates the sheet ordering in that
  orientation}
The description of the nonabelianization map for molecule II is similar to that for
molecule I. We can for instance determine it using the data in
Figure~\ref{fig:molecule-I-Swall-base-paths}. The off-diagonal
components of the $\cS$-matrices are functions of the $GL(1)$ holonomies
$\cX_{1,2,3}$ along closed
paths $\gamma_{1,2,3}$ homotopic to the respective boundary
components. If we redefine these parameters through the
equations $\cX_1 = \cX_a$, $\cX_2 =  \cX_c$ and $\cX_3 =  1/(\cX_a
\cX_b^2 \cX_c)$, the off-diagonal components are a series in the
parameters $\cX_{a,b,c}$. Crossing the
$21$-wall first and the $12$-wall last, we have
\begin{align}\label{eq:S-matrix-3puncturedsphere}
\cS_{a,b,c} &=  \begin{pmatrix} 1 &0 \\
      S(w_{a,b,c}^{12}) & 1 \end{pmatrix} \begin{pmatrix} 1 & S(w_{a,b,c}^{21}) \\ 0 &
    1 \end{pmatrix} ,
\end{align}
where
\begin{align}\label{eq:S-matrix-moleculeII}
\begin{array}{ll}
 S(w_a^{12}) = \cX(a^{12}_a)   \frac{1}{\left(
  {1-\cX_a^2}\right)}  &
 S(w_a^{21})  = \cX(a^{21}_a)  \frac{\left(1+\cX_b^2\right)
   \left(1+\cX_c^2 \cX_b^2\right)}{1-\cX_a^2 \cX_b^4
   \cX_c^2 }  \\
S(w_b^{12})  = \cX(a^{12}_b)  \frac{1+\cX_a^2
    \left(1+\left(1+\cX_c^2\right)
      \cX_b^2\right)}{1-\cX_a^2  \cX_b^4
    \cX_c^2}  &
S(w_b^{21})  = \cX(a^{21}_b)  \frac{1+\cX_c^2
    \left(1+\left(1+\cX_a^2\right)
      \cX_b^2\right)}{1-\cX_a^2 \cX_b^4
    \cX_c^2} \\
S(w_c^{12})  = \cX(a^{12}_c)  \frac{1}{ \left(
  1-\cX_c^2\right)}  &
S(w_c^{21}) = \cX(a^{21}_c)  \frac{\left(1+\cX_b^2\right)
  \left(1+\cX_a^2 \cX_b^2\right)}{\left( 1-\cX_a^2 \cX_b^4
  \cX_c^2 \right)} . \\
 \end{array}
 \end{align}
As the reader can verify, all monomials in the $\cS$-matrices can be interpreted as
abelian holonomies along auxiliary paths in the limiting network
(which is illustrated in
Figure~\ref{fig:molecule-I-british-punct}). 

\insfigscaled{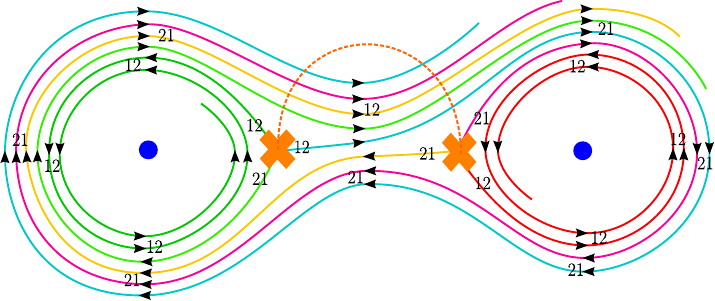}{0.7}{Limiting spectral
  network to molecule I. }

\subsubsection{Four-holed sphere}

Given a pair of pants decomposition of the four-holed sphere, we can
build a Fenchel-Nielsen type spectral network choosing either molecule
I or II on each pair of pants. Let us for instance determine the
nonabelianization map for the spectral network illustrated in
Figure~\ref{fig:four-holed-sphere-Swall}, built by gluing two pairs
of pants containing molecules of type II.

\insfigscaled{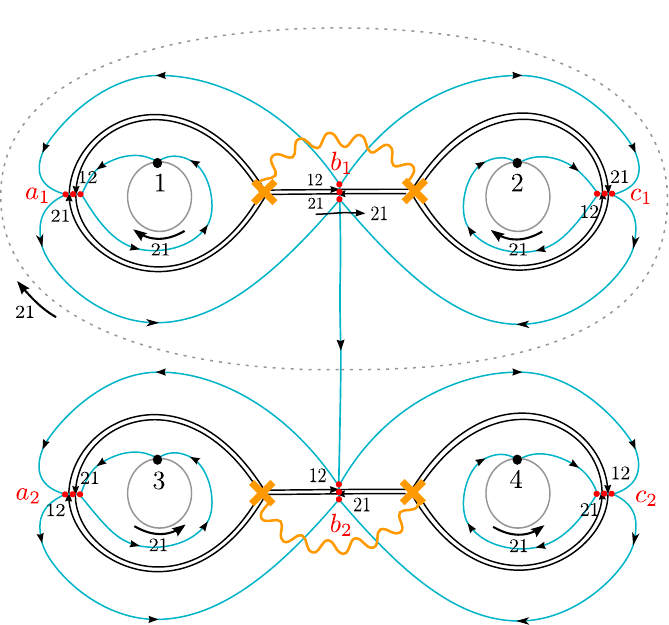}{0.7}{Fenchel-Nielsen spectral
  network on the four-holed sphere
  with a basepoint on each boundary (black dot), three
  basepoints for every double wall (red dots) and a collection of simple paths
  (blue lines). Additionally, an arrow denotes the sheet ordering in
  that direction on every boundary
  component and on the curve that defines the pants decomposition.}

We fix basepoints and choose simple paths $\wp_{1 \ldots 13}$ as illustrated in
Figure~\ref{fig:four-holed-sphere-Swall}. The matrices $\cD_{\wp}$ thus
contain 13 variables in total. The bulk constraints consist of
boundary conditions at the boundaries and at infinity, and reduce the
number of variables to 8. After fixing  6 abelian gauge degrees of
freedom (one for each double wall), we are left with two gauge-invariant
$GL(1)$ holonomies. Let
us call these holonomies $\cX_A$ and
$\cX_B$.

The cover $\Sigma$ of the 4-holed sphere is an 8-holed torus. The
$GL(1)$ holonomies $\cX_{1,2,3,4}^{\pm 1}$ around the boundary components
are part of the boundary conditions that we
have just imposed. The
  coordinates $\cX_A$ and $\cX_B$ are abelian
holonomies along a choice of 1-cycles $A$  and $B$ on the cover $\Sigma$ that intersect
each other once. To be
explicit, let the 1-cycle $A$ be a lift to sheet 1 of $\Sigma$ of the
1-cycle $\alpha$ that defines the pants
decomposition (see
Figure~\ref{fig:four-holed-sphere-cycles}). We take $B$ as
the lift of a 1-cycle $\beta$ in
Figure~\ref{fig:four-holed-sphere-cycles}, lifted and oriented in such
a way that $A \# B =1$.\footnote{We furthermore impose the
  constraints that $\cX_{1,2,3,4,A} \neq 0, \pm 1$, and that when
  restricted to each pair of pants $\cX_{i_1} \neq \cX_ {i_2} \cX_{i_3}$,
$\cX_{i_2} \neq \cX_ {i_3} \cX_{i_1}$, $\cX_{i_3} \neq \cX_ {i_1}
\cX_{i_2}$ and $\cX_{i_1} \cX_ {i_2}
\cX_{i_3} \neq 1$. These constraint are equivalent to the constraints that
we imposed when abelianizing with Fenchel-Nielsen networks \cite{Kabaya}.}

\insfigscaled{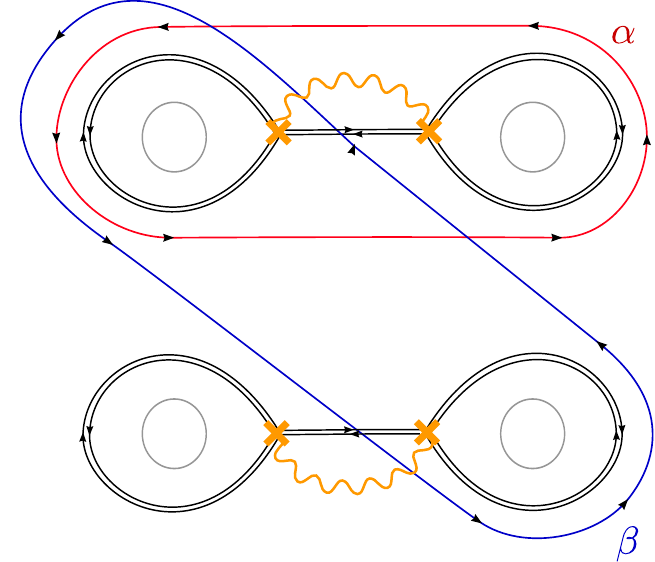}{0.5}{Choice of 1-cycles $\alpha$
  and $\beta$ on the four-holed sphere.}

Having parametrized all equivariant $GL(1)$-representations
on the cover $\Sigma$, we construct the corresponding nonabelian representations
 on the four-holed sphere by splicing in
$\cS$-matrices. These $\cS$-matrices are completely
determined by the non-abelian branch point equations. As can be
expected from the gluing construction, we find that the
$\cS$-matrices are similar to those for the molecule II network on the
three-holed sphere.

More precisely, the
off-diagonal components of the
$\cS$-matrices for the top (bottom) molecule are equal to those in
equation~(\ref{eq:S-matrix-moleculeII})  after a change of
variables $$(\cX_a,\cX_b,\cX_c) \mapsto (\cX_{a_1},\cX_{b_1},\cX_{c_1})$$
(or $(\cX_a,\cX_b,\cX_c) \mapsto  (\cX_{a_2},\cX_{b_2},\cX_{c_2})$),
with $\cX_1 =
 \cX_{a_1}, \cX_2 = \cX_{c_1},   \cX_3 = \cX_{a_2} ,
 \cX_4 = \cX_{c_2}$ and
$\cX_A =  1/(\cX_{a_1} \cX_{b_1}^2 \cX_{c_1})
 =   1/(\cX_{a_2} \cX_{b_2}^2 \cX_{c_2})$.

Monodromy representations of
the resulting $SL(2)$-representations can be expressed purely
in terms of abelian data: the abelian holonomies $\cX_{1,2,3,4}$ 
around the boundary components, the parallel transport coefficients along open
paths that run from one boundary component to another and the abelian
holonomies $\cX_A$ and $\cX_B$ around the 1-cycles $A$ and $B$. If we
divide out by diagonal gauge transformations at the  marked
points, the dependence on the parallel transport coefficients drops
out.

\subsubsection*{Length-twist coordinates}

The $SL(2)$ monodromy matrices around the loops $\alpha$ and
$\cP_{1,2,3,4}$ are all diagonal, of the form
\begin{equation}
 \rho(\alpha,\cP_{1,2,3,4}) = \begin{pmatrix} \cX_{A,1,2,3,4} &0 \\
     0 & \cX_{A,1,2,3,4}^{-1} \end{pmatrix}.
\end{equation}
The abelian holonomies
$\cX_{1,2,3,4,A}$  are equal to square-roots of the exponentiated,
complexified Fenchel-Nielsen length parameters.

In particular,  we have found that the spectral coordinate $\cX_A$
is the square-root of the exponentiated, complexified
Fenchel-Nielsen length coordinate.
Much less trivial from this construction is that the spectral coordinate
$\cX_B$ is the exponential of a complexified Fenchel-Nielsen twist
parameter.

We verify this claim by comparing the traces of the
monodromy matrices for $\rho$ with a complex version of Okai's
formula, which computes these traces in terms of the Fenchel-Nielsen
twist parameter \cite{Okai}.  (We also found it helpful to look at
Kabaya's analysis in \cite{Kabaya}, in which
a set of matrix generators of $SL(2)$ flat connections is
systemetically found in terms of a complexified version of the
Fenchel-Nielsen coordinates and compared to Okai's formula for
Fuchsian representations.) Let us consider the monodromy matrix $M_{\beta}$ along the loop $\beta$,
illustrated in Figure~\ref{fig:four-holed-sphere-cycles}. 

On the one
hand, Okai's formula computes the trace of this monodromy  as
 \begin{equation}
 \Tr M_{\beta} =  \sqrt{c}\left( \frac{t}{\cX_A} +
     \frac{\cX_A}{t} \right) + d,
 \end{equation}
in terms of the Fenchel-Nielsen twist parameter $t$. The coefficients $c$ and $d$ are given by
 \begin{align*}
 c &= \frac{1}{(\cX_A - \cX_{A}^{-1})^4} \left( \cZ_A^2 +
     \cZ_1^2 + \cZ_2^2 - \cZ_A \cZ_1 \cZ_2 - 4 \right) \left( \cZ_A^2 +
     \cZ_3^2 + \cZ_4^2 - \cZ_A \cZ_3 \cZ_4 - 4 \right)
   \\
 d &= \frac{1}{(\cX_A - \cX_{A}^{-1})^2}  \left( \cZ_A (\cZ_2 \cZ_4 + \cZ_1
   \cZ_3) - 2 (\cZ_1 \cZ_4 + \cZ_2 \cZ_3) \right),
 \end{align*}
 and $\cZ_A=\cX_A+\cX_A^{-1}$ and $\cZ_{p} = \cX_p + \cX_p^{-1}$.

On the other hand, the trace of the monodromy matrix $M_{\beta}$ of
the $SL(2)$-representation $\rho$ is given by
 \begin{align}\label{eq:FNtwist-sphere}
 \Tr M_{\beta} = &-\left( \cX_{B} +
     \frac{c}{\cX_{B}} \right) + d.
 \end{align}

 Comparing the two formulas we conclude that
 \begin{equation}\label{eq:twistspectralcompare}
 t = -\left( \frac{ \cX_A }{ \sqrt{c}} \right) \cX_B.
 \end{equation}
 In other words, the spectral coordinate $\cX_B$ is the
 exponential of a complexified version of a Fenchel-Nielsen twist
 coordinate.

 Even though the precise form of the multiplicative factor
 depends on characteristics of the chosen Fenchel-Nielsen
 network, we find that the spectral coordinate
 $\cX_{B}$  is always proportional to the Fenchel-Nielsen parameter $t$.

Of course,  in this comparison we have to be careful to match the
$SL(2)$ monodromies around the same 1-cycle $\beta$. If we compute
the $SL(2)$ monodromy along any other choice of 1-cycle $\beta'$ that
intersects the 1-cycle $\alpha$ once, equation~(\ref{eq:FNtwist-sphere}) holds in a
modified form
\begin{align}
\Tr M_{\beta'} = &-\left( \cX_{B'} +
    \frac{c'}{\cX_{B'}} \right) + d',
\end{align}
where $B'$ is a lift of $\beta'$ to the covering (following the same
rules as before). For any such choice we find that the abelian holonomy
$\cX_{B'}$ is a twist coordinate (by comparing to formulae in \cite{Okai,Kabaya}).

We can further reduce to $PSL(2)$-representations by dividing
out the transformations
\begin{align}
(\cX_1, \cX_2, \cX_a, \cX_4, \cX_5) \mapsto (\epsilon_1 \cX_1,
\epsilon_2 \cX_2, \epsilon_a \cX_a, \epsilon_4 \cX_4, \epsilon_5
\cX_5),
\end{align}
where $\epsilon_I = \pm 1$ with the constraints that $\epsilon_1 \epsilon_2
\epsilon_a = 1$ as well as $\epsilon_a \epsilon_4 \epsilon_5 = 1$
\cite{Kabaya}.\footnote{Here one has to be careful that the action of
  the group $(\Z/2\Z)^3$ on $\cX_B$ is not trivial (since $B$ is not
  an even cycle). Instead $\cX_B$
  transforms in the same way as the factor $d$. The
action of $(\Z/2\Z)^3$ \emph{is} trivial on the twist parameter
$\tilde{t} = \sqrt{\frac{(\cX_a \cX_1 - \cX_2)(1-\cX_1 \cX_2 \cX_a) (\cX_a \cX_5 -
    \cX_4)(1-\cX_a \cX_4 \cX_5)}{(\cX_a \cX_2 - \cX_1)( \cX_1
    \cX_2-\cX_a) (\cX_a \cX_4 - \cX_5)(\cX_4 \cX_5-\cX_a)}} \, t$.}

\subsubsection{Genus $2$ surface}

Finally we consider the case of a genus $2$
surface.  We choose the Fenchel-Nielsen network and the
auxiliary data as shown in Figure~\ref{fig:genustwo-aux-data}.
Our goal is to check directly that the $6$ Fenchel-Nielsen length and
twist coordinates occur as spectral coordinates.

Let us first verify that there is indeed a $6$-parameter family of equivariant
$GL(1)$-representations $\rho^{\ab}$. There are 15 simple paths $\wp$ and 3
nontrivial flatness conditions. Furthermore, there are 6 gauge
redundancies (one for each double wall). This leaves us with 6
gauge-invariant abelian holonomies. These can be associated to
holonomies along six 1-cycles $A_{1,2,3}$ and $B_{1,2,3}$ on
the cover $\Sigma$ of the genus two surface, with $A_i \# B_j =
\delta_{ij}$.\footnote{There are ten 1-cycles on
  the genus 5
  cover, but the holonomies of an equivariant $GL(1)$-connection around these
  1-cycles are not all independent: $\cX(\sigma^* \gamma) =
  \cX(\gamma)^{-1}$.}
To be explicit, these 1-cycles are lifts of 1-cycles
$\alpha_{1,2,3}$ and $\beta_{1,2,3}$ on the genus two surface, which
we choose in a way analogous
to what we did in the four-holed sphere example. That is, we identify each tube in
the pants decomposition of the genus $2$ surface with the internal
tube of a four-holed sphere, and then choose the
1-cycles $\alpha_k$ and $\beta_k$ for this tube as shown in
Figure~\ref{fig:four-holed-sphere-cycles}.

As in previous examples, we can solve for the $\cS$-matrices uniquely (by imposing the
branch point equations as well fixing the eigenvalues of the
non-abelian monodromies around the tubes), and, as can be expected
from the gluing construction, we find that they are
of the form~(\ref{eq:S-matrix-3puncturedsphere}). It is immediate that the 
complexified Fenchel-Nielsen length coordinates match the spectral coordinates
$\cX_{A_{1,2,3}}$. To find the Fenchel-Nielsen
twist coordinates we compute the $SL(2)$ monodromies around the cycles
$\beta_{1,2,3}$. We have carried out this computation, and find that the spectral
coordinates $\cX_{B_{1,2,3}}$ are indeed
 complexified Fenchel-Nielsen twist coordinates. (Similar to
 formula~(\ref{eq:twistspectralcompare}) they agree with the
 Fenchel-Nielsen twist parameter in section 12.4 of
 \cite{Kabaya} up to a multiplicative term). 

\insfigscaled{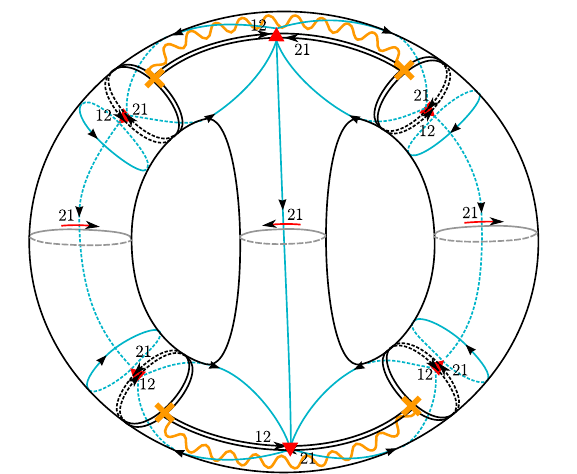}{0.9}{Example of a Fenchel-Nielsen
  network on the genus two surface, together with a choice of
  auxiliary data: branch-cut (orange wigly lines), basepoints (red
  triangles), short paths (light blue lines) and a sheet ordering for
  every orientation of the
  pants curves.}

\subsection{Mixed networks}

Just as for Fenchel-Nielsen and Fock-Goncharov type networks, we can
apply the nonabelianization map to any mixed type spectral network. In
particular, for the examples discussed in
\S\ref{sec:spectralnetworks}, we find unique solutions for the
$\cS$-matrices. The off-diagonal elements in these solutions have a
similar interpretation in the limiting network  as abelian holonomies
along auxiliary paths on $\Sigma$.

Using this explicit form of the $\cS$-matrices, we obtain
$SL(2)$-representations in terms of purely abelian data. We can
compute a 
monodromy representation in terms of the spectral coordinates, and
compare to known coordinate systems. We have done this, and found
that not all spectral
coordinates for mixed spectral networks can be expressed as
Fenchel-Nielsen or Fock-Goncharov coordinates. We leave
the further exploration of these coordinates for future work.

\section{Some consequences}\label{sec:physcons}

\subsection{Asymptotics}\label{sec:asymptotics}

One of the most interesting properties of the spectral coordinates \cite{Gaiotto2012}
is that they have simple asymptotics when evaluated along certain
natural one-parameter families of $SL(2)$-connections.  This asymptotic
property was described in \cite{Gaiotto:2009hg} in the case of Fock-Goncharov spectral networks.
One consequence of the constructions described in this paper is that a similar asymptotic property
should hold for Fenchel-Nielsen coordinates.
In this section we briefly describe that story.

We will consider
one-parameter families of flat $SL(2)$-connections, $\{\nabla(\zeta)\}_{\zeta \in \C^\times}$.
We will not consider arbitrary families, but rather only the ones which come from
solutions of Hitchin equations. 
These equations concern a tuple $(E,D,\varphi)$ where $E$ is
a Hermitian rank $2$ bundle over $C$, $D$ a (generally non-flat) 
unitary connection in $E$, and $\varphi$
a section of $\Omega^{1,0}(\End E)$. The equations are:
\begin{align}
F_D + [\varphi,\varphi^\dagger] &= 0, \\
\bar\partial_D \varphi &= 0.
\end{align}
The moduli of solutions of these equations was studied in \cite{MR89a:32021},
and in \cite{hbnc} this analysis was extended to the case
where $D$ and $\varphi$ are allowed to have first-order poles.
This is the case relevant for us.

Given a solution of Hitchin's equations, 
the corresponding family of flat (non-unitary) connections is
\begin{equation} \label{eq:wkb-family}
 \nabla(\zeta) = \zeta^{-1} \varphi + D + \zeta \varphi^\dagger.
\end{equation}
Indeed Hitchin's equations say precisely that
the connection $\nabla(\zeta)$ is flat for every $\zeta \in \C^\times$.

Now we ask:
how does $\nabla(\zeta)$ behave asymptotically as $\zeta \to 0$?
The answer to this question
is somewhat intricate; in particular, it turns out to depend on precisely
how $\zeta$ approaches $0$.  One can get simpler answers if one restricts
$\zeta$ to lie in a half-plane, say the open half-plane $\HH_\vartheta$ centered
on the ray $e^{\I \vartheta} \R_+$, and if one chooses
the correct coordinate system on the moduli of flat connections.  Namely,
given the Higgs field $\varphi$, we consider the quadratic differential
\begin{equation}
\varphi_2 = \Tr (e^{-2 \I \vartheta} \varphi^2).
\end{equation}
Because $\varphi$ has only first-order poles, $\varphi_2$ has at
most second-order poles. Moreover, if $\varphi$ is chosen generically,
then $\varphi_2$ has only simple zeroes; we assume that from now on.
As described in \S\ref{sec:WKB}, the quadratic differential $\varphi_2$ corresponds to
a particular spectral curve $\Sigma \subset T^* C$ and spectral network $\cW(\varphi_2)$.
This spectral network
induces spectral coordinate functions $\cX_\gamma$.  The $\zeta \to 0$ asymptotics of the
spectral coordinates of the family $\nabla(\zeta)$
were studied in \cite{Gaiotto:2009hg,Gaiotto2012}, where it was argued that
they are controlled by the periods of the spectral curve:
\begin{equation} \label{eq:asymp}
\cX_\gamma(\nabla(\zeta)) \sim c_\gamma \exp \left( \zeta^{-1} \oint_\gamma \lambda \right)
\end{equation}
where $c_\gamma$ is $\zeta$-independent, and $\lambda$ denotes the tautological (Liouville)
1-form on $T^* C$.

In particular, suppose $\varphi$ is a generic Higgs field, so that $\varphi_2$ is a
generic quadratic differential.  In this case the spectral network $\cW(\varphi_2)$
is a Fock-Goncharov network, and the corresponding coordinates are Fock-Goncharov coordinates.
Thus \eqref{eq:asymp} tells us the asymptotic behavior of the Fock-Goncharov coordinates
of the connections $\nabla(\zeta)$. These asymptotics were derived using the
WKB method in \cite{Gaiotto:2009hg}, and a more general method of obtaining them
was described in \cite{Gaiotto2012}.

On the other hand, suppose that $\varphi$ is a special Higgs field, such that $\varphi_2$
is a Strebel differential.  In this case $\cW(\varphi_2)$ is a Fenchel-Nielsen network,
and \eqref{eq:asymp} tells us the asymptotic behavior of the Fenchel-Nielsen coordinates
of $\nabla(\zeta)$.  (More precisely, the spectral coordinates associated to $\cW(\varphi_2)$
depend on whether we take the British or American resolution; the asymptotics \eqref{eq:asymp}
in the open half-plane $\HH_\vartheta$
apply to either.  We believe that
for the British resolution the asymptotics extend also to one of the rays
on the \ti{boundary}
of the half-plane, and for the American resolution they extend to the other boundary.
It would be useful to verify this by direct application of the WKB method.)

Remember that the periods
$\oint_{\gamma} \lambda$ are real-valued when $\varphi_2$ is a
Strebel differential and $\gamma$ is the lift to $\Sigma$ of a pants
curve on $C$. The Fenchel-Nielsen length coordinates thus behave asymptotically as the
exponentials of a set of real period integrals. This makes precise the
identification between Coulomb parameters and length
coordinates that follows from the AGT correspondence.

\subsection{Line defects} \label{sec:defects}

A second interesting consequence of our discussion concerns the physics of line defects.

Recall from \cite{Drukker:2009tz,Gaiotto:2010be} that, in the theory $S[A_1, C]$ associated to a punctured
Riemann surface $C$, there is a class of supersymmetric line defects $L(\wp,\zeta)$ labeled by closed
loops $\wp$ on $C$ and phases $\zeta \in \C^\times$.

According to \cite{Gaiotto:2010be}, the spectrum of framed BPS states
attached to $L(\wp,\zeta)$ can be read off from the abelianization map, as follows.  Let $M$ denote
the holonomy of an $SL(2)$-connection around $\wp$.  Given any spectral network $\cW$,
the $\cW$-abelianization map expresses $\Tr M$
as a sum of $GL(1)$ holonomies
\begin{equation}
\Tr M = \sum_\gamma \fro(\wp,\gamma,\cW) \cX_\gamma,
\end{equation}
where the coefficients $\fro(\wp,\gamma, \cW) \in \Z$.  If $\cW$ is the WKB spectral network
corresponding to a quadratic differential $\varphi_2 / \zeta^2$,
then these coefficients have a physical interpretation:  $\fro(\wp, \gamma, \cW)$ counts the framed BPS
states of charge $\gamma$ attached to the line defect $L(\wp,\zeta)$,
at the point of the Coulomb branch determined by $\varphi_2$.

If $\varphi_2$ is generic, and $\wp$ is not a cycle contractible to a puncture,
this always leads to at least three framed BPS states.
However, if $\cW$ happens to be a Fenchel-Nielsen network, or more generally a mixed
network whose complement contains at least one annulus $A$, then we may take $\wp$ to be a path
going around $A$.  In this case, $\cW$-abelianization for $\wp$ is very simple:
$\wp$ crosses no walls at all, so
\begin{equation}
\Tr M = \cX_\gamma + \cX_{-\gamma}
\end{equation}
where $\gamma$ and $-\gamma$ are the two lifts of $\wp$ to the double cover $\Sigma$.
Thus, in this situation the supersymmetric line defect $L(\gamma, \zeta)$ supports just two
framed BPS states, carrying charges $\gamma$, $-\gamma$.

Physically we would interpret this in the following way.  Let us choose a
particular way of looking at the theory $S[A_1, C]$, in which this theory is
obtained by gauging a particular $SU(2)$ symmetry in another theory $S[A_1,C']$,
with $C'$ obtained from $C$ by cutting along the annulus $A$.
From this point of view, the line defect $L(\wp,\zeta)$ is a Wilson line for the
new $SU(2)$ gauge symmetry, in the fundamental representation.
Naive classical reasoning would suggest that in the IR this
defect should support two framed BPS states, corresponding to the
decomposition
of the fundamental representation into two weight spaces under $U(1) \subset SU(2)$.
At generic $(\varphi_2, \zeta)$ this classical reasoning is not exactly correct:  we do get these
two states but we get
additional states as well.  What we have found here is that, if $(\varphi_2, \zeta)$ are
chosen specially, the classical picture is precisely correct.  (See also some related discussion
in \cite{Cordova2013}).

It would be very interesting to have a direct physical understanding of \ti{why}
this simplification occurs.

\bibliographystyle{utphys}
\bibliography{fn-paper}

\end{document}